\documentclass{amsart}

\usepackage{amsfonts,amsthm,amsmath,amssymb,latexsym}
\usepackage{amscd}

\theoremstyle{plain}

\theoremstyle{definition}

\theoremstyle{remark}

\begin{document}

\begin{flushleft}
{\Large\sf 7 June, 2009}
\end{flushleft}

\bigskip
\bigskip

\begin{center}
   {\LARGE\bf  Brou{\'e}'s abelian defect group conjecture holds\\
           for the Harada-Norton sporadic simple group {\bf\sf HN}}
\end{center}
\bigskip
\begin{center}
   {\large
        {\bf Shigeo Koshitani}$^{\text a,*}$, 
        {\bf J{\"u}rgen M{\"u}ller}$^{\text b}$ 
   }
\end{center}
\bigskip
\begin{center}
  {\it\small
${^{\mathrm{a}}}$Department of Mathematics, Graduate School of Science, \\
Chiba University, Chiba, 263-8522, Japan \\
${^{\mathrm{b}}}$Lehrstuhl D f{\"u}r Mathematik,
RWTH Aachen University, 52062, Aachen, Germany}
\end{center}

\bigskip
\bigskip
\bigskip
\bigskip
\bigskip

\footnote
{$^*$ Corresponding author. \\
\indent {\it E-mail addresses:} koshitan@math.s.chiba-u.ac.jp (S.Koshitani),
\\
juergen.mueller@math.rwth-aachen.de (J.M{\"u}ller).
}

%\begin{center}
%Dedicated to Professor Hiroyuki Tachikawa on his eightieth birthday
%\end{center}
%
\begin{center}
Dedicated to Professor Tetsuro Okuyama on his sixtieth birthday
\end{center}

\bigskip
\bigskip
\bigskip
\bigskip

\hrule 
{\small
\noindent{\bf Abstract}

\smallskip
In representation theory of finite groups, there is a well-known and
important conjecture due to M.~Brou{\'e}. He conjectures that, for
any prime $p$, if a $p$-block $A$ of a finite group $G$ has an abelian
defect group $P$, then $A$ and its Brauer corresponding block $B$
of the normaliser $N_G(P)$ of $P$ in $G$ are derived equivalent
(Rickard equivalent).
This conjecture is called 
{\it Brou{\'e}'s abelian defect group conjecture}.
We prove in this paper that Brou{\'e}'s abelian defect group conjecture
is true for a non-principal $3$-block $A$ with an elementary abelian
defect group $P$ of order $9$ of the Harada-Norton
simple group {\sf HN}.
It then turns out that Brou{\'e}'s abelian defect group conjecture
holds for all primes $p$ and for all $p$-blocks 
of the Harada-Norton simple group {\sf HN}.

\bigskip\noindent
{\it Keywords:} Brou{\'e}'s conjecture; abelian defect group;
Harada-Norton simple group}
\bigskip
\hrule

\bigskip\bigskip

\begin{flushleft}
{\bf 1. Introduction and notation}
\end{flushleft}

\bigskip\noindent
In representation theory of finite groups, one of the
most important and interesting problems
is to give an affirmative answer to a conjecture,
which was introduced by M.~Brou{\'e} around 1988
\cite{Broue1990},
and is nowadays called {\it Brou{\'e}'s Abelian
Defect Group Conjecture}.
He actually conjectures the following:

\bigskip\noindent
{\bf 1.1.Conjecture} 
(Brou\'e's Abelian Defect Group Conjecture)
(\cite[6.2.Question]{Broue1990} and 
\cite[Conjecture in p.132]{KoenigZimmermann}){\bf.}
{\it
Let $p$ be a prime, and let $(\mathcal K, \mathcal O, k)$ be a
splitting $p$-modular system for all subgroups of a
finite group $G$. Assume that $A$ is a block algebra of
$\mathcal OG$ with a defect group $P$ and that $B$ is a
block algebra of $\mathcal ON_G(P)$ such that $B$ is the
Brauer correspondent of $A$, where $N_G(P)$ is the normaliser of
$P$ in $G$. Then, $A$ and $B$ should be derived equivalent
(Rickard equivalent)
provided $P$ is abelian.
}

\bigskip\noindent
In fact, a stronger conclusion than {\bf 1.1} is expected.
If $G$ and $H$ are finite groups and if $A$ and $B$ are block algebras of 
$\mathcal OG$ and $\mathcal OH$ (or $kG$ and $kH$) respectively, we say that
$A$ and $B$ are 
{\it splendidly Rickard equivalent} in the sense
of Linckelmann (\cite{Linckelmann1998}, \cite{Linckelmann2001}), where
he calls it a {\it splendid derived equivalence},
see the end of {\bf 1.8}.
Note that this is the same as that given by Rickard in
\cite{Rickard1996} when $A$ and $B$ are the principal block algebras,
which he calls a {\it splendid equivalence}.

\bigskip
\noindent
{\bf 1.2.Conjecture}
(Rickard \cite{Rickard1996}, 
\cite[Conjecture 4, in p.193]{Rickard1998}){\bf .}
{\it
Keep the notation, and suppose that $P$ is abelian as in {\bf 1.1}.
Then, there should be a splendid Rickard equivalence between
the block algebras $A$ of $\mathcal OG$ and $B$ of $\mathcal ON_G(P)$.
}

\bigskip\noindent
There are several cases where the conjectures of
Brou{\'e} {\bf 1.1} and Rickard {\bf 1.2} are checked.
For example we prove that {\bf 1.1} and {\bf 1.2} are true for 
the principal block algebra $A$ of an arbitrary finite group $G$
when the defect group $P$ of $A$ is elementary abelian
of order $9$ (and hence $p = 3$), see
\cite[(0.2)Theorem]{KoshitaniKunugi2002}.
Then, it may be natural to ask what about the case of
non-principal block algebras with the same defect group
$P = C_3 \times C_3$. Namely, this paper should be 
considered as a continuation of such a project, which has
already been accomplished for several cases in our previous papers 
for the O'Nan simple group and the Higman-Sims simple
group in \cite[0.2.Theorem]{KoshitaniKunugiWaki2002},
for the Held simple
group and the sporadic simple Suzuki group in 
\cite[Theorem]{KoshitaniKunugiWaki2004},
and for the Janko's simple group $J_4$
\cite[Theorem~1.3]{KoshitaniKunugiWaki2008},
see also 
\cite{MuellerSchaps} and \cite{Kunugi}.
That is to say, our main theorem of this paper is the following:

\bigskip\noindent
{\bf 1.3.Theorem.} {\it
Let $G$ be the Harada-Norton simple group ${\sf HN}$, and 
let $(\mathcal K, \mathcal O, k)$ be a splitting $3$-modular
system for all subgroups of $G$, 
see the definition {\bf 1.8} below. 
Suppose that $A$ is a non-principal block algebra of $\mathcal OG$
with a defect group $P$ which is an elementary abelian group
$C_3 \times C_3$ of order $9$, and that $B$ is a block algebra
of $\mathcal ON_G(P)$ such that $B$ is the Brauer correspondent
of $A$. Then, $A$ and $B$ are splendidly Rickard
equivalent, and hence the conjectures 
{\bf 1.1} and {\bf 1.2} of
Brou{\'e} and Rickard hold.
}

\bigskip\noindent
As a matter of fact, the main result {\bf 1.3} above is obtained
by proving the following:

\bigskip\noindent
{\bf 1.4.Theorem.} {\it
Keep the notation and the assumption as in {\bf 1.3}.
Then, the non-principal block algebra $A$ of $\mathcal OG$
with a defect group $P = C_3 \times C_3$ and the principal
block algebra $A'$ of $\mathcal O \mathcal{\sf HS}$ 
of the Higman-Sims simple group are 
Puig equivalent, that is 
$A$ and $A'$ are Morita equivalent
which is realized by a $\Delta P$-projective 
$p$-permutation $\mathcal O[G \times {\sf HS}]$-module,
in other words, $A$ and $A'$ have isomorphic source
algebras as interior $P$-algebras.
}

\bigskip\noindent
Then, it turns out that, as a corollary to the main result (1.3),
we eventually can prove that

\bigskip\noindent
{\bf 1.5.Corollary.} {\it 
Brou\'e's abelian defect group conjecture {\bf 1.1} and 
even Rickard's splendid equivalence conjecture {\bf 1.2} are true
for all primes $p$ and for all block algebras of $\mathcal OG$
when $G = {\sf HN}$.
}

\bigskip\noindent
{\bf 1.6.Starting point and strategy.}
A story of the birth of this paper is actually very similar
to that of the Janko's simple group $J_4$
which is given in \cite[{\bf 1.6}]{KoshitaniKunugiWaki2008}.
Namely, relatively recently
%JM: Add a citation
G.~Hiss, J.~M{\"u}ller, F.~Noeske and J.G.~Thackray 
\cite{HissMuellerNoeskeThackray}
have determined the $3$-decomposition matrix
of the group ${\sf HN}$ with defect group $C_3 \times C_3$,
see {\bf 4.1}.
Our starting point for this work 
was actually to realize that
the $3$-decomposition matrix for the non-principal
block of ${\sf HN}$ 
with an elementary abelian defect group of
order $9$ is exactly the same as that
for the principal $3$-block of the
Higman-Sims simple group $\sf HS$.
Furthermore, the generalised 3-decomposition 
matrices of these two blocks are the same.
Therefore, it is natural to suspect 
whether these two $3$-block algebras 
would be Morita equivalent not only over
an algebraically closed field $k$ of 
characteristic $3$ but also over
a complete discrete valuation ring
$\mathcal O$ whose residue field is $k$,
and we might expect even that they are {\it Puig equivalent}
(we shall give a precise definition of Puig equivalence 
in {\bf 1.8} below).
Anyhow, since the two conjectures of Brou\'e and Rickard
in {\bf 1.1} and {\bf 1.2} respectively
have been solved for the principal
$3$-block of ${\sf HS}$ in a 
paper of Okuyama \cite{Okuyama1997} 
it turns out
that Brou{\'e}'s abelian defect group conjecture {\bf 1.1}
and Rickard's splendid equivalence conjecture {\bf 1.2}
shall be solved also for the non-principal 
$3$-block of ${\sf HN}$ with the same defect
group $C_3 \times C_3$.

\bigskip\noindent
{\bf 1.7.Contents.}
In \S2, we shall give several fundamental lemmas, which are
useful and powerful to prove our main results.
In \S\S 3 and 4, we shall investigate $3$-modular representations
for ${\sf HN}$ and we shall get trivial source ($p$-permutation)
modules which are in the non-principal $3$-block $A$
of ${\sf HN}$ with a defect group $P = C_3 \times C_3$.
In \S 5, we shall list data on Green correspondents
of simples in the principal $3$-block
$A'$ of ${\sf HS}$,
which are known by a result of 
\cite[Theorem]{Waki1993}, see
\cite[Example 4.8]{Okuyama1997}.
Finally, in \S\S 6-8, we shall give complete proofs of
our main results {\bf 1.3}, {\bf 1.4} and {\bf 1.5}.

%JM: more comments on computations
%It should be noted that calculations of computers,
%say {\sf GAP} are used for having character tables,
%induced characters and also structure of subgroups
%of finite groups.

To achieve our results, next to theoretical reasoning
we have to rely on fairly heavy computations.
As tools, we use the computer algebra system {\sf GAP} \cite{GAP},
to calculate with permutation groups as well as with ordinary 
and Brauer characters. We also make use of the data library 
\cite{CTblLib}, in particular allowing for easy access to the data
compiled in \cite{Atlas}, \cite{ModularAtlas} and \cite{ModularAtlasProject},
and of the interface \cite{AtlasRep} to the data library \cite{ModAtlasRep}.
Moreover, we use the computer algebra system {\sf MeatAxe} \cite{MA}
to handle matrix representations over finite fields,
as well as its extensions to compute 
submodule lattices \cite{LuxMueRin},
radical and socle series \cite{LuxWie},
homomorphism spaces and endomorphism rings \cite{LuxSzoke},
and direct sum decompositions \cite{LuxSzokeII}.
We give more detailed comments on the relevant computations
in the spots where they enter the picture.

\bigskip\noindent
{\bf 1.8.Notation.}
Throughout this paper, we use the following notation and
terminology. Let $A$ be a ring. We denote by $1_A$, $Z(A)$ and $A^\times$
for the unit element of $A$, the centre of $A$ 
and the set of all units in $A$, respectively.
We denote by $\mathrm{rad}(A)$ the Jacobson radical of $A$ and 
by $\mathrm{rad}^i(A)$ the $i$-th power $(\mathrm{rad}(A))^i$
for any positive integer $i$ while we define
$\mathrm{rad}^0(A) = A$.
We write $\mathrm{Mat}_n(A)$ 
for the matrix ring of all $n \times n$-matrices 
whose entries are in $A$.
Let $B$ be another ring.
We denote by $\mathrm{mod}{\text -}A$,
$A{\text -}\mathrm{mod}$ and
$A{\text -}\mathrm{mod}{\text -}B$
the categories of finitely generated right $A$-modules,
left $A$-modules and $(A,B)$-bimodules, respectively.
We write $M_A$, $_AM$ and $_AM_B$ when $M$ is 
a right $A$-module,
a left $A$-module and an $(A,B)$-bimodule. 
However, by a module we mean 
a finitely generated right module unless otherwise stated.
Let $M$ and $N$ be $A$-modules. We write $N | M$ if $N$
is (isomorphic to) a direct summand of $M$ as an $A$-module.

%JM: `From' should not be at the beginning of a line.
\mbox{}From now on, let $k$ be a field and assume that $A$ is a
finite dimensional $k$-algebra. Suppose that $M$ is an $A$-module.
Then, we denote by $\mathrm{soc}(M)$ the socle of $M$.
We define $\mathrm{soc}_0(M) = 0$ and 
$\mathrm{soc}_1(M) = \mathrm{soc}(M)$.
Then, we define $\mathrm{soc}_i(M)$ by
$\mathrm{soc}_i(M)/\mathrm{soc}_{i-1}(M)
 = \mathrm{soc}(M/\mathrm{soc}_{i-1}(M))$ for any
integer $i \geqslant 2$. 
Similarly, we write $\mathrm{rad}^i(M)$ for
$M{\cdot}\mathrm{rad}^i(A)$ for any integer $i \geqslant 0$.
By using this, we define $L_i(M)$ by
${\mathrm{rad}}^{i-1}(M)/{\mathrm{rad}}^{i}(M)$
for $i = 1, 2, \cdots$. We call $L_i(M)$ the $i$-th Loewy
layer of $M$. We denote by $j(M)$ the Loewy length of $M$,
namely $j(M)$ is the least positive integer $j$ satisfying
${\mathrm{rad}}^j(M) = 0$. 
We write $P(M)$ and $I(M)$ for the projective cover and the
injective hull (envelope) of $M$, respectively, and
we write $\Omega$ for the Heller operator (functor), namely,
$\Omega M$ is the kernel of the projective cover
$P(M) \twoheadrightarrow M$.
Dually, $\Omega^{-1}M$ is the cokernel of the
injective hull 
$M \rightarrowtail I(M)$.
For simple $A$-modules $S_1, \cdots, S_n$ (some of which are
possibly isomorphic) 
%JM: Replaced this notation by a `boxed' business.
%we denote by $U(S_1, \cdots, S_n)$ 
%a uniserial $A$-module $M$ such that
%$\mathrm{rad}^{i-1}(M)/\mathrm{rad}^i(M) \cong S_i$ for
%$i = 1, \cdots, n$. 
%Note that $A$-modules of such a form
%$U(S_1, \cdots, S_n)$ possibly are not uniquely determined
%up to isomorphism in general.
%For the same $S_1, \cdots, S_n$,
we write that $M = a_1 \times S_1 + \cdots + a_n \times S_n$,
as composition factors for positive integers $a_1, \cdots, a_n$
when the set of all composition factors are
$a_1$ times $S_1$, $\cdots$, $a_n$ times $S_n$.
For an $A$-module $M$ and a simple $A$-module $S$,
we denote by $c_M(S)$ the multiplicity of
all composition factors of $M$ which are isomorphic to $S$.
We write $c(S,T)$ for $c_{P(S)}(T)$ for simple
$A$-modules $S$ and $T$, namely, this is 
so-called the Cartan invariant with respect to $S$ and $T$. 

%JM: Comments on pictures
%JM: now a bit more cautious with comments on diagrams
To describe the structure of an $A$-module, we either indicate
the radical and socle series, in cases where these series coincide
and are sufficient for our analysis,
or we draw an Alperin diagram \cite{Alperin}. 
An $A$-module need not have an Alperin diagram, but if it does then
it is a compact way to give a more detailed 
%The latter are only defined
%if all Loewy layers contain any composition factor with multiplicity
%at most one, but in this case are a compact way to give a complete 
structural description of the module under consideration;
note that the Alperin diagram is closely related to the Hasse diagram
of the incidence relation amongst the local submodules
in the sense of \cite{MueLatt}, hence for explicit examples 
is easily determined using
the techniques described in \cite{LuxMueRin}.
Note, however, that by giving any kind of diagram an $A$-module
in general is not uniquely determined up to isomorphism.

Let $N$ be another $A$-module. Then,
$\mathrm{Hom}_A(M,N)$ is the set of all 
right $A$-module-homomorphisms from $M$ to $N$, 
which canonically is a $k$-vector space, and we denote by
$\mathrm{PHom}_A(M,N)$ the set of all (relatively) projective
homomorphisms in $\mathrm{Hom}_A(M,N)$, which is a $k$-subspace
of $\mathrm{Hom}_A(M,N)$. Hence, we can define the factor space,
that is, we write $\underline{\mathrm{Hom}}_A(M,N)$ 
for the factor space 
$\mathrm{Hom}_A(M,N)/\mathrm{PHom}_A(M,N)$.
By making use of this, as well-known, we can construct
the stable module category $\underline{\mathrm{mod}}{\text -}A$,
which is a quotient category of $\mathrm{mod}{\text -}A$
such that the set of all morphisms is given by
$\underline{\mathrm{Hom}}_A(M,N)$.

In this paper, $G$ is always a finite group and we fix a
prime number $p$. Assume that $(\mathcal K, \mathcal O, k)$ is a
splitting $p$-modular system for all subgroups of $G$, that is
to say, $\mathcal O$ is a complete discrete valuation ring of
rank one such that its quotient field is $\mathcal K$ which is
of characteristic zero and its residue field
$\mathcal O/\mathrm{rad}(\mathcal O)$ is $k$ which is of
characteristic $p$, and that $\mathcal K$ and $k$ are
splitting fields for all subgroups of $G$.
We mean by an $\mathcal OG$-lattice a finitely generated
right $\mathcal OG$-module which is a free $\mathcal O$-module.
We sometimes call it just an $\mathcal OG$-module.
Let $X$ be a $kG$-module. Then, we write $X^\vee$ for the
$k$-dual of $X$, namely, $X^\vee = \mathrm{Hom}_k(X,k)$ which
is again a right $kG$-module via
$(x)(\varphi g) = (xg^{-1})\varphi$ for $x \in X$,
$\varphi \in X^\vee$ and $g \in G$.
Similarly, we write $\chi^{\vee}$ for the dual 
(complex conjugate) of $\chi$ for an ordinary
character $\chi$ of $G$.
Let $H$ be a subgroup of $G$, and let $M$ and $N$ be
an $\mathcal OG$-lattice and an $\mathcal OH$-lattice, respectively.
Then, let ${M}{\downarrow}^G_H = {M}{\downarrow}_H$ be the
restriction of $M$ to $H$, and let 
${N}{\uparrow}_H^G = {N}{\uparrow}^G$ be the induction
(induced module) of $N$ to $G$, that is,
${N}{\uparrow}^G = (N \otimes_{\mathcal OH}\mathcal OG)_{\mathcal OG}$.
Similar for $kG$- and $kH$-modules. 

We denote by $\mathrm{Irr}(G)$ and $\mathrm{IBr}(G)$ the sets of
all irreducible ordinary and Brauer characters of $G$,
respectively.
Let $A$ be a block algebra ($p$-block) of $\mathcal OG$.
Then, we write $\mathrm{Irr}(A)$ and $\mathrm{IBr}(A)$ for
the sets of all characters in $\mathrm{Irr}(G)$ and $\mathrm{IBr}(G)$
which belong to $A$, respectively.
We often mean by $\mathrm{IBr}(A)$ the set
of all non-isomorphic simple $kG$-modules belonging to $A$.
We sometimes denote by $A^*$ the block algebra of $kG$
corresponding to $A$.
But, we usually abuse $A$ and $A^*$, namely, we often
mean the block algebra of $kG$ by $A$ as well
when it is clear from the context.
For ordinary characters $\chi$ and $\psi$ of $G$, we denote by
$(\chi, \psi)^G$ the inner product of $\chi$ and $\psi$ in
usual sense. 
Let $X$ and $Y$ be $kG$-modules.
Then, we write $[X, Y]^G$ for 
$\mathrm{dim}_k[\mathrm{Hom}_{kG}(X,Y)]$.
We denote by $k_G$ the trivial $kG$-module.
Similar for $\mathcal O_G$. 
For $A$-modules $M$ and $N$ we write $[M,N]^A$ for 
$\mathrm{dim}_k[\mathrm{Hom}_{A}(M,N)]$.

We say that $M$ is
a {\it trivial source} ($p$-{\it permutation}) $kG$-module 
if $M$ is an indecomposable $kG$-module whose source is $k_Q$,
where $Q$ is a vertex of $M$.
Let $G'$ be another finite group, and let $V$ be an 
$(\mathcal OG, \mathcal OG')$-bimodule. 
Then we can regard $V$ as a right
$\mathcal O[G \times G']$-module via
$v(g, g') = g^{-1}vg'$ for $v \in V$, $g \in G$ and $g' \in G'$.
Similar for $(kG, kG')$-bimodules.
We denote by $\Delta G$ the diagonal copy of $G$ in $G \times G$,
namely, $\Delta G = \{ (g,g) \in G \times G \, | \, g \in G \}$.
Let $A$ and $A'$ be block algebras of $\mathcal OG$ and $\mathcal OG'$, 
respectively. Then, we say that $A$ and $A'$ are {\it Puig equivalent}
if $A$ and $A'$ have a common defect group $P$ 
(and hence $P \subseteq G \cap G'$)
and if there is a Morita equivalence between $A$ and $A'$
which is induced by an $(A,A')$-bimodule $\mathfrak M$ 
such that, as a right
$\mathcal O[G \times G']$-module, 
$\mathfrak M$ is a $p$-permutation (trivial source)
module and $\Delta P$-projective. Similar for blocks of $kG$ and $kG'$.
Due to a result of Puig (and independently of Scott), 
see \cite[Remark 7.5]{Puig1999},
this is equivalent to a condition that
$A$ and $A'$ have source algebras which are isomorphic as
interior $P$-algebras,
see \cite[Theorem 4.1]{Linckelmann2001}.
For an $(\mathcal OG, \mathcal OG')$-bimodule $V$ and a common
subgroup $Q$ of $G$ and $G'$, we set
$V^{Q} 
%= \{ v \in V \ | \ v(q,q) = v, \forall q \in Q \}
= \{ v \in V \ | \ qv = vq, \forall q \in Q \}$.
If $Q$ is a $p$-group, the Brauer construction is defined to be
a quotient
$V(Q) = V^{Q}/ [\sum_{R \lneqq Q}
 {\mathrm{Tr}}{\uparrow}_R^Q (V^R) 
 + {\mathrm{rad}}{\mathcal O}{\cdot}V^{Q}]$
where ${\mathrm{Tr}}{\uparrow}_R^Q$ is the usual trace map.
The Brauer homomorphism
${\mathrm{Br}}_{Q}: (\mathcal OG)^{Q} 
 \rightarrow  kC_G(Q)$
is obtained from composing the canonical epimorphism
$(\mathcal OG)^{Q} \twoheadrightarrow
 (\mathcal OG)(Q)$
and a canonical isomorphism
$(\mathcal OG)(Q) 
 \overset{\approx}{\rightarrow} kC_G(Q)$.

We say that $A$ and $A'$ are {\it stably equivalent of 
Morita type} if there exists an $(A, A')$-bimodule 
$\mathfrak M$ such that
$_A(\mathfrak M \otimes_{A'} \mathfrak M^\vee)_A 
  \cong {_A}{A}{_A} \oplus
 ({\mathrm{projective}} \ (A,A){\text{-}}{\mathrm{bimodule}})$
and
$_{A'} (\mathfrak M^\vee \otimes_A \mathfrak M)_{A'} 
\cong {_{A'}}{A'}{_{A'}} \oplus
 ({\mathrm{projective}} \ (A' ,A'){\text{-}}{\mathrm{bimodule}})$.
We say that $A$ and $A'$ are {\it splendidly stably equivalent of 
Morita type}
if $A$ and $A'$ have a common defect group $P$ and 
the stable equivalence of Morita type is induced by
an $(A,A')$-bimodule $\mathfrak M$ 
which is a $p$-permutation (trivial source)
$\mathcal O[G \times G']$-module and is $\Delta P$-projective,
see \cite[Theorem 3.1]{Linckelmann2001}.
We say that $A$ and $A'$ are {\it Rickard equivalent}
if $A$ and $A'$ are derived equivalent, namely,
${\mathrm{D}}^b({\mathrm{mod}}{\text{-}}A)$ and 
${\mathrm{D}}^b({\mathrm{mod}}{\text{-}}A')$ are
equivalent as triangulated categories.
We say that $A$ and $A'$ are 
{\sl splendidly Rickard equivalent}
if $A$ and $A'$ are derived equivalent
by a complex 
$M^{\bullet} \in
{\mathrm{C}}^b (A{\text{-}}{\mathrm{mod}}{\text{-}}A')$
and its dual
$(M^{\bullet})^{\vee}$ such that
each term $M^n$ of $M^{\bullet}$ is
a $\Delta(P)$-projective and 
$p$-permutation module as an
$\mathcal O[G \times G']$-module,
where 
${\mathrm{C}}^b (A{\text{-}}{\mathrm{mod}}{\text{-}}A')$
is the category of bounded complexes
of finitely generated $(A,A')$-bimodules.

For a positive integer $n$,
%JM: symmetric group
$\mathcal A_n$ and $\mathcal S_n$ denote the alternating and
symmetric group
on $n$ letters, $M_n$ denotes the Mathieu group,
and $C_n$, $D_n$ and $SD_n$ denote 
the cyclic group, the dihedral group
and the semi-dihedral group
of order $n$, respectively.
For a subgroup $E$ of $\mathrm{Aut}(G)$, $G \rtimes E$ denotes a
semi-direct product such that $G$ is normal in
$G \rtimes E$ and $E$ acts on $G$ canonically.
For $g \in G$ and a subset $S$ of $G$, 
we denote $g^{-1}Sg$ by $S^g$, and similarly,
$x^g = g^{-1}xg$ for $x \in G$.
For non-empty subsets $S$ and $T$ of $G$,
we write $S =_G T$ if $T = S^g$ for an
element $g \in G$.

For other notation and terminology, see the books of
Nagao-Tsushima \cite{NagaoTsushima} and Th\'evenaz \cite{Thevenaz}.

\bigskip

%\newpage

\begin{flushleft}
{\bf 2. Preliminaries}
\end{flushleft}

\bigskip\noindent
In this section we list many lemmas, some of which are
theorems due to other people.
These lemmas are so useful and powerful to prove
our main results.

\bigskip\noindent
{\bf 2.1.Lemma} (\cite[(1.1)Lemma]{Koshitani1985}){\bf.}
{\it 
Let $A$ be a finite-dimensional algebra over a field and
$X$ an $A$-module. Assume that $Y$ is a non-zero uniserial 
$A$-submodule of $X$ with Loewy layers
$$
  {\mathrm{rad}}^{i-1}(Y) /{\mathrm{rad}}^{i}(Y) \cong S_i
\qquad {\text{for}} \ i = 1, \cdots , n
$$
where $S_i$ is a simple $A$-module.
Set $\bar{X} = X/Y$. Then, we get the following:
\begin{enumerate}
  \renewcommand{\labelenumi}{(\roman{enumi})}
    \item
For each $ j = 1, \cdots , j(X)$, \  
$
{\mathrm{rad}}^{j-1}(X) /{\mathrm{rad}}^{j}(X) \cong
{\mathrm{rad}}^{j-1}(\bar X) /{\mathrm{rad}}^{j}(\bar X)
\quad {\text{or}}
\\
{\mathrm{rad}}^{j-1}(X) /{\mathrm{rad}}^{j}(X) \cong
{\mathrm{rad}}^{j-1}(\bar X) /{\mathrm{rad}}^{j}(\bar X) 
\bigoplus S_i
$
for some $S_i$.
\item
For each $i = i, \cdots , n$, there is a positive integer $m_i$
such that $m_1 < m_2 < \cdots < m_n$ and
that
${\mathrm{rad}}^{{m_i}-1}(X) /{\mathrm{rad}}^{m_i}(X) \cong
 \Big({\mathrm{rad}}^{{m_i}-1}(\bar X) 
      /{\mathrm{rad}}^{m_i}(\bar X) \Big)
\bigoplus S_i$.
\end{enumerate}
}

\bigskip\noindent
{\bf 2.2.Lemma} (Okuyama \cite[Lemma 2.2]{Okuyama1981}){\bf.}
{\it
Let $S$ be a simple $kG$-module with vertex $P$, and
let $f$ 
be the Green correspondence
with respect to $(G, P, N_G(P))$.
If $S$ is a trivial source module, then its Green
correspondent $f(S)$ is again simple as $kN_G(P)$-module. 
}

\bigskip\noindent
{\bf 2.3.Lemma} (Scott 
\cite[II Theorem 12.4 and I Proposition 14.8]{Landrock} and 
\cite[Corollary 3.11.4]{Benson}){\bf.}
{\it
\begin{enumerate}
  \renewcommand{\labelenumi}{(\roman{enumi})}
    \item If $M$ is a trivial source $kG$-module,
then $M$ uniquely (up to isomorphism) 
lifts to a trivial source
$\mathcal O G$-lattice $\widehat M$.
    \item If $M$ and $N$ are both trivial source $kG$-modules,
then $[M,N]^G = (\chi_{\widehat M}, \chi_{\widehat N})^G$.
\end{enumerate}
}

\bigskip\noindent
{\bf 2.4.Lemma} (Fong-Reynolds){\bf.}
{\it
Let $H$ be a normal subgroup of $G$, and let 
$A$ and $B$ be block algebras of $\mathcal OG$ and
$\mathcal OH$, respectively, such that $A$ covers $B$.
Let $T = T_G(B)$ be the inertial subgroup (stabiliser)
of $B$ in $G$. Then, there is a block algebra
$\tilde A$ of $\mathcal OT$ such that
$\tilde A$ covers $B$,
$1_A 1_{\tilde A} = 1_{\tilde A} 1_A = 1_{\tilde A}$, 
$A = \tilde A ^G$ (block induction), and the block
algebras $A$ and $\tilde A$ are Morita equivalent via a
pair 
$(1_A{\cdot}{\mathcal OG}{\cdot}1_{\tilde A}, \, 
  1_{\tilde A}{\cdot}{\mathcal OG}{\cdot}1_A)$, 
that is, the Morita equivalence is a Puig equivalence
and induces a bijection
$$
{\mathrm{Irr}}(\tilde A) \rightarrow
{\mathrm{Irr}}(A),   \quad
\tilde{\chi} \mapsto \tilde{\chi}{\uparrow}^G;
\qquad
{\mathrm{Irr}}(A) \rightarrow
{\mathrm{Irr}}(\tilde A),   \quad
\chi \mapsto {\chi}{\downarrow}_T{\cdot}1_{\tilde A}
$$
between $\mathrm{Irr}(\tilde A)$ and $\mathrm{Irr}(A)$,
and a bijection
$$
{\mathrm{IBr}}(\tilde A) \rightarrow
{\mathrm{IBr}}(A),   \quad
\tilde{\phi} \mapsto \tilde{\phi}{\uparrow}^G;
\qquad
{\mathrm{IBr}}(A) \rightarrow
{\mathrm{IBr}}(\tilde A),   \quad
\phi \mapsto {\phi}{\downarrow}_T{\cdot}1_{\tilde A}
$$
between $\mathrm{IBr}(\tilde A)$ and $\mathrm{IBr}(A)$,
}

\bigskip\noindent
{\bf Proof.} See \cite[1.5.Theorem]{KoshitaniKunugiWaki2004} and
\cite[Chapter 5  Theorem 5.10]{NagaoTsushima}. \quad$\blacksquare$

\bigskip\noindent
{\bf 2.5.Lemma.}
{\it
Let $A$ be a block algebra of $\mathcal OG$ with a defect group $P$,
let $N = N_G(P)$, and let $A_N$ be a block algebra of $\mathcal ON$
which is the Brauer correspondent of $A$. Moreover,
let $(P,e)$ be a maximal $A$-Brauer pair, 
$H = N_G(P,e)$, the normaliser of $(P,e)$ in $N_G(P)$,
and let $B$ be a block algebra of $\mathcal OH$ which is the
Fong-Reynolds correspondent of $A_N$, see {\rm\bf{2.4}}.
Then, 
$A{\downarrow}^{G \times G}_{G \times H}{\cdot}1_B
= 1_A{\cdot}\mathcal OG{\cdot}1_B$,
as a right $\mathcal O [G \times H]$-module,
has a unique (up to isomorphism)
indecomposable direct summand with
vertex $\Delta P$.
}

\bigskip\noindent
{\bf Proof.} See \cite[Lemma 2.4]{KoshitaniKunugiWaki2008} and
\cite[Chapter 5, Theorem 5.10]{NagaoTsushima}. \quad$\blacksquare$

\bigskip\noindent
{\bf 2.6.Lemma.}
{\it
Assume that $G \geqslant H$, and let $A$ and $B$
respectively be block algebras of $\mathcal OG$ 
and $\mathcal OH$ with a common
defect group $P$, and hence $P \leqslant H$.
Suppose, moreover, that a pair $(M, M^\vee)$
induces 
a splendid stable equivalence of Morita type 
between $A$ and $B$,
where $M$ is an $(A,B)$-bimodule such that
$M \Big| {1_A}{\cdot}{\mathcal OG}{\cdot}1_B$
as $(A,B)$-bimodules.
\begin{enumerate}
\renewcommand{\labelenumi}{(\roman{enumi})}
    \item
If $X$ is a non-projective trivial source $kG$-module in $A$, 
then 
$(X \otimes_A M)_B = Y \oplus
 ({\mathrm{proj}})$
for a non-projective indecomposable $kH$-module $Y$
such that $Y$ has a trivial source.
    \item
If $X$ is a non-projective indecomposable $kG$-module in $A$, 
then
$(X \otimes_A M)_B = Y \oplus
 ({\mathrm{proj}})$
for a non-projective indecomposable $kH$-module $Y$
such that
there is a $p$-subgroup $Q$ of $H$ such that
$Q$ is a common vertex of $X$ and $Y$.
\end{enumerate}
}

\bigskip\noindent
{\bf Proof.} See \cite[Lemma 2.7]{KoshitaniKunugiWaki2004}.
\quad$\blacksquare$

\bigskip\noindent
{\bf 2.7.Lemma.} 
{\it
Let $k$ be a field, and let $A$ be a finite-dimensional
symmetric $k$-algebra. Moreover, suppose that
$S$ is a simple $A$-module and $M$ is a projective-free
$A$-module. Then, we have
$\underline{\mathrm{Hom}}_A(S,M) \cong {\mathrm{Hom}}_A(S,M)$
and 
$\underline{\mathrm{Hom}}_A(M,S) \cong {\mathrm{Hom}}_A(M,S)$
as $k$-spaces.
}

\bigskip\noindent
{\bf Proof.}
Follows by \cite[(3.2), (3.2*), (3.3)]{Green},
see \cite[II, Lemma 2.7, Corollary 2.8]{Landrock}.
\quad$\blacksquare$

\bigskip\noindent
{\bf 2.8.Lemma.} {\it
Let $k$ be an algebraically closed field,
and let $A$ and $B$ be finite-dimensional symmetric $k$-algebras.
Suppose that $M$ is an $(A,B)$-bimodule such that
$_A M$ and $M_B$ are both projective modules.
Then a functor 
$F: {\mathrm{mod}}{\text{-}}A \rightarrow
    {\mathrm{mod}}{\text{-}}B$
defined by $F(X') = X' \otimes_A M$ for ${X'}_A$,
is additive and exact. Assume, furthermore, that
$F$ induces a stable equivalence between $A$ and $B$.
\begin{enumerate}
\renewcommand{\labelenumi}{(\roman{enumi})}
    \item
Let $X$ be a projective-free $A$-module such that
$X$ has a simple $A$-submodule $S$. Set $T = F(S)$.
Then, we can write $F(X) = Y \oplus R$ for a projective-free
$B$-module $Y$ and a projective $B$-module $R$.
Now, if $T$ is a simple $B$-module, then we may assume that
$Y$ contains $T$ and
that $F(X/S)= Y/T \oplus ({\mathrm{proj}})$.
    \item
{\rm{(dual of (i))}} \ 
Let $X$ be a projective-free $A$-module such that
$X$ has an $A$-submodule $X'$ satisfying that
$X/X'$ is simple.
Set $T = F(X/X')$.
Then, we can write $F(X) = Y \oplus R$ 
for a projective-free $B$-module $Y$ and a projective
$B$-module $R$.
Now, if $T$ is a simple $B$-module, then we may assume that
$T$ is an epimorphic image of $Y$ and
that 
${\mathrm{Ker}}(F(X) \twoheadrightarrow T) 
= {\mathrm{Ker}}(Y \twoheadrightarrow T) \oplus ({\mathrm{proj}})$. 
\end{enumerate}
}

\bigskip\noindent
{\bf Proof.}
We get (i) from {\bf 2.7} and 
\cite[1.11.Lemma]{KoshitaniKunugiWaki2004},
just as in the proof of
\cite[3.25.Lemma and 3.26.Lemma]{KoshitaniKunugiWaki2004},
see \cite[Proposition 2.2]{Kunugi}.
(ii) is just the dual of (i).
\quad$\blacksquare$

\bigskip\noindent
{\bf 2.9.Lemma} (Linckelmann
\cite[Theorem 2.1(ii)]{Linckelmann1996MathZ}){\bf .} 
{\it
Let $A$ and $B$ be finite-dimensional $k$-algebras
for a field $k$ such that $A$ and $B$ are both
self-injective and indecomposable as algebras, and
none of them are simple algebras.
Suppose that there is an indecomposable $(A,B)$-bimodule $M$ such that
a pair $(M, M^\vee)$ induces a stable equivalence
between $A$ and $B$.
If $S$ is a simple $A$-module, then
$(S \otimes_A M)_B$ is a non-projective indecomposable
$B$-module.
}

\bigskip\noindent
The next lemma is a new result due to Kunugi and the first
author. This is actually so useful and convenient
when we want to apply so-called "Rouquier's glueing"
to our inductive argument in order to get 
a stable equivalence between two block algebras
which we are looking at.

\bigskip\noindent
{\bf 2.10.Lemma} 
(Koshitani-Kunugi
\cite[Theorem 1.2]{KoshitaniKunugi2006}){\bf .}
{\it
Let $A$ be a block algebra of $\mathcal OG$ 
with a cyclic defect group  $P {\not=} \; 1$. 
Let $H = N_G(P)$, and let $B$ be a block algebra of
$\mathcal OH$ such that $B$ is the Brauer correspondent
of $A$. 
Then, we get the following:
\begin{enumerate}
  \renewcommand{\labelenumi}{(\roman{enumi})}
\item
The following {\rm{(1)}} and {\rm{(2)}} are equivalent:
  \begin{enumerate}
  \renewcommand{\labelenumii}{(\arabic{enumii})}
    \item
The Brauer tree of $A$ is a star with exceptional vertex
in the centre, and there exists a non-exceptional 
irreducible ordinary character $\chi$ of $G$ in $A$ such that
$\chi (u) > 0$ for any element $u \in P$.
    \item
The block algebras $A$ and $B$ are Puig equivalent.
  \end{enumerate}
\item
If one of the conditions {\rm{(1)}} and {\rm{(2)}} in {\rm{(i)}}
holds (and hence both hold), then all simple $kG$-modules in $A$
are trivial source modules.
\item
If one of the conditions {\rm{(1)}} and {\rm{(2)}} in {\rm{(i)}}
holds (and hence both hold),
then there is an indecomposable
$(A, B)$-bimodule $\mathfrak M$ such that
$1_A{\cdot}{\mathcal OG}{\cdot}1_{B}$
$ = \mathfrak M \oplus ({\mathrm{proj}})$
and $\mathfrak M$, as an $\mathcal O[G \times H]$-module, 
has $\Delta P$ as its vertex, and $\mathfrak M$ 
realizes a Puig equivalence
between $A$ and $B$.
\end{enumerate}
}

\bigskip\noindent
{\bf 2.11.Lemma.} 
{\it
Let $A$ be a block algebra of $\mathcal OG$ with defect
group $P$. Set $H = N_G(P)$, and let $B$ be a block
algebra of $\mathcal OH$ such that $B$ is the Brauer
correspondent of $A$.
Assume that $Q$ is a subgroup of $P$ with
$Q \subseteq Z(G)$.
Set $\bar G = G/Q$, $\bar H = H/Q$ and $\bar P = P/Q$.
It is well-known that there exist block algebras
$\bar A$ and $\bar B$ of $\mathcal O{\bar G}$ and
$\mathcal O{\bar H}$, respectively, such that
$\bar A$ and $\bar B$ dominate $A$ and $B$, namely
${\mathrm{Irr}}(\bar A) \subseteq {\mathrm{Irr}}(A)$
and
${\mathrm{Irr}}(\bar B) \subseteq {\mathrm{Irr}}(B)$,
and that
both $\bar A$ and $\bar B$ have $\bar P$ as defect groups,
see 
{\rm{ \cite[Chapter 5 Theorems 8.10 and 8.11]{NagaoTsushima}}}.
\begin{enumerate}
  \renewcommand{\labelenumi}{(\roman{enumi})}
\item
It holds that $\bar H = N_{\bar G}(\bar P)$ 
and that $\bar B$ is the
Brauer correspondent of $\bar A$.

\smallskip\noindent
In the rest of the lemma, assume in addition that $P$ is
elementary abelian of order $p^2$, namely,
$P = Q \times R$ with $Q \cong R \cong C_p$.

\smallskip\noindent
\item
It holds that
$$
\bar A \otimes_{\mathcal O{\bar H}} \bar B
=
{_{\bar A}}{(\bar A{\cdot}1_{\bar B})_{\bar B}}
=
{_{\bar A}}{X_{\bar B}} \oplus ({\mathrm{proj}})
$$
for an indecomposable $(\bar A, \bar B)$-bimodule $X$
with vertex $\Delta \bar P$.
\item
In particular, if $X$ realizes a Morita equivalence 
between $\bar A$ and $\bar B$, then there exists an
$(A,B)$-bimodule $M$ such that $M$ is an indecomposable
direct summand of 
${_{A}}{(A{\cdot}1_{B})_{B}}$
with vertex $\Delta P$, and hence $M$ induces a
Puig equivalence between $A$ and $B$.
\end{enumerate}
}

\bigskip\noindent
{\bf Proof.}
(i) The first part is easy. The second part follows
from \cite[(3.2)Lemma]{Navarro2004}, see
\cite[$\ell$.10 on p.1314]{Linckelmann2009}.

(ii) This follows by \cite[Proposition 6.1]{Linckelmann1996Inv}
since $\bar P \cong C_p$.

(iii) This is obtained from {\bf (ii)} and
\cite[Theorem]{KoshitaniKunugi2005},
see 
\cite[$\ell. -7 \sim \ell. -4$ on p.1314]{Linckelmann2009}
and \cite[Theorem 4.1]{Linckelmann2001}.
\quad$\blacksquare$

\bigskip\noindent
{\bf 2.12.Lemma.} 
{\it
Suppose that $p = 3$ and $G = \mathcal A_9$.
\begin{enumerate}
  \renewcommand{\labelenumi}{(\roman{enumi})}
\item
There uniquely exists a non-principal block algebra
$A$ of $\mathcal OG$ with defect group $P \cong C_3$.
In addition we can write
${\mathrm{Irr}}(A) = \{ \chi_5, \chi_{17}, \chi_{18} \}$
such that
$\chi_5(1) = 27$, $\chi_{17}(1) = 189$, $\chi_{18}(1) = 216$
and 
$\chi_5(u) = \chi_{17}(u) = 9$
for any element $u \in P - \{1\}$,
and that a part of the $3$-decomposition matrix is
\begin{center}
{
{\rm
\begin{tabular}{l|cc}
            & $27$ & $189$   \\
\hline
$\chi_{5}$  & 1  &  0 \\
$\chi_{17}$ & 0  &  1 \\
$\chi_{18}$ & 1  &  1 \\
\end{tabular} 
}}
\end{center}
where the indices of $\chi_i$ are the same as in
\cite[p.37]{Atlas}.
(In the following, we use the notation $A$ and $P$ 
as in {\bf (i)}).
\item
%
%JM: Replaced , by ., and making the statement precise.
Set $H = N_G(P)$. Then
$H = (P \times \mathcal A_6).C_2$, where the action on
$P \times \mathcal A_6$ by $C_2$ is the diagonal one,
extending $\mathcal A_6$ to $\mathcal S_6$.
\item
Let $H$ be as in {\bf (ii)}, and let $B$ be a block algebra
of $\mathcal OH$, which is the Brauer correspondent of $A$.
Then, $A$ and $B$ are Morita equivalent via an $(A,B)$-bimodule
$M$ such that $M$ is (up to isomorphism) the unique
indecomposable direct summand of 
${_A}{(A{\cdot}1_B)_B}$ with vertex $\Delta P$,
and hence it holds that $M$ induces a Puig equivalence
between $A$ and $B$, and that the simples
$27$ and $189$ in $A$ are both trivial source
$kG$-modules.
\end{enumerate}
}

\bigskip\noindent
{\bf Proof.}
(i) This follows from \cite[p.37]{Atlas},
\cite[$A_9$ (mod 3)]{ModularAtlasProject}
and
\cite{ModularAtlas}.

%JM: this is obvious
%Easy by direct calculations by {\sf GAP} \cite{GAP}.
(ii) Easy by inspection.

(iii) This is obtained from {\bf (i)} and {\bf 2.10}.
\quad$\blacksquare$

\bigskip\noindent
{\bf 2.13.Lemma.} 
{\it
Let $A$ and $B$ be finite dimensional $k$-algebras.
Assume that there exists a functor
$F: {\underline{\mathrm{mod}}}{\text{-}}A \rightarrow
    {\underline{\mathrm{mod}}}{\text{-}}B$
realizing a stable equivalence between $A$ and $B$.
Assume, in addition, that there is 
a simple $A$-module $S_0$ such that $S_0$ is 
sent to a simple $B$-module $T_0$, namely, 
$F(S_0) = T_0$.
Then, for any simple $A$-module $S$ with
$S \not\cong S_0$, it holds
$[F(S), T_0]^B = [T_0, F(S)]^B = 0$.
}

\bigskip\noindent
{\bf Proof.}
We get by {\bf 2.7} and the assumptions that
\begin{align*}
   0   &= {\mathrm{Hom}}_A(S, S_0) 
        \cong {\underline{\mathrm{Hom}}}_A(S, S_0)       \\
       &\cong {\underline{\mathrm{Hom}}}_B(F(S), F(S_0)) 
        =     {\underline{\mathrm{Hom}}}_B(F(S), T_0)    \\
        &\cong                {\mathrm{Hom}}_B(F(S), T_0).
\end{align*}
Hence $[F(S), T_0]^B = 0$. The rest is similar.
\quad$\blacksquare$               

%\bigskip\noindent
%{\bf 2.15??.Lemma.}
%{\it
%Let $X$ be an $A$-module such that
%$X = P_1 \oplus \cdots \oplus P_n \oplus Y$
%for a positive integer $n$ where each $P_i$ is an
%injective indecomposable $A$-module and $Y$ is an
%$A$-module.
%If $Z$ is an $A$-submodule of $X$ satisfying that
%$S = {\mathrm{soc}}(Z)$ is simple and 
%$I(S) \cong P_1$, then we may assume that $Z \subseteq P_1$.
%}
%
%\bigskip\noindent
%{\bf Proof.} 
%Set $X' = I(Z) \oplus P_2 \oplus \cdots \oplus P_n \oplus Y$.
%Then we get $X' \cong X$ since $I(Z) = I(S) \cong P_1$.
%Hence, we obtain the assertion by replacing $X$ by $X'$.
%\quad$\blacksquare$

\bigskip\noindent
{\bf 2.14.Lemma.} 
{\it
Let $A$ be a finite-dimensional $k$-algebra, and assume that
$X$ is an $A$-module satisfying that
$X = P_1 \oplus \cdots \oplus P_n \oplus Y \supseteq Z$
for an integer $n \geqslant 2$, 
$A$-submodules $P_1, \cdots, P_n$ and $Y$ of $X$ such that
$P_1 \cong \cdots \cong P_n \cong P(S)$ 
for a simple $A$-module $S$,
and ${\mathrm{soc}}(Z) \cong S$, and $Y_A$ 
is projective-free. Assume moreover that $j(P(S)) = j(A)$.
Then, $X/Z$ has a direct summand isomorphic to $P(S)$.
}

\bigskip\noindent
{\bf Proof.}
Set $j = j(A) = j(P(S))$, and $\bar X = X/Z$.
Assume that $\bar X$ is projective-free.
Then, ${\mathrm{rad}}^{j-1}(\bar X) = 0$, and hence
${\mathrm{rad}}^{j-1}
(P_1 \oplus \cdots \oplus P_n \oplus Y) \subseteq Z$.
Clearly,
${\mathrm{rad}}^{j-1}
(P_1 \oplus \cdots \oplus P_n \oplus Y) 
 = {\mathrm{soc}}(P_1) \oplus \cdots 
\oplus {\mathrm{soc}}(P_n) \cong S \oplus \cdots \oplus S$ 
($n$ times)
since ${\mathrm{rad}}^{j-1}(Y) = 0$.
This means that 
$(S \oplus S \oplus \cdots \oplus S) \, | \, {\mathrm{soc}}(Z)$, 
a contradiction.
Thus, $P(T) | {\bar X}$
for a simple $A$-module $T$. 
This implies that there is an
epimorphism $X \twoheadrightarrow P(T)$, and hence
$P(T) | X$. Then, we get $P(T) \cong P(S)$ by Krull-Schmidt
theorem.
\quad$\blacksquare$

\bigskip\noindent
{\bf 2.15.Lemma.} 
{\it
Let $G$, $H$ and $L$ be finite groups such that all of them
contain a common subgroup $P$, namely, 
$P \subseteq G \cap H \cap L$.
Let $M$ be a $k[G \times H]$-module such that
$M \, \Big| \, {k_{\Delta P}}{\uparrow^{G \times H}}$,
and 
let $N$ be a $k[H \times L]$-module such that
$N \, \Big| \, {k_{\Delta P}}{\uparrow^{H \times L}}$.
Then, it follows that
\linebreak
$M \otimes_{kH} N \, \Big| \, {k_{\Delta P}}{\uparrow^{G \times L}}$.
}

\bigskip\noindent
{\bf Proof.}
This is a special case of 
\cite[2.5.Proposition]{HarrisKoshitani}.
\quad$\blacksquare$

\bigskip\noindent
{\bf 2.16.Lemma.} 
{\it
Let $A$ be a finite-dimensional $k$-algebra,
and assume that $X$ is an indecomposable non-simple
$A$-module. Then, it holds 
${\mathrm{soc}}(X) \subseteq {\mathrm{rad}}(X)$.
}

\bigskip\noindent
{\bf Proof.}
Assume that 
${\mathrm{soc}}(X) \not\subseteq {\mathrm{rad}}(X)$.
Then, $X$ has a simple $A$-submodule $S$ with
$S \not\subseteq {\mathrm{rad}}(X)$.
Hence, $X$ has a maximal $A$-submodule $M$ with
$S \not\subseteq M$.
These imply that $S \cap M = 0$ and 
$S + M = X$. Namely, $X = S \oplus M$.
Since $M$ is indecomposable, $X = S$.
This is a contradiction.
\qquad$\blacksquare$

\bigskip
%\newpage

\begin{flushleft}
{\bf 3. $3$-Local structure for {\sf HN}}
\end{flushleft}

\bigskip\noindent
{\bf 3.1.Notation and assumption.} From 
now on, we assume that $G$ is the Harada-Norton simple
group {\sf HN}, and hence
$|G| = 2^{14}{\cdot}3^6{\cdot}5^6{\cdot}7{\cdot}11{\cdot}19
 \ \fallingdotseq \ 2.7 \times 10^{14}$,
see \cite[p.164--166]{Atlas} and \cite{Harada}.

\bigskip\noindent
{\bf 3.2.Lemma.}
{\it
\begin{enumerate}
  \renewcommand{\labelenumi}{(\roman{enumi})}
    \item
In order to prove 
Brou{\'e}'s abelian defect group conjecture for $G = {\sf HN}$,
it suffices to prove it for the case $p = 3$.
   \item
There exists a unique $3$-block $A$  with non-cyclic abelian 
defect group $P$, and $P$ is elementary abelian of order $9$.
%JM: identify P
   \item
$P$ is the Sylow $3$-subgroup of the second largest maximal subgroup
$2.{\sf HS}.2$ of $G$, a two-fold cover of the 
automorphism group of the Higman-Sims simple group {\sf HS}.
\end{enumerate}
}

\bigskip\noindent
{\bf Proof.}
(i)
We may assume $p \in \{ 2,3,5 \}$ by {\rm{\bf 3.1}} just as
in the proof of \cite[Lemma 3.2]{KoshitaniKunugiWaki2008}.
Assume that $p = 2$. Then, $G$ has only two $2$-blocks
$B_0$ and $B_1$ with positive defect by
\cite{ModularAtlasProject},
where $B_0$ is the principal $2$-block.
Then, the non-principal $2$-block $B_1$ has a defect group
$D$ with $D \cong SD_{16}$, see
\cite[Lemma 4.2(c)]{AnObrien}.
Thus, $B_0$ and $B_1$ both have non-abelian defect groups.
Next, suppose $p = 5$. By \cite{ModularAtlasProject}, 
$G$ has only a unique $5$-block $B_0$ which has defect
$\geqslant 2$, and hence $B_0$ is the principal $5$-block.
Then, $B_0$ has non-abelian defect group $5_{+}^{1+4}.5$
by \cite[p.164--166]{Atlas}.

(ii)
Finally, assume $p = 3$. Sylow $3$-subgroups of $G$ are
non-abelian by \cite[p.164--166]{Atlas}. Thus, $G$ has a
unique non-principal $3$-block $A$ such that
$A$ has a defect group $P$
with $|P| \geqslant 3^2$, 
and actually $P \cong C_3 \times C_3$,
see \cite[Lemma 4.2(b)]{AnObrien}.

(iii)
Using the character table of $G$,
calculations with {\sf GAP} \cite{GAP} show that 
the conjugacy class {\rm{2A}} of $G$ is a defect class of $A$,
where we follow the notation in \cite[p.164--166]{Atlas}.
Hence $P$ is a Sylow $3$-subgroup of the centralizer
$C_G(2A)\cong 2.{\sf HS}.2$.
\quad$\blacksquare$

\bigskip\noindent
%JM: `From' should not be at the beginning of a line.
{\bf 3.3.Notation.} From now on, we assume $p = 3$, and we use the notation
$A$ and $P$ as in {\rm{\bf} 3.2}, namely, $A$ is a block
algebra of $kG$ with defect group $P \cong C_3 \times C_3$.
Set $H = N_G(P)$, and let $B$ be a block algebra of $kH$
that is the Brauer correspondent of $A$.
Let $(P, e)$ be a maximal $A$-Brauer pair in $G$, that it,
$e$ is a block idempotent of $k C_G(P)$ such that
${\mathrm{Br}}_P (1_A){\cdot}e = e$,
see \cite{AlperinBroue}, \cite{BrouePuig1980} and
\cite[\S 40]{Thevenaz}.
Set ${\widetilde H} = N_G(P, e)$, namely,
${\widetilde H} = \{ g \in N_G(P) \ | \ e^g = e \}$,
where $e^g = g^{-1}eg$.
%Let $\widetilde B$ be a block algebra of $k\widetilde H$
%which is the Fong-Reynolds correspondent of $B$,
%see {\rm{\bf 2.4}}.
Finally set $E = H/C_G(P)$, and let $Q$ be a subgroup of $P$
of order $3$.

\bigskip\noindent
{\bf 3.4.Lemma.}
{\it
It holds the following:
\begin{enumerate}
  \renewcommand{\labelenumi}{(\roman{enumi})}
    \item
$H = \widetilde H = (P \times {\mathcal A}_6).SD_{16}$.
    \item
$C_G(P) = C_H(P) = P \times {\mathcal A}_6$.
    \item
$E = {\widetilde H}/C_G(P) \cong SD_{16}$,
where the action of $E$ on $P$ is given by the embedding of 
$SD_{16}$ as a Sylow $2$-subgroup of $\mathrm{Aut}(P)\cong\mathrm{GL}_2(3)$.
    \item
All elements in $P - \{ 1 \}$ are conjugate in $H$, and hence in $G$.
    \item
$P - \{ 1 \} \subseteq {\mathrm{3A}}$, 
where {\rm{3A}} is a conjugacy class of $G$ following the notation
in \cite[p.164--166]{Atlas}.
    \item
All subgroups of $P$ of order $3$ are conjugate in $H$, and hence
in $G$.
    \item
Recall the subgroup $Q$ of $P$ in {\rm{\bf 3.3}}. Then, we have
$C_G(Q) = Q \times {\mathcal A}_9$,
$N_G(Q) = (Q \times {\mathcal A}_9).2 \leqslant {\mathcal A}_{12}$,
$C_H(Q) = (P \times {\mathcal A}_6).2$, and
$N_H(Q) = (P \times {\mathcal A}_6).2^2$.
    \item
$C_G(Q)/Q \cong {\mathcal A}_9$,
$C_H(Q)/Q \cong (C_3 \times {\mathcal A}_6).2$, 
%    = \Sigma_3.{\mathcal A}_6$,
$N_G(Q)/Q \cong {\mathcal A}_9.2$, and
$N_H(Q)/Q \cong (C_3 \times {\mathcal A}_6).2^2$.
\end{enumerate}
}

\bigskip\noindent
{\bf Proof.}
%JM: more comments on computations, revised the proof
This is found using explicit computation with {\sf GAP} \cite{GAP}.
The starting point is the smallest faithful permutation representation 
of $G$ on $1140000$ points, available in terms of so-called
standard generators \cite{Wilson} in \cite{ModAtlasRep}.
The associated one-point stabiliser is the largest maximal 
subgroup $\mathcal A_{12}$ of $G$, which hence can be found explicitly
by a randomised Schreier-Sims technique. Having completed that,
all the following computations can be done using this permutation
representation of $G$.

Actually, one of the standard generators is an element of the $2A$
conjugacy class of $G$, where we use the notation in \cite[p.164--166]{Atlas}.
Hence the second largest maximal subgroup $2.{\sf HS}.2\cong C_G(2A)$
can be found be a centraliser computation. In turn, by {\bf 3.2(iii)}
$P$ can be computed explicitly as a Sylow $3$-subgroup of $2.{\sf HS}.2$.

(i)--(ii)
The normaliser $H=N_G(P)$ and the centraliser $C_G(P)$ of $P$ can be
computed explicitly, and as these are fairly small groups 
their structure is easily revealed.
 
(iii)
It follows from \cite[$A_6$ (mod 3)]{ModularAtlasProject} 
and \cite{ModularAtlas} that
$\mathcal A_6$ has exactly two $3$-blocks.
Let $\beta$ be the non-principal block algebra
of $k{\mathcal A_6}$, and hence $\beta$ is of
defect zero. Then, $e = 1_{\beta}$.
Since $\beta$ is a unique block algebra of
$k{\mathcal A_6}$ of defect zero, this shows $H = \widetilde H$.

%(ii) We get this from {\bf (i)} and its proof.
%(iii) Easy by {\bf (i)} and {\bf (ii)}.

(iv) Easy by {\bf (iii)} and inspection.

(v) We use the notation $3A$ and $3B$ as in \cite[p.164--166]{Atlas}.
By {\bf (iv)}, $P - \{1\} \subseteq {\mathrm{3A}}$ or ${\mathrm{3B}}$.
Assume $P - \{1\} \subseteq {\mathrm{3B}}$.
Then, $\chi (u) = 0$ for any $\chi \in {\mathrm{Irr}}(A)$ and
any $u \in {\mathrm{3A}}$ by
\cite[Chapter~5 Corollary 1.10(i)]{NagaoTsushima}.
But we know that $\chi_8 \in {\mathrm{Irr}}(A)$ by 
%using {\sf GAP} \cite{GAP} (see also 
\cite[Lemma 4.2(b)]{AnObrien}, see also {\bf 4.1},
and that
$\chi_8(u) = 27$ for any $u \in {\mathrm{3A}}$.
%by \cite[p.164--166]{Atlas}. 
This is a contradiction.

(vi) Easy by {\bf (iv)}.

(vii)--(viii) 
It is easy to see that $N_G(Q)<\mathcal A_{12}$, the largest maximal
subgroup of $G$, which is the one-point stabiliser in the
given permutation representation of $G$. Hence again the normaliser 
$N_{\mathcal A_{12}}(Q)$ and the centraliser $C_{\mathcal A_{12}}(P)$ 
of $P$ can be computed explicitly and their structure determined.
\quad$\blacksquare$

\bigskip\noindent
{\bf 3.5.Lemma.}
{\it We get the following diagram:}
%JM: The index of N_G(Q) in A12 indeed is 220. 
%    It was just a mistake in my old notes.
%
\bigskip
{
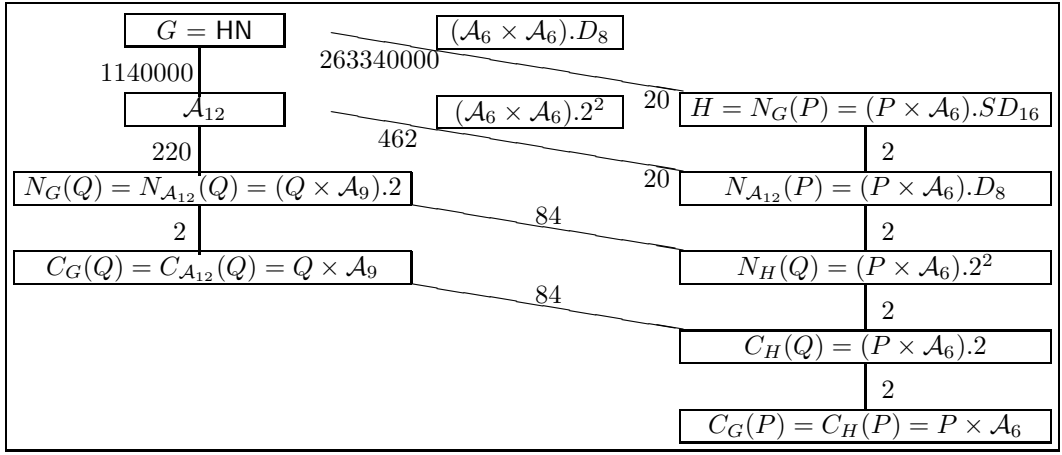
\begin{center}
{
\boxed{
\begin{picture}(392,163)(-90,-150)
%(200,100)(-100,-100)
%JM: Rearranged boxes slightly
%
\put(-48,0){\framebox(60,12){$G$ = {\sf HN}}}
\put(-20,0){\line(0,-1){19}}
\put(-48,-30){\framebox(60,12){$\mathcal A_{12}$}}
\put(-20,-30){\line(0,-1){19}}
\put(-90,-60){\framebox(150,12){$N_G(Q) = N_{\mathcal A_{12}}(Q) 
                                        = (Q \times {\mathcal A_9}).2$}}
\put(-20,-60){\line(0,-1){19}}
%JM: Corrected a subscript.
\put(-90,-90){\framebox(150,12){$C_G(Q) = C_{\mathcal A_{12}}(Q) 
                                        = Q \times {\mathcal A_9}$}}
\put(162,-30){\framebox(140,12){$H = N_G(P) 
                                   = (P \times {\mathcal A_6}).SD_{16}$}}
\put(232,-30){\line(0,-1){18}}
\put(162,-60){\framebox(140,12){$N_{\mathcal A_{12}}(P) 
                           = (P \times {\mathcal A_6}).D_8$}}
\put(232,-60){\line(0,-1){18}}
\put(162,-90){\framebox(140,12){$N_H(Q) = (P \times {\mathcal A_6}).2^2$}}
\put(232,-90){\line(0,-1){18}}
\put(162,-120){\framebox(140,12){$C_H(Q) = (P \times {\mathcal A_6}).2$}}
\put(232,-120){\line(0,-1){18}}
\put(162,-150){\framebox(140,12){$C_G(P) = C_H(P) = P \times {\mathcal A_6}$}}
\put(30,5){\line(6,-1){133}}
\put(30,-25){\line(6,-1){133}}
\put(60,-60){\line(6,-1){103}}
\put(60,-90){\line(6,-1){103}}
\put(70,-1)
  {\framebox(70,12){$(\mathcal A_{6} \times \mathcal A_6).D_8$}}
\put(70,-31)
  {\framebox(70,12){$(\mathcal A_{6} \times \mathcal A_6).2^2$}}
\put(97,-69){\makebox(30,10){$84$}}
\put(97,-99){\makebox(30,10){$84$}}
\put(-65,-15){\makebox(50,10){1140000}}
\put(-38,-45){\makebox(15,10){220}}
\put(-30,-75){\makebox(5,10){2}}
\put(238,-45){\makebox(5,10){2}}
\put(238,-75){\makebox(5,10){2}}
\put(238,-105){\makebox(5,10){2}}
\put(238,-135){\makebox(5,10){2}}
\put(148,-25){\makebox(10,10){20}}
\put(148,-55){\makebox(10,10){20}}
\put(18,-10){\makebox(60,10){263340000}}
\put(45,-40){\makebox(20,10){462}}
\end{picture}
}}
\end{center}
}

%\vspace{5.5cm}
\bigskip\noindent
{\bf Proof.} This follows from \cite[p.164--166]{Atlas}, {\bf 3.4}
and calculations with {\sf GAP} \cite{GAP}.
\quad$\blacksquare$

\bigskip\noindent
{\bf 3.6.Lemma.}
{\it
The following holds:
\begin{enumerate}
  \renewcommand{\labelenumi}{(\roman{enumi})}
    \item
$B \cong {\mathrm{Mat}}_9({\mathcal O}[P \rtimes SD_{16}])$
as $\mathcal O$-algebras, 
    \item
The block algebra $B$ has a source algebra
$jBj \cong {\mathcal O}[P \rtimes SD_{16}]$,
as interior $P$-algebras,
where $j$ is a source idempotent of $B$ with 
respect to $P$, namely, $j$ is a primitive idempotent
of $B^P$ such that ${\mathrm{Br}}_{P}(j) \not= 0$
for the Brauer homomorphism  ${\mathrm{Br}}_P$ for $P$,
see \cite[\S\S 19 and 27]{Thevenaz}.
\smallskip
    \item
We can write
$$
{\mathrm{Irr}}(B) = 
\{ \chi_{9a}, \chi_{9b}, \chi_{9c}, \chi_{9d}, 
\chi_{18a}, \chi_{18b}, \chi_{18c}, \chi_{72a}, \chi_{72b} \}
$$
and 
$$
{\mathrm{IBr}}(B) = \{ 9a, 9b, 9c, 9d, 18a, 18b, 18c \},
$$
where the numbers mean the degrees of characters and
the dimensions of simples, respectively.
Note that $\chi_{18b}$ and $\chi_{18c}$ are dual each other,
and so are $18b$ and $18c$.
The other characters and simples are self-dual.
\smallskip
    \item
The $3$-decomposition matrix and the Cartan matrix of $B$
are %respectively 
the following:
\medskip
\begin{center}
{\rm
\begin{tabular}{l|ccccccc}
  & $9a$ & $9b$ & $9c$ & $9d$ & $18a$ & $18b$ & $18c$   \\
\hline
$\chi_{9a}$  & 1  &  .    & .     &  .    &  .  & . & .  \\
$\chi_{9b}$  & .  &  1    & .     &  .    &  .  & . & .  \\
$\chi_{9c}$  & .  &  .    & 1     &  .    &  .  & . & .  \\
$\chi_{9d}$  & .  &  .    & .     &  1    &  .  & . & .  \\
$\chi_{18a}$ & .  &  .    & .     &  .    &  1  & . & .  \\
$\chi_{18b}$ & .  &  .    & .     &  .    &  .  & 1 & .  \\
$\chi_{18c}$ & .  &  .    & .     &  .    &  .  & . & 1  \\
$\chi_{72a}$ & 1  &  1    & .     &  .    &  1  & 1 & 1  \\
$\chi_{72b}$ & .  &  .    & 1     &  1    &  1  & 1 & 1  \\
\end{tabular} 
}
\end{center}
\medskip
\begin{center}
{\rm
\begin{tabular}{r|ccccccc}
      & $P(9a)$ & $P(9b)$  & $P(9c)$  & $P(9d)$ 
      &$P(18a)$ & $P(18b)$ & $P(18c)$ \\
\hline
$9a$  & 2 & 1 & 0 & 0 & 1 & 1 & 1 \\ 
$9b$  & 1 & 2 & 0 & 0 & 1 & 1 & 1 \\ 
$9c$  & 0 & 0 & 2 & 1 & 1 & 1 & 1 \\ 
$9d$  & 0 & 0 & 1 & 2 & 1 & 1 & 1 \\ 
$18a$ & 1 & 1 & 1 & 1 & 3 & 2 & 2 \\ 
$18b$ & 1 & 1 & 1 & 1 & 2 & 3 & 2 \\
$18c$ & 1 & 1 & 1 & 1 & 2 & 2 & 3 \\
\end{tabular}
}
\end{center}
\medskip
    \item
There are unique conjugacy classes $4A$ and $4B$ of $H$, consisting of 
elements of order $4$, and having centralisers of order $40$
and $48$, respectively. A part of the character table of
${\mathrm{Irr}}(B)$ then is the following:
\begin{center}
{\rm
%JM: I also do not know why class 3A could be needed.
%    Seems to be a misunderstanding, sorry for that.
%    Deleted the fourth column again, old version kept below.
\begin{tabular}{l|rrr}

{conjugacy class} & ${4A}$ & ${4B}$ & ${12A}$ \\
{centraliser}& $40$ & $48$ & $24$ \\
\hline
\qquad $\chi_{9a}$  & $1$  &  $-1$    & $-1$ \\
\qquad $\chi_{9b}$  & $-1$ &  $-1$    & $-1$ \\
\qquad $\chi_{9c}$  &  $1$ &   $1$    & $1$  \\
\qquad $\chi_{9d}$  & $-1$ &   $1$    & $1$  \\
\qquad $\chi_{18a}$ &  $0$ &   $0$    & $0$  \\
\qquad $\chi_{18b}$ &  $0$ &   $0$    & $0$  \\
\qquad $\chi_{18c}$ &  $0$ &   $0$    & $0$  \\
\qquad $\chi_{72a}$ &  $0$  &  $-2$   & $1$  \\
\qquad $\chi_{72b}$ &  $0$  &  $2$    & $-1$ \\
\end{tabular} 
%\begin{tabular}{l|rrrr}
%
%{\small{conjugacy class}} & ${4B}$ & ${4C}$ 
%             & ${12B}$ & ${3A}$ \\
%{\small{size of centraliser}}& $40$ & $48$ & $24$ & $6480$  \\
%\hline
%\qquad $\chi_{9a}$  & $1$  &  $-1$    & $-1$  & $9$  \\
%\qquad $\chi_{9b}$  & $-1$ &  $-1$    & $-1$  & $9$  \\
%\qquad $\chi_{9c}$  &  $1$ &   $1$    & $1$   & $9$   \\
%\qquad $\chi_{9d}$  & $-1$ &   $1$    & $1$   & $9$   \\
%\qquad $\chi_{18a}$ &  $0$ &   $0$    & $0$   & $18$   \\
%\qquad $\chi_{18b}$ &  $0$ &   $0$    & $0$   & $18$   \\
%\qquad $\chi_{18c}$ &  $0$ &   $0$    & $0$   & $18$   \\
%\qquad $\chi_{72a}$ &  $0$  &  $-2$   & $1$   & $-9$   \\
%\qquad $\chi_{72b}$ &  $0$  &  $2$    & $-1$  & $-9$   \\
%\end{tabular} 
}
\end{center}
Note that this identifies the characters 
$\chi_{9a}$, $\chi_{9b}$, $\chi_{9c}$, $\chi_{9d}$,
$\chi_{72a}$ and  $\chi_{72b}$ uniquely.

\medskip
   \item 
The radical and socle series of PIMs in $B$ 
are the following: 
%JM: introduced boxes
{\rm
$$
\boxed{
\begin{matrix}
9a \\
18b \\
9b \ 18a \\
18c\\
9a 
\end{matrix}}
\qquad
\boxed{
\begin{matrix}
9b \\
18c \\
9a \ 18a \\
18b\\
9b 
\end{matrix}}
\qquad
\boxed{
\begin{matrix}
9c \\
18c \\
9d \ 18a \\
18b\\
9c 
\end{matrix}}
\qquad
\boxed{
\begin{matrix}
9d \\
18b \\
9c \ 18a \\
18c\\
9d
\end{matrix}}
$$
\smallskip
$$
\boxed{
\begin{matrix}
18a\\
18b \ 18c \\
9b \ 9c \ 18a \ 9a \ 9d\\
18c \ 18b\\
18a
\end{matrix}}
\qquad
\boxed{
\begin{matrix}
18b\\
9b \  18a \ 9c \\
18c \ 18b \ 18c\\
9a  \ 18a \ 9d\\
18b
\end{matrix}}
\qquad
\boxed{
\begin{matrix}
18c\\
9a \ 18a \ 9d \\
18b \ 18c \ 18b\\
9b \ 18a \ 9c\\
18c
\end{matrix}}
$$
}
Note that this identifies the simples $18b$ and $18c$ uniquely.

\medskip
   \item 
%JM: For the identification in 8.2
%JM: more cautious...
An Alperin diagram of the PIM $P(18a)$ is given as follows:
$$ P(18a) \ = \ \begin{matrix} \boxed
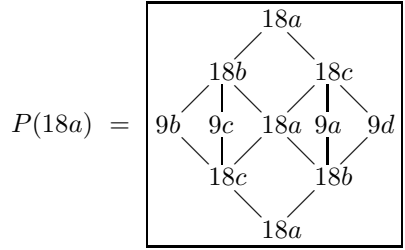
{
\begin{picture}(90,86)(0,0)
\put(40,80){$18a$}
\put(20,60){$18b$}
\put(60,60){$18c$}
\put(0,40){$9b$}
\put(20,40){$9c$}
\put(40,40){$18a$}
\put(60,40){$9a$}
\put(80,40){$9d$}
\put(20,20){$18c$}
\put(60,20){$18b$}
\put(40,0){$18a$}
\put(30,18){\line(1,-1){10}}
\put(50,8){\line(1,1){10}}
\put(10,38){\line(1,-1){10}}
\put(30,28){\line(1,1){10}}
\put(50,38){\line(1,-1){10}}
\put(70,28){\line(1,1){10}}
\put(10,48){\line(1,1){10}}
\put(30,58){\line(1,-1){10}}
\put(50,48){\line(1,1){10}}
\put(70,58){\line(1,-1){10}}
\put(50,78){\line(1,-1){10}}
\put(30,68){\line(1,1){10}}
\put(25,28){\line(0,1){10}}
\put(25,48){\line(0,1){10}}
\put(65,28){\line(0,1){10}}
\put(65,48){\line(0,1){10}}
\end{picture} } \end{matrix} $$

\end{enumerate}
}

\bigskip\noindent
{\bf Proof.} 
%JM: More comments on computations, revised the proof
This again relies on computations with {\sf GAP} \cite{GAP}.
Starting with the explicit restriction of the permutation
representation of $G$ to $H$ obtained in {\bf 3.4},
we find a faithful permutation representation of $H$ on a
small number of points. This then is used to compute the
conjugacy classes of $H$, and its ordinary character table
using the Dixon-Schneider algorithm. 

(i) Since the Schur multiplier of
$SD_{16}$ is trivial, see e.g. \cite[Proof of Lemma 1.3]{Koshitani1986},
we get the assertion by {\bf 3.4(i)-(iii)},
\cite[A.Theorem]{Kuelshammer1985}.
%and the degrees of $\chi \in {\mathrm{Irr}}(B)$
%which are obtained by {\sf GAP} \cite{GAP}.

(ii) This follows by a result of
%SK
Puig \cite[Proposition 14.6]{Puig1988JAlg} and {\bf (i)},
see \cite[Theorem 13]{AlperinLinckelmannRouquier}
and \cite[(45.12)Theorem]{Thevenaz}.

(iii)--(v) Easy from the character table of $H$.
%Later on, for instance in {\bf 5.4},
%we shall need to distinguish the four characters
%$\chi_{9a}, \cdots, \chi_{9d}$ explicitly.
%Hence we put this here.

(vi) The radical and socle series have been determined in \cite{Waki1989}.

(vii) 
To find the %full
structure of $P(18a)$,
we have used the {\sf MeatAxe} \cite{MA} to construct $P(18a)$
explicitly as a matrix representation, from the 
permutation representation of $H$ obtained above,
and subsequently we have used the method described
in \cite{LuxMueRin} to compute the whole submodule lattice
of $P(18a)$, from which the result follows easily.
%As a by-product, we have used the extension described
%in \cite{LuxWie} to cross-check its radical and socle series.
\qquad$\blacksquare$

\bigskip\noindent
{\bf 3.7.Notation.}
We use the notation 
$\chi_{9a}, \chi_{9b}, \chi_{9c}, \chi_{9d}$,
$\chi_{18a}$, $\chi_{18b}$, $\chi_{18c}$, $\chi_{72a}$, $\chi_{72b}$,
$9a$, $9b$, $9c$, $9d$, $18a$, $18b$, $18c$,
and also the source idempotent $j$
as in {\rm{\bf 3.6}}.

\bigskip\noindent
{\bf 3.8.Lemma.}
{\it
The block algebra $B$ and its source algebra
$k[P \rtimes SD_{16}]$ have exactly 18
trivial source modules.
In fact, it holds the following:
\begin{enumerate}
  \renewcommand{\labelenumi}{(\roman{enumi})}
    \item 
Seven PIMs: $P(9a)$, $P(9b)$, $P(9c)$, $P(9d)$, 
             $P(18a)$, $P(18b)$, $P(18c)$.
    \item 
Seven trivial source modules with a vertex $P$ : 
$9a$, $9b$, $9c$, $9d$, $18a$, $18b$, $18c$.
%
%

%\pagebreak
%
   \item
Four trivial source modules with vertex $Q \cong C_3$: 

%JM: introduced boxes

\bigskip\bigskip

\begin{center}
$\chi_{9a}+\chi_{9b}+\chi_{18a}+\chi_{72a} 
\quad
\leftrightarrow 
\quad 
V_1 \ = \  
\begin{matrix}\boxed{
\begin{picture}(90,46)(0,0)
\put(40,40){$18a$}
\put(20,20){$18b$}
\put(40,0){$18a$}
\put(60,20){$18c$}
\put(80,40){$9b$}
\put(0,40){$9a$}
\put(0,0){$9b$}
\put(80,0){$9a$}
\put(10,8){\line(1,1){10}}
\put(30,18){\line(1,-1){10}}
\put(50,8){\line(1,1){10}}
\put(70,18){\line(1,-1){10}}
\put(10,38){\line(1,-1){10}}
\put(30,28){\line(1,1){10}}
\put(50,38){\line(1,-1){10}}
\put(70,28){\line(1,1){10}}
\end{picture}}
\end{matrix} $

\end{center}

\bigskip\bigskip

\begin{center}
$\chi_{9c}+\chi_{9d}+\chi_{18a}+\chi_{72b}
\quad 
\leftrightarrow 
\quad V_2 \ = \ 
\begin{matrix}\boxed{
\begin{picture}(90,46)(0,0)
\put(40,40){$18a$}
\put(20,20){$18c$}
\put(40,0){$18a$}
\put(60,20){$18b$}
\put(80,40){$9d$}
\put(0,40){$9c$}
\put(0,0){$9d$}
\put(80,0){$9c$}
\put(10,8){\line(1,1){10}}
\put(30,18){\line(1,-1){10}}
\put(50,8){\line(1,1){10}}
\put(70,18){\line(1,-1){10}}
\put(10,38){\line(1,-1){10}}
\put(30,28){\line(1,1){10}}
\put(50,38){\line(1,-1){10}}
\put(70,28){\line(1,1){10}}
\end{picture}}\end{matrix}$

\end{center}

\bigskip\bigskip

\begin{center}
{
$\chi_{18b}+\chi_{18c}+\chi_{72a}
\quad \leftrightarrow \quad
V_3 \ = \  
\begin{matrix}\boxed{
\begin{picture}(90,46)(0,0)
\put(60,40){$18c$}
\put(20,40){$18b$}
\put(80,20){$9a$}
\put(40,20){$18a$}
\put(0,20){$9b$}
\put(60,0){$18b$}
\put(20,0){$18c$}
\put(10,18){\line(1,-1){10}}
\put(50,18){\line(1,-1){10}}
\put(30,8){\line(1,1){10}}
\put(70,8){\line(1,1){10}}
\put(10,28){\line(1,1){10}}
\put(50,28){\line(1,1){10}}
\put(30,38){\line(1,-1){10}}
\put(70,38){\line(1,-1){10}}
\end{picture}} \end{matrix}$
}\end{center}

%
%
%
%\begin{matrix}
%V_3 =
%&\boxed{
%  \begin{matrix}
%  18b \ \ 18c  \\
%  9b \ \ 18a \ \ 9a \\
%  18c \ \ 18b
%  \end{matrix}
%  }
%            \\
% &  \updownarrow \\
% & \chi_{18b}+\chi_{18c}+\chi_{72a}
%\end{matrix},
%
%\qquad
%

\bigskip\bigskip

\begin{center}
{
$\chi_{18b}+\chi_{18c}+\chi_{72b}
\quad \leftrightarrow \quad V_4 \ = \ 
\begin{matrix}\boxed{
\begin{picture}(90,46)(0,0)
\put(60,40){$18c$}
\put(20,40){$18b$}
\put(80,20){$9d$}
\put(40,20){$18a$}
\put(0,20){$9c$}
\put(60,0){$18b$}
\put(20,0){$18c$}
\put(10,18){\line(1,-1){10}}
\put(50,18){\line(1,-1){10}}
\put(30,8){\line(1,1){10}}
\put(70,8){\line(1,1){10}}
\put(10,28){\line(1,1){10}}
\put(50,28){\line(1,1){10}}
\put(30,38){\line(1,-1){10}}
\put(70,38){\line(1,-1){10}}
\end{picture}} \end{matrix} $
}
\end{center}

\bigskip\bigskip\bigskip

%$$
%\begin{matrix}
%V_4 =
%&\boxed{
%  \begin{matrix}
%  18c \ \ 18b  \\
%  9d \ \ 18a \ \ 9c \\
%  18b \ \ 18c
%  \end{matrix}
%  }
%            \\
% &  \updownarrow \\
% & \chi_{18b}+\chi_{18c}+\chi_{72b}
%\end{matrix},
%$$
%
%
%
\end{enumerate}
and all characters $\chi_{V_i}$ realized by $V_i$
has values
$\chi_{V_i} (u) = 27$
for any $u \in 3A$,
where $3A$ is the unique conjugacy
class of $H$ of elements of order $3$ on which
$\chi_{V_i}$ does not vanish,
see \cite[Chapter~5, Corollary 1.10(i)]{NagaoTsushima}.
}

\bigskip\noindent
{\bf Proof.}
These follow from {\bf 3.4}, a theorem of Green
\cite[Chapter 4, Problem 10, p.302]{NagaoTsushima}
As for (iii), starting again with the permutation representation
of $H$, using {\sf GAP} \cite{GAP} we compute $N_H(Q)$,
use the {\sf MeatAxe} \cite{MA} and the methods described 
in \cite{LuxSzokeII} to find the PIMs of $N_H(Q)/Q$
as direct summands of its regular representation,
induce them to $H$, and find the submodule structure 
of the induced modules using the methods described in \cite{LuxMueRin}.
\qquad$\blacksquare$

\bigskip\noindent
{\bf 3.9.Notation.}
We use the notation $V_1$, $V_2$, $V_3$, $V_4$
as in {\rm{\bf 3.8}}.

\bigskip\noindent
{\bf 3.10.Lemma.}
{\it
There are no $kH$-modules in $B$
whose radical and socle series are the same and
which have the following structure:
%JM: layers are sufficient, diagram not needed
$$
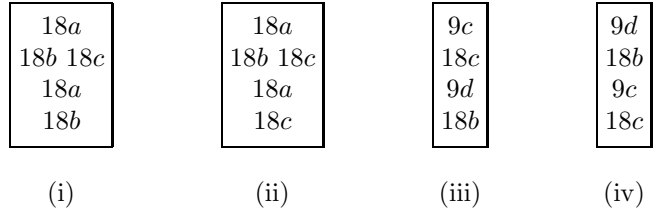
\begin{matrix}
   \boxed{
  \begin{matrix}
     18a  \\
     18b \ 18c \\
     18a  \\
     18b
  \end{matrix}
        }
\\
\\
    {\mathrm{(i)}}
\end{matrix}
\qquad\qquad
 \begin{matrix}
  \boxed{
    \begin{matrix} 
     18a  \\
     18b \ 18c \\
     18a  \\
     18c
    \end{matrix}     
 }
\\
\\
   {\mathrm{(ii)}}
\end{matrix}
\qquad\qquad
\begin{matrix}
\boxed{
   \begin{matrix}
      9c \\
     18c \\
      9d \\
     18b
   \end{matrix}      
      }
\\
\\
    {\mathrm{(iii)}}
\end{matrix}
\qquad\qquad
   \begin{matrix}
\boxed{
   \begin{matrix}
      9d \\
     18b \\
      9c \\
     18c
   \end{matrix}
      }
\\
\\
   {\mathrm{(iv)}}
  \end{matrix}
$$
}

\bigskip\noindent
{\bf Proof.}
(i) Assume that such a $kH$-module, which we call $M$, exists.
There is an epimorphism $\pi : P(18a) \twoheadrightarrow M$.
Set $K = {\mathrm{Ker}}(\pi)$. 
Then, {\bf 3.6(vi)} and {\bf 1.1} imply that
%JM: added `radical'
$K$ has radical and socle series
$$  \boxed
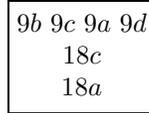
{
 \begin{matrix}     9b \ 9c \ 9a \ 9d \\
                          18c   \\
                          18a
 \end{matrix}
  } 
$$
Since there does not exist a $kH$-module
$\boxed{\begin{matrix} 9a \\
                       18c 
        \end{matrix}
       }
$ by {\bf 3.6(vi)}, we have a contradiction.

(ii) Similar to (i).

(iii) Assume that such a $kH$-module, which we call $M$, exists.
There is an epimorphism $\pi : P(9c) \twoheadrightarrow M$.
Set $K = {\mathrm{Ker}}(\pi)$. 
Then, by {\bf 3.6(vi)} we get
$K = \boxed{\begin{matrix} 18a  \\ 9c
            \end{matrix}   }$.
This contradicts the structure of $P(9c)$ in {\bf 3.6(vi)}.

(iv) Similar to (iii).
\quad$\blacksquare$

\bigskip

%\newpage

\begin{flushleft}
{\bf 4. $3$-Modular representations of {\sf HN}}
\end{flushleft}

%JM: This comment is not strictly needed.
%\bigskip\noindent
%The next result is very important for our purpose
%in this paper. This was relatively recently obtained
%by four people including the second author of this
%paper,
%see also \cite{LuxNoeskeRyba} for the group {\sf HN}
%but for $p = 5$.

\bigskip\noindent
{\bf 4.1.Theorem {\rm (Hiss-M{\"u}ller-Noeske-Thackray
\cite{HissMuellerNoeskeThackray})}.}
{\it
The $3$-decomposition matrix and the Cartan matrix of $A$ are
the following:
}
\bigskip

\begin{center}
\begin{tabular}{r|c|ccccccc}
degree & \cite[p.164--166]{Atlas} & 
$S_1$ & $S_2$ & $S_3$ & $S_4$ & $S_5$ &$S_6$  &$S_7$  \\
\hline
%JM: deleted the commas, because we do not do it elsewise
       8910 & $\chi_{8} $ & 1  &  .    & .     &  .    &  .  & .   & .      \\
      16929 & $\chi_{10}$ & .  &  1    & .     &  .    &  .  & .   & .      \\
     270864 & $\chi_{19}$ & .  &  .    & 1     &  .    &  .  & .   & .       \\
   1185030 & $\chi_{32}$ & 1  &  1    & .     &  1    &  .  & .   & .      \\
   1354320 & $\chi_{33}$ & 1  &  .    & .     &  .    &  1  & 1   & .      \\
   1575936 & $\chi_{37}$ & .  &  .    & 1     &  .    &  .  & 1   & .      \\
   2784375 & $\chi_{43}$ & 1  &  .    & 1     &  1    &  1  & 1   & .      \\
   4561920 & $\chi_{49}$ & .  &  .    & .     &  1    &  1  & .   & 1      \\
   4809375 & $\chi_{50}$ & .  &  1    & 1     &  1    &  .  & .   & 1      \\
\end{tabular} 
\end{center}
\medskip
\begin{center}
\begin{tabular}{r|ccccccc}
      & $P(S_1)$ & $P(S_2)$ & $P(S_3)$ & $P(S_4)$ & $P(S_5)$ & $P(S_6)$ & $P(S_7)$ \\
\hline
$S_1$ & 4 & 1 & 1 & 2 & 2 & 2 & 0   \\ 
$S_2$ & 1 & 3 & 1 & 2 & 0 & 0 & 1   \\ 
$S_3$ & 1 & 1 & 4 & 2 & 1 & 2 & 1   \\ 
$S_4$ & 2 & 2 & 2 & 4 & 2 & 1 & 2   \\ 
$S_5$ & 2 & 0 & 1 & 2 & 3 & 2 & 1   \\
$S_6$ & 2 & 0 & 2 & 1 & 2 & 3 & 0   \\ 
$S_7$ & 0 & 1 & 1 & 2 & 1 & 0 & 2   \\
\end{tabular}
\end{center}

\bigskip\noindent
{\it
where $S_1, \cdots , S_7$ are non-isomorphic simple $kG$-modules
in $A$ whose degrees respectively are
$8910$, $16929$, $270864$, $1159191$, $40338$, $1305072$, $3362391$.
}

%\bigskip\noindent
%{\bf Proof.} See \cite{HissMuellerNoeskeThackray}.
%\quad$\blacksquare$

\bigskip\noindent
{\bf 4.2.Notation.} We use the notation
$\chi_8, \chi_{10}, \chi_{19}, \chi_{32}, \chi_{33}, \chi_{37}, 
 \chi_{43},  \chi_{49}, \chi_{50}$ and
$S_1, \cdots , S_7$ as in {\rm{\bf 4.1}}.

\bigskip\noindent
{\bf 4.3.Lemma.}
{\it
\begin{enumerate}
  \renewcommand{\labelenumi}{(\roman{enumi})}
    \item
All simples $S_1, \cdots, S_7$ are self-dual.
    \item
{\rm{(Kn\"orr)}} \ All simples $S_1, \cdots, S_7$ 
have $P$ as vertices.
\end{enumerate}
}

\bigskip\noindent
{\bf Proof.}
%The group $G$ has exactly three $3$-blocks with defect
%groups of order $3^6$, $3^2$ and $3$, respectively.
%by {\sf GAP} \cite{GAP}, see
%\cite[Lemma 4.2(b)]{AnObrien}.
%Hence we get (i) by {\bf 4.1}.
(i) Easy from {\bf 4.1}. 

(ii) This is a result of Kn{\"o}rr
\cite [3.7.Corollary]{Knoerr}.
\quad$\blacksquare$

\bigskip\noindent
{\bf 4.4.Lemma.}
{\it
\begin{enumerate}
  \renewcommand{\labelenumi}{(\roman{enumi})}
    \item
The heart 
$\mathcal H (P(S_i)) = 
{\mathrm{rad}}(P(S_i))/{\mathrm{soc}}(P(S_i))$
is indecomposable as a $kG$-module for any $i = 1, \cdots , 7$.
    \item
${\mathrm{Ext}}_{kG}^1(S_i, S_j) = 0$
for any pair 
$(i,j) \in \{ (1,2), (1,3), (2,1), (2,2), (2,3), (2,7), (3,1)$, 
 $(3,2), (3,5), (3,7), (4,6), (5,3), (5,5), (5,7), 
 (6,4), (6,6), (7,2), (7,3), (7,5), (7,7) \}$.
\end{enumerate}
}

\bigskip\noindent
{\bf Proof.} (i) This follows by the Cartan matrix of $A$
in {\bf 4.1} and results of Erdmann and Kawata,
see \cite[Theorem 1]{Erdmann}, \cite[Theorem 1.5]{Kawata}
and \cite[1.9.Lemma]{KoshitaniKunugiWaki2002}.

(ii) If ${\mathrm{Ext}}_{kG}^1(S_1, S_2) \not= 0$,
then $S_2 | \mathcal H (P(S_1))$ since $c_{12} = 1$
by {\bf 4.1}, which contradicts to {\bf (i)}.
Similar for the others.
\quad$\blacksquare$

%\pagebreak
\bigskip\noindent
{\bf 4.5.Lemma.}
{\it
%The simple modules $S_1$, $S_2$, $S_3$ are all
%trivial source modules corresponding to
%$S_1 \leftrightarrow \chi_8$,
%$S_2 \leftrightarrow \chi_{10}$,
%$S_3 \leftrightarrow \chi_{19}$.
%}
\begin{enumerate}
  \renewcommand{\labelenumi}{(\roman{enumi})}
    \item
The simple $S_1$ is a trivial source module 
with $S_1 \leftrightarrow \chi_8$.
   \item
The simple $S_2$ is a trivial source module 
with $S_2 \leftrightarrow \chi_{10}$.
   \item
The simple $S_3$ is a trivial source module 
with $S_3 \leftrightarrow \chi_{19}$.
\end{enumerate}
}

\bigskip\noindent
{\bf Proof.}
(i)--(ii)
%First we prove (i) and (ii) together.
Let $M = 2.{\mathsf{HS}}.2$, where ${\mathsf{HS}}$
is the Higman-Sims simple group, be the second largest 
maximal subgroups of $G$, see by \cite[p.164--166]{Atlas}.
%JM: Slightly more detailed on the computations
Now, a calculation with {\sf GAP} \cite{GAP}, using the
character tables of $M$ and $G$, shows that
$1_M{\uparrow}^G{\cdot}1_A = \chi_8 + \chi_{10}$.
Set $X = {k_M}{\uparrow}{^G}{\cdot}{1_A}$. We then get
$X = S_1 + S_2$ (as composition factors) by {\bf 4.1}.
Since $X$, $S_1$, $S_2$ are all self-dual by 
{\bf 4.3(i)}, we obtain $X = S_1 \oplus S_2$.

(iii) 
Let $M$ be the same as above. 
%JM: That's obvious
%It follows from 
%\cite[HS.2]{ModularAtlasProject}
%or \cite{ModularAtlas} that
There uniquely exists a non-trivial linear character
$\chi$ of $M$.
%There is a block
%algebra $\beta$ of $\mathcal OM$, which is called
%"Block~2" in \cite[HS.2]{ModularAtlasProject},
%see \cite{ModularAtlas}.
%There is a character $\chi \in {\mathrm{Irr}}(\beta)$
%with $\chi(1) = 1$. 
Then, a calculation with {\sf GAP} \cite{GAP} shows that
$\chi{\uparrow}{_M^G}{\cdot}{1_A} = \chi_{19}$.
Hence, by {\bf 4.1}, $S_3$ is a trivial source module.
\quad$\blacksquare$

\bigskip\noindent
{\bf 4.6.Lemma.}
{\it
There is a trivial source $kG$-module in $A$ whose vertex is $Q$ 
and whose structure is} 

\begin{center}
$\begin{matrix}\boxed{
\begin{picture}(70,46)(0,0)
\put(30,40){$S_4$}
\put(30,0){$S_4$}
\put(60,20){$S_7$}
\put(40,20){$S_5$}
\put(20,20){$S_2$}
\put(0,20){$S_1$}
\put(10,18){\line(2,-1){20}}
\put(40,8){\line(2,1){20}}
\put(10,28){\line(2,1){20}}
\put(40,38){\line(2,-1){20}}
\put(25,18){\line(2,-3){7}}
\put(45,18){\line(-2,-3){7}}
\put(25,28){\line(2,3){7}}
\put(45,28){\line(-2,3){7}}
\end{picture}} \end{matrix}$
\quad
$\leftrightarrow \quad\chi_{32} + \chi_{49}$.
\end{center}

\bigskip

%and whose radical and socle series is
%$$
%\boxed
%{
%\begin{matrix}
%      S_4  \\
%S_1 \ \ S_2 \ \ S_5 \ \ S_7  \\
%      S_4
%\end{matrix}
%}
%\ \ 
%\leftrightarrow
%\ \ 
%\chi_{32} + \chi_{49}.
%$$
%}

\bigskip\noindent
{\bf Proof.}
It follows from \cite[p.164--166]{Atlas} that the fourth largest
maximal subgroup of $G$ is of the form
$M = 2_+^{1+8}.(\mathcal A_5 \times \mathcal A_5).2$.
%Let $\tilde M$ be a normal subgroup of $M$ of index $2$,
Let $P_M \in {\mathrm{Syl}}_3(M)$. Then
$P_M \cong C_3 \times C_3$, but a calculation with {\sf GAP} \cite{GAP},
using the character tables of $G$ and $M$, shows that $P_M$ 
contains elements belonging to the $3B$ conjugacy class of $G$,
hence ${P_M} \, {\not=_G} \, {P}$ by {\bf 3.4(v)}.
Clearly, there is a non-trivial $kM$-module $T$ with 
$\dim_k T = 1$. Set
$X = T{\uparrow}_M^G{\cdot}1_A$.
Then, $X$ is a direct sum of trivial source $kG$-modules, 
and a calculation with {\sf GAP} \cite{GAP} shows that
$X \leftrightarrow \chi_{32} + \chi_{49}$.
Since $P$ is a defect group of $A$, any indecomposable
$kG$-module $Y$ with $Y|X$ does not have $P$ as its vertex.

Suppose that $X$ is decomposable. Then, {\bf 2.3(i)} implies
that $X = Y \oplus Z$ such that
$Y \leftrightarrow \chi_{32}$ and
$Z \leftrightarrow \chi_{49}$.
Hence, {\bf 4.1} yields that $Y = S_1 + S_2 + S_4$
(as composition factors). We know by {\bf 4.3(ii)} and
{\bf 4.1} that $S_1$, $S_2$, and $Y$ are all self-dual.
If $[Y, S_1]^G \not= 0$, then the self-dualities imply
$S_1 | Y$, and hence
$0 \not= [S_1, Y]^G = ( \chi_8, \chi_{32})^G$ from
{\bf 2.3(ii)} and {\bf 4.5(i)}, a contradiction.
Hence, $[Y, S_1]^G = [S_1, Y]^G = 0$. Similarly, we obtain
$[Y, S_2]^G = [S_2, Y]^G = 0$ by {\bf 2.3(ii)} and {\bf 4.5(ii)}.
This is a contradiction since $Y$ has only three composition 
factors $S_1$, $S_2$ and $S_4$.
  
Thus, $X$ is indecomposable. By the decomposition matix of $A$
in {\bf 4.1}, $X$ is not a PIM. Thus, the order of a vertex of $X$
is $3$, and hence $Q$ is a vertex of $X$ by {\bf 3.4(vi)}.
Clearly, $X$ is a trivial source $kG$-module in $A$.
We know by {\bf 4.1} that
$X = S_1 + S_2 + 2 \times S_4 + S_5 + S_7$ 
(as composition factors). Note that
$X$, $S_1$, $S_2$, $S_4$, $S_5$, $S_7$ are all self-dual
from {\bf 4.3(i)}. Then,
$[S_i, X]^G = [X, S_i]^G = 0$ for any $i = 1,2,5,7$
since $X$ is indecomposable. Thus,
$X/{\mathrm{rad}}(X) \cong {\mathrm{Soc}}(X) \cong S_4$.
Therefore, again by the self-dualities, it holds that
${\mathrm{rad}}(X) / {\mathrm{Soc}}(X)  \cong
S_1 \oplus S_2 \oplus S_5 \oplus S_7$.
\quad$\blacksquare$

\bigskip\noindent
{\bf 4.7.Lemma.}
{\it
There is a trivial source $kG$-module in $A$ which has $Q$
as a vertex and has radical and socle series %the form
$\boxed{\begin{matrix} S_3 \\ S_6 \\ S_3 \end{matrix}}
%U(S_3, S_6, S_3)
\
\leftrightarrow
\  
\chi_{19} + \chi_{37}.
$

%JM: It is safer to do it like this ...
\bigskip
\bigskip\noindent
{\rm{(}}{\bf Note:} We can prove that this module has $Q$
as its vertex, but only later on in {\bf 7.2(ii)}{\rm{)}}.
}

\bigskip\noindent
{\bf Proof.}
First, the third largest maximal subgroup of $G$ is of shape
$M = U_3(8).3$, see \cite[p.164--166]{Atlas}.
Then, a calculation with {\sf GAP} \cite{GAP},
using the character tables of $G$ and $M$, shows that
\begin{align}
  {1_M}{\uparrow}^G{\cdot}1_A = 
  \chi_{19} + \chi_{32} + \chi_{37} + \chi_{49}.
\end{align}
Set $X = k_M{\uparrow}^G{\cdot}1_A$, hence $X$ is
self-dual and is a direct sum of trivial source
$kG$-modules. Then, by the decomposition matrix
in {\bf 4.1}, we know
\begin{align}
 X = S_1 + S_2 + 2 \times S_3 + 2 \times S_4 + S_5 + S_6 + S_7
\qquad
({\text{as \ composition \ factors}}).
\end{align}
If $[X, S_1]^G \not= 0$, then
{\bf 2.3(ii)} and {\bf 4.5(i)} imply that
$(\chi_{\hat X}, \chi_8)^G = [X, S_1]^G \not= 0 $
where $\chi_{\hat X}$ is a  character afforded
by $X$ (see {\bf 2.3(i)}), 
which is a contradiction by (1).
Hence, it holds $[X, S_1]^G = [S_1, X]^G = 0$ by the self-dualities
in {\bf 4.3(i)}. Similarly, by {\bf 4.5(ii)-(iii)} and
{\bf 2.3(ii)}, we know also
$[X, S_2]^G = [S_2, X]^G = 0$ and 
$[X, S_3]^G = [S_3, X]^G = 1$.
If $[X, S_5]^G \not= 0$, then (2) and the self-dualities
imply that $S_5 | X$, and hence $S_5$ is liftable
by {\bf 2.3(i)}, which contradicts {\bf 4.1}.
Hence, $[X, S_5]^G = [S_5, X]^G = 0$ by the self-dualities.
Similarly, it holds also that
$[X, S_i]^G = [S_i, X]^G = 0$ for $i = 6, 7$.
If $[X, S_4]^G = 2$, then it follows from (2) and the
self-dualities that
$(S_4 \oplus S_4) | X$, and hence $S_4$ is liftable
by {\bf 2.3(i)}, which contradicts {\bf 4.1}.
This shows
$[X, S_4]^G = [S_4, X]^G \not= 2$.
Namely, 
\begin{align}
  [S_3, X]^G = [X, S_3]^G = 1,
\\
  [S_4, X]^G =  [X, S_4]^G \not= 2, 
\\ 
  [S_i, X]^G = [X, S_i]^G = 0
   \quad {\text{for}} \ i = 1, 2, 5, 6, 7.
\end{align}

Now, {\bf 4.6} says that there is a trivial source
$kG$-module $Y$ that has radical and socle series 
\begin{align}
Y \ = \ 
\begin{matrix}
\boxed{
\begin{picture}(70,46)(0,0)
\put(30,40){$S_4$}
\put(30,0){$S_4$}
\put(60,20){$S_7$}
\put(40,20){$S_5$}
\put(20,20){$S_2$}
\put(0,20){$S_1$}
\put(10,18){\line(2,-1){20}}
\put(40,8){\line(2,1){20}}
\put(10,28){\line(2,1){20}}
\put(40,38){\line(2,-1){20}}
\put(25,18){\line(2,-3){7}}
\put(45,18){\line(-2,-3){7}}
\put(25,28){\line(2,3){7}}
\put(45,28){\line(-2,3){7}}
\end{picture}}
%      S_4  \\
%S_1 \ \ S_2 \ \ S_5 \ \ S_7  \\
%      S_4
\end{matrix}
\ \ 
\leftrightarrow
\ \ 
\chi_{32} + \chi_{49}.
\end{align} 
in $A$. Then, by (1), (6) and {\bf 2.3(ii)}, we have
\begin{align}
 [Y, X]^G \ = \ [X, Y]^G \ = \ 2.
\end{align}
Then, by (4) and (2), we know that
\begin{align}
  [S_4, X]^G =  [X, S_4]^G \leqslant 1.
\end{align}

Next, we want to claim that there is a homomorphism
$\varphi \in {\mathrm{Hom}}_{kG}(Y, X)$ 
with $0 \not= {\mathrm{Im}}(\varphi) \not\cong  S_4$.
Suppose that any non-zero 
$\varphi \in {\mathrm{Hom}}_{kG}(Y, X)$
satisfies that ${\mathrm{Im}}(\varphi) \cong S_4$.
By (7), let $\{ \varphi_1, \varphi_2 \}$ be a $k$-basis
of ${\mathrm{Hom}}_{kG}(Y, X)$.
Then, it follows from Schur's lemma that
${\mathrm{Im}}(\varphi_1) \not= {\mathrm{Im}}(\varphi_2)$,
and hence that there exists a direct sum
${\mathrm{Im}}(\varphi_1) 
 \oplus {\mathrm{Im}}(\varphi_2) \subseteq X.$
This means that $[S_4, X]^G \geqslant 2$,
contradicting (8).

Therefore, there is a homomorphism
$\varphi \in {\mathrm{Hom}}_{kG}(Y, X)$ 
with $0 \not= {\mathrm{Im}}(\varphi)  \not\cong  S_4$.
Then, by (6), we know ${\mathrm{Ker}}(\varphi) = 0$
since $S_i \, {\not|} \, {\mathrm{soc}}(X)$
for $i = 1, 2, 5, 7$ by (5).
That is, there is a monomorphism
$\varphi : Y \rightarrowtail X$ of $kG$-modules.

Then, just by the dual argument,
we know also that there is an
epimorphism $\psi : X \twoheadrightarrow Y$
of $kG$-modules.
It follows then by (2) and (6) that there is a direct sum
${\mathrm{Im}}(\varphi) \oplus {\mathrm{Ker}}(\psi)
 \subseteq X$,
and hence 
${\mathrm{Im}}(\varphi) \oplus {\mathrm{Ker}}(\psi) = X$.
Set $Z = {\mathrm{Ker}}(\psi)$.
We can write $X = Y \oplus Z$.
Since $Z = 2 \times S_3 + S_6$
(as composition factors), we get by (5) that
$Z = \boxed{\begin{matrix} S_3 \\ S_6 \\ S_3 \end{matrix}}$. 
%U(S_3, S_6, S_3)$.
Hence, it is easy to know from (1) and (6) that
$Z$ is a trivial source $kG$-module with
$Z \leftrightarrow \chi_{19} + \chi_{37}$.
%JM: Actually, I cannot not verify the old GAP computation,
%    probably it was just nonsense.
%Since $Z$ is not projective by {\bf 4.1}, and since
%$Z$ is in $A$, it holds that $P$ or $Q$ is a vertex
%of $Z$ by {\bf 3.4(vi)}. 
%By the definition of $Z$,
%$Z$ is relatively $M$-projective.
%A calculation with {\sf GAP} \cite{GAP} ????? shows that
%$P \not\subseteq_G M$. 
%Thus $P$ is not a vertex of $Z$, but so is $Q$.
\quad$\blacksquare$

\bigskip\noindent
{\bf 4.8.Lemma.}
{\it
There is a trivial source $kG$-module in $A$ 
whose structure is

\bigskip
\begin{center}
$\begin{matrix}\boxed{
\begin{picture}(51,46)(0,0)
\put(0,0){$S_1$}
\put(40,0){$S_2$}
\put(20,20){$S_4$}
\put(0,40){$S_1$}
\put(40,40){$S_2$}
\put(10,8){\line(1,1){10}}
\put(10,38){\line(1,-1){10}}
\put(30,18){\line(1,-1){10}}
\put(30,28){\line(1,1){10}}
\end{picture}}
\end{matrix}$
\quad
$\leftrightarrow\quad
\chi_{8} + \chi_{10} + \chi_{32}$.
\end{center}
%
%
%which has the radical and socle series
%$$
%\boxed{
%\begin{matrix}
%    S_1 \ S_2 \\
%       S_4    \\
%    S_1 \ S_2 
%\end{matrix}
%}
%\
%\leftrightarrow
%\  
%\chi_{8} + \chi_{10} + \chi_{32}.
%$$

\bigskip
\bigskip\noindent
{\rm{(}}{\bf Note:} We can prove that this module has $Q$
as its vertex, but only later on in {\bf 7.2(i)}{\rm{)}}.
}

\bigskip\noindent
{\bf Proof.}
By \cite[p.91]{Atlas}, we have
$1_{\mathcal A_{11}}{\uparrow}^{\mathcal A_{12}}
 = 1_{\mathcal A_{12}} + {\widetilde \chi}_{11}$,
where 
${\widetilde \chi}_{11} \in {\mathrm{Irr}}(\mathcal A_{12})$
is of degree $11$. It follows from the $3$-decomposition matrix
of $\mathcal A_{12}$ in
\cite[$A_{12}$ (mod $3$)]{ModularAtlasProject} and \cite{ModularAtlas} that
\begin{equation}
  k_{\mathcal A_{11}}{\uparrow}^{\mathcal A_{12}} \ = \
  \boxed{\begin{matrix} k \\ 10 \\ k \end{matrix}} 
  %U(k, 10, k) 
  \quad \leftrightarrow \quad
  1_{\mathcal A_{12}} + {\widetilde \chi}_{11},
\end{equation}
where $10$ is a simple $k{\mathcal A_{12}}$-module
of dimension $10$. 
Set 
$X = (k_{\mathcal A_{11}}{\uparrow}^{\mathcal A_{12}})
     {\uparrow}^G{\cdot}1_A$.
Note that $X$ is a direct sum of trivial source
$kG$-modules. Then, we know from a calculation with
{\sf GAP} \cite{GAP},
using the character tables of $G$ and $\mathcal A_{12}$, that
\begin{equation}
  X \ \leftrightarrow \ 
 \chi_{8} + \chi_{10} + \chi_{32}
\end{equation}
and
\begin{equation}
  X = 2 \times S_1 + 2 \times S_2 + S_4
  \qquad
  {\text{(as \ composition \ factors)}}.
\end{equation}
By (7), {\bf 2.3(ii)} and {\bf 4.5(ii)}, we obtain
$[X, S_1]^G = (\chi_{\hat X}, \chi_8)^G =  1$.
Hence,
$[X, S_1]^G = [S_1, X]^G = 1$ by the self-dualities.
Similarly, we have 
$[X, S_2]^G = [S_2, X]^G = 1$. 
Since $S_4$ is not liftable by {\bf 4.1}, $S_4$ is not
a trivial source module by {\bf 2.3(i)}. This implies that
$[X, S_4]^G = [S_4, X]^G = 0$ by (8).
These yield
\begin{equation}
 X/{\mathrm{rad}}(X) \ \cong \ {\mathrm{Soc}}(X) 
 \ \cong S_1 \oplus S_2.
\end{equation}

Next, we want to claim that $X$ is indecomposable.
Suppose that $X$ is decomposable.
By (12), we can write $X = X_1 \oplus X_2$ for 
$A$-submodules $X_1$ and $X_2$ of $X$ with
${\mathrm{soc}}(X_i) \cong S_i$ for $i = 1, 2$.
If $X_1/{\mathrm{rad}}(X_1) \not\cong S_1$, then
%JM: replaced 11 by 12
(12) shows that $X_1/{\mathrm{rad}}(X_1) \cong S_2$,
and hence we get by (12) and (11) that
$X = X_1 \oplus X_2 = \boxed{\begin{matrix} S_2 \\ S_4 \\ S_1 \end{matrix}}
%U(S_2, S_4, S_1) 
\oplus \boxed{\begin{matrix} S_1 \\ S_2 \end{matrix}}$ %U(S_1, S_2)$
or
$X = X_1 \oplus X_2 = \boxed{\begin{matrix} S_2 \\ S_1 \end{matrix}} 
%U(S_2, S_1) 
\oplus \boxed{\begin{matrix} S_1 \\ S_4 \\ S_2 \end{matrix}}$,
% U(S_1, S_4, S_2)$,
which is a contradiction by 
the self-dualities of $X$ and each $S_i$ in {\bf 4.4(i)}.
This means that $X_i/{\mathrm{rad}}(X_i) \cong S_i$
for $i = 1,2$ by (12).
If $X_1$ is simple, then we get by (12) that
%JM: replaced `Loewy' by `radical' 
$X_2$ has radical and socle series
which is one of the following three cases:
%JM: introduced boxes 
$$
\boxed
{\begin{matrix}   S_2 \\ S_1 \ S_4 \\ S_2 \end{matrix}}
\qquad
\boxed
{\begin{matrix}   S_2 \\ S_1 \\ S_4 \\ S_2 \end{matrix}}
\qquad
\boxed
{\begin{matrix}   S_2 \\ S_4 \\ S_1 \\ S_2 \end{matrix}}
$$
So we have a contradiction by {\bf 4.4(ii)}.
Thus, $X_1$ is not simple. Similarly, we know that
$X_2$ is not simple. Hence, {\bf 2.16} yields that
${\mathrm{soc}}(X_i) \subseteq {\mathrm{rad}}(X_i)$
for $i = 1, 2$.
Thus,
$X = X_1 \oplus X_2 = \boxed{\begin{matrix} S_1 \\ S_4 \\ S_1 \end{matrix}}
%U(S_1, S_4, S_1) 
\oplus \boxed{\begin{matrix} S_2 \\ S_2 \end{matrix}}$ %U(S_2, S_2)$
or
$X = X_1 \oplus X_2 = \boxed{\begin{matrix} S_1 \\ S_1 \end{matrix}} 
%U(S_1, S_1) 
\oplus \boxed{\begin{matrix} S_2 \\ S_4 \\ S_2 \end{matrix}}$. 
% U(S_2, S_4, S_2)$.
This is a contradiction by (10), {\bf 2.3(i)} and {\bf 4.1}.

Therefore $X$ is indecomposable. Hence, we get
by (11), (12) and {\bf 2.16} that 
${\mathrm{soc}}(X) \subseteq {\mathrm{rad}}(X)$.
Thus we get the structure of $X$ as desired.
\qquad$\blacksquare$

\bigskip\noindent
{\bf 4.9.Notation.}
In the rest of paper let $f$ be the Green correspondence
from $G$ to $H$ with respect to $P$, 
see \cite[Chapter 4 \S 4]{NagaoTsushima}.

\bigskip\noindent
{\bf 4.10.Lemma.}
{\it
It holds that $f(S_1) = 9a$.
}

%JM: This does not follow uniquely from characters!
\bigskip\noindent
{\bf Proof.}
It follows from {\bf 4.5(i)}, {\bf 4.3(ii)} and {\bf 2.1} 
that $f(S_1)$ is a simple $kH$-module in $B$,
see {\bf 3.4(i)}.
Using the ordinary characters afforded by the %information on 
trivial source $kH$-modules in $B$, see {\bf 3.8},
we get the following possible decompositions of
$S_1{\downarrow}_H{\cdot}1_B$, by
a calculation with {\sf GAP} \cite{GAP}
using the character tables of $G$ and $H$:
%In fact, we get the decomposition of
%$S_1{\downarrow}_H{\cdot}1_B$ by using the character tables
%of characters in $A$ and $B$ depending on {\sf GAP} \cite{GAP}.
%Actually, we can check it by hand as well by using the
%information on trivial source $kH$-modules in $B$, see {\bf 3.8}.
$$
S_1{\downarrow}_H{\cdot}1_B \ = \
9a \ \bigoplus  \ \Big(
7 \times P(9a) \oplus 7 \times P(9b) \oplus 
5 \times P(18a) \oplus P(18b) \oplus P(18c) \Big) 
$$
or
$$
S_1{\downarrow}_H{\cdot}1_B \ = \
9b \ \bigoplus  \ \Big(
8 \times P(9a) \oplus 6 \times P(9b) \oplus 
5 \times P(18a) \oplus P(18b) \oplus P(18c) \Big).
$$
In particular, $f(S_1) = 9a$ or $f(S_1) = 9b$,
and we have to decide which case actually occurs.

To this end, let $M=2.{\sf HS}.2$ be the second largest maximal 
subgroup of $G$, see {\bf 4.5}. By \cite[$HS$ (mod 3)]{ModularAtlasProject}
and \cite{ModularAtlas}, let $A^-$ be the block algebra of ${\mathcal O}M$
containing the unique non-trivial linear character $\chi$ of $M$.
Hence letting $A^+$ and $A'$, see {\bf 5.1}, be the principal block algebras 
of ${\mathcal O}M$ and of ${\mathcal O}{\sf HS}$, respectively,
we have $A^+\cong A'$ and an isomorphism 
$-\otimes\chi\colon A^+\rightarrow A^-$.
Moreover, $P$ being a Sylow $3$-subgroup of $M$, it is the 
block defect group of $A^-$, and hence let $B^-$ be the Brauer 
correspondent of $A^-$ in $N_{M}(P)$. 

Using the smallest faithful permutation representation
of $M$ on $1408$ points, available in \cite{ModAtlasRep},
the normaliser $N_{M}(P)$ and the centraliser $C_{M}(P)$ 
of the Sylow $3$-subgroup $P$
is easily computed explicitly with {\sf GAP} \cite{GAP}
and their structure determined,
we find $N_{M}(P)=(P\times D_8).SD_{16}$ and
$C_{M}(P)=P\times D_8$.
Now the conjugacy classes of $N_{M}(P)$ can be computed,
its ordinary character table is found using 
the Dixon-Schneider algorithm, from that its blocks
are determined and $B^-$ is identified.

Then a computation with {\sf GAP} \cite{GAP},
using the character tables of $G$ and $M$, shows that
$S_1{\downarrow}_M{\cdot}1_{A^-}= 22^-$,
where the latter denotes the unique simple $A^-$-module of that dimension.
Moreover, using the character tables of $M$ and $N_{M}(P)$,
{\sf GAP} \cite{GAP} shows that 
$(22^-){\downarrow}_{N_{M}(P)}{\cdot}1_{B^-}=\lambda$,
where $\lambda$ is a certain linear character;
actually, $\lambda$ is the Green correspondent of $22^-$ with 
respect to $(M,P,N_{M}(P))$, which must be linear in view of {\bf 5.7}.
Hence 
$\lambda=(S_1{\downarrow}_M{\cdot}1_{A^-})
 {\downarrow}_{N_{M}(P)}{\cdot}1_{B^-}$
is a direct summand of
\begin{align*}
 (S_1{\downarrow}_H){\downarrow}_{N_{M}(P)}{\cdot}1_{B^-}
&= \ \Big( f(S_1)\oplus({\mathcal R}{\text{-}}{\mathrm{proj}}) \Big)
 {\downarrow}_{N_{M}(P)}{\cdot}1_{B^-} \\
&= \ f(S_1){\downarrow}_{N_{M}(P)}{\cdot}1_{B^-} 
 \oplus (Q{\text{-}}{\mathrm{proj}}) \oplus ({\mathrm{proj}}) ,
\end{align*}
where ${\mathcal R}$ consists of elementary-abelian $3$-subgroups of $H$,
of order at most $9$, not $H$-conjugate to $P$, 
and $Q\cong C_3$ is as in {\bf 3.3}.
Since $\lambda$ has $P$ as a vertex, we conclude
that $\lambda$ is a direct summand of 
$f(S_1){\downarrow}_{N_{M}(P)}{\cdot}1_{B^-}$.

Now a computation with {\sf GAP} \cite{GAP},
using the character tables of $H$ and $N_{M}(P)$, shows that
$f(S_1){\downarrow}_{N_{M}(P)}{\cdot}1_{B^-}
=(9x){\downarrow}_{N_{M}(P)}{\cdot}1_{B^-}$, where $x\in\{a,b\}$,
already is linear, where 
$$ (9a){\downarrow}_{N_{M}(P)}{\cdot}1_{B^-}=\lambda
   \not=(9b){\downarrow}_{N_{M}(P)}{\cdot}1_{B^-} .$$
This shows that $f(S_1)=9a$.
\quad$\blacksquare$

\bigskip\noindent
{\bf 4.11.Lemma.}
{\it
It holds that $f(S_2) = 9b$. 
}

\bigskip\noindent
{\bf Proof.}
It follows from {\bf 4.5(i)}, {\bf 4.3(ii)} and {\bf 2.1} 
that $f(S_1)$ is a simple $kH$-module in $B$, see {\bf 3.4(i)}.
Using the ordinary characters afforded by the 
trivial source $kH$-modules in $B$, see {\bf 3.8},
we get the following possible decompositions of
$S_1{\downarrow}_H{\cdot}1_B$, by
a calculation with {\sf GAP} \cite{GAP}
using the character tables of $G$ and $H$:
$$
S_2{\downarrow}_H{\cdot}1_B \ = \ 
9b \ \bigoplus  \Big(
9 \times P(9a) \oplus 8 \times P(9b) \oplus 
7 \times P(18a) \oplus 5 \times P(18b) \oplus 5 \times P(18c) \Big)
$$
or
$$
S_2{\downarrow}_H{\cdot}1_B \ = \ 
9a \ \bigoplus  \Big(
8 \times P(9a) \oplus 9 \times P(9b) \oplus 
7 \times P(18a) \oplus 5 \times P(18b) \oplus 5 \times P(18c) \Big).
$$
In particular, $f(S_2) = 9b$ or $f(S_2) = 9a$,
hence the assertion follows from {\bf 4.10}.
%Similar to that of {\bf 4.10}.
\quad$\blacksquare$

\bigskip\noindent
{\bf 4.12.Lemma.}
{\it
It holds that $f(S_3) = 9c$. 
}

\bigskip\noindent
{\bf Proof.}
It follows from {\bf 4.5(i)}, {\bf 4.3(ii)} and {\bf 2.1}
that $f(S_1)$ is a simple $kH$-module in $B$, see {\bf 3.4(i)}.
Using the ordinary characters afforded by the
trivial source $kH$-modules in $B$, see {\bf 3.8},
we get the following possible decompositions of
$S_1{\downarrow}_H{\cdot}1_B$, by
a calculation with {\sf GAP} \cite{GAP}
using the character tables of $G$ and $H$:
\begin{align*}
S_3{\downarrow}_H{\cdot}1_B \ = \
9c \ \bigoplus  \Big(
54 \times P(9a) &\oplus 54 \times P(9b) \oplus 40 \times P(9c) 
                 \oplus 41 \times P(9d)                       \\
                &\oplus 94 \times P(18a) \oplus 95 \times P(18b) 
                 \oplus 95 \times P(18c) \Big)
   \ \bigoplus \ V_3
\end{align*}
or
\begin{align*}
S_3{\downarrow}_H{\cdot}1_B \ = \
9c \ \bigoplus  \Big(
54 \times P(9a) &\oplus 54 \times P(9b) \oplus 39 \times P(9c) 
                 \oplus 40 \times P(9d)                       \\
                &\oplus 93 \times P(18a) \oplus 96 \times P(18b) 
                 \oplus 96 \times P(18c) \Big)
   \ \bigoplus \ V_2
\end{align*}
or
\begin{align*}
S_3{\downarrow}_H{\cdot}1_B \ = \
9d \ \bigoplus  \Big(
54 \times P(9a) &\oplus 54 \times P(9b) \oplus 41 \times P(9c) 
                 \oplus 40 \times P(9d)                       \\
                &\oplus 94 \times P(18a) \oplus 95 \times P(18b) 
                 \oplus 95 \times P(18c) \Big)
   \ \bigoplus \ V_3
\end{align*}
or
\begin{align*}
S_3{\downarrow}_H{\cdot}1_B \ = \
9d \ \bigoplus  \Big(
54 \times P(9a) &\oplus 54 \times P(9b) \oplus 40 \times P(9c) 
                 \oplus 39 \times P(9d)                       \\
                &\oplus 93 \times P(18a) \oplus 96 \times P(18b) 
                 \oplus 96 \times P(18c) \Big)
   \ \bigoplus \ V_2,
\end{align*}
where $V_3$ and $V_2$ are the trivial source $kH$-modules in $B$ 
with vertex $Q$ given in {\bf 3.8}.
In particular, $f(S_3) = 9c$ or $f(S_3) = 9d$,
and we have to decide which case actually occurs.

Keeping the notation from {\bf 4.10}, we by the proof of {\bf 4.5(iii)}
have $S_3=\chi{\uparrow}^G{\cdot}1_A$, hence
$(\chi{\uparrow}^G{\cdot}1_A){\downarrow}_H{\cdot}1_B 
=S_3{\downarrow}_H{\cdot}1_B
=f(S_3)\oplus(Q{\text{-}}{\mathrm{proj}})\oplus({\mathrm{proj}})$.
Hence $f(S_3)$ is a direct summand of
$$
(\chi{\uparrow}^G){\downarrow}_H{\cdot}1_B
=\bigoplus_{g} 
 \Big((\chi^g){\downarrow}_{M^g \cap H}\Big){\uparrow}^H{\cdot}1_B ,$$
where $g$ runs through a set of representatives of
the $M$-$H$ double cosets in $G$.
Since $f(S_3)$ has $P$ as a vertex, and $P$ is normal in $H$,
we only have to look at summands coming from $g\in G$
such that $P\leq M^g \cap H$. But for these $g$ we have
$P,P^{g^{-1}}\leq M$, which since $P$ is a Sylow $3$-subgroup of $M$
implies the existence of $m\in M$ such that $P^m=P^{g^{-1}}$,
hence $h:=mg\in H$, and thus $g=m^{-1}h\in MH$, that is, we may
assume $g=1$.

Thus we conclude that $f(S_3)$ is a direct summand of 
$(\chi{\downarrow}_{M\cap H}){\uparrow}^H{\cdot}1_B
=(\chi{\downarrow}_{N_M(P)}){\uparrow}^H{\cdot}1_B$.
Now a computation with {\sf GAP} \cite{GAP},
using the character tables of $N_M(P)$ and $H$,
shows that 
$(\chi{\downarrow}_{N_M(P)}){\uparrow}^H{\cdot}1_B=9c$
is indecomposable, showing that $f(S_3)=9c$.

We just remark that it is possible, using {\sf GAP} \cite{GAP}
and specially tailored programs to deal efficiently with permutations
on millions of points, to construct the transitive permutation 
representation of 
$G$ on $3078000$ points, that is, the action of $G$ on the cosets
of $2.\mathsf{HS}$ in $G$, where $2.\mathsf{HS}$ is the derived 
subgroup of $M$, and to use the restriction of this representation 
to $H$ to show that the first of the four possible
decompositions of $S_3{\downarrow}_H{\cdot}1_B$ listed above
actually occurs. But we will not need this fact.
\quad$\blacksquare$

\bigskip
%\newpage

\begin{flushleft}
{\bf 5. Green correspondence for {\sf HS}} 
\end{flushleft}

\bigskip\noindent
{\bf 5.1.Notation and assumption.}
In the rest of this paper, we use the following notation, too.
Let $G'$ be the Higman-Sims simple group {\sf HS}.
Since Sylow $3$-subgroups of $G'$ are 
isomorphic to $C_3 \times C_3$,
%JM: I think, we do not need an explicit identification of 
%    the block defect groups of A and A', although of course
%    P < 2.HS < G, and the principal blocks of HS and 2.HS coincide.
%we can assume that $P$ is a Sylow $3$-subgroup of {\sf HS} as well,
%which is originally a defect group of the block algebra $A$ of
%$kG$ where $G = {\sf HN}$, see {\bf 3.3}. 
we by abuse of notation let $P$ denote a
Sylow $3$-subgroup of {\sf HS} as well.
There is exactly one conjugacy class of $G'$ which contain
elements of order $3$, that is, $P$ has exactly one
$G'$-conjugacy class of subgroups of order $3$,
see \cite[p.81]{Atlas}.
%JM: changed C_2 to 2
Let $H' = N_{G'}(P)$, and hence $H' = (P \rtimes SD_{16}) \times 2$,
where the action of $SD_{16}$ on $P$ is given by the embedding
of $SD_{16}$ as a Sylow $2$-subgroup of $\mathrm{Aut}(P)\cong\mathrm{GL}_2(3)$.
Let $A'$ and $B'$, respectively, be the principal block algebras
of $\mathcal OG'$ and $\mathcal OH'$. 
%JM: see La.5.2
%Thus $B' = \mathcal O[P \rtimes SD_{16}]$. 

\bigskip\noindent
{\bf 5.2.Lemma.} {\it
\begin{enumerate}
  \renewcommand{\labelenumi}{(\roman{enumi})}

    \item
The character table of $P \rtimes SD_{16}$
%JM: This is clear anyway.
%\cong (C_3 \times C_3)\rtimes SD_{16}$
is given as follows:
\smallskip

\begin{center}
\begin{tabular}{c|rrrrrrrrr}
{\rm{conjugacy class}} & $1A$& $2A$& $2B$& $3A$& $4A$& $4B$& $6A$& $8A$& $8B$ \\
{\rm{centraliser}} & $144$& $16$& $12$& $18$& $8$& $4$& $6$& $8$& $8$ \\
\hline
$\chi_{1a}$ & $1$& $1$& $1$& $1$& $1$& $1$& $1$& $1$& $1$ \\
$\chi_{1b}$ & $1$& $1$& $1$& $1$& $1$& $-1$& $1$& $-1$& $-1$ \\
$\chi_{1c}$ & $1$& $1$& $-1$& $1$& $1$& $1$& $-1$& $-1$& $-1$ \\
$\chi_{1d}$ & $1$& $1$& $-1$& $1$& $1$& $-1$& $-1$& $1$& $1$ \\
$\chi_{2a}$ & $2$& $2$& $0$& $2$& $-2$& $0$& $0$& $0$& $0$ \\
$\chi_{2b}$ & $2$& $-2$& $0$& $2$& $0$& $0$& $0$& $\sqrt{-2}$& $-\sqrt{-2}$\\
$\chi_{2c}$ & $2$& $-2$& $0$& $2$& $0$& $0$& $0$& $-\sqrt{-2}$& $\sqrt{-2}$\\
$\chi_{8a}$ & $8$& $0$& $2$&  $-1$& $0$& $0$& $-1$& $0$& $0$ \\
$\chi_{8b}$ & $8$& $0$& $-2$& $-1$& $0$& $0$& $1$& $0$& $0$ \\
\end{tabular}
\end{center}
\medskip
Note that this identifies the characters
$\chi_{1a}$, $\chi_{1b}$, $\chi_{1c}$, $\chi_{1d}$,
$\chi_{8a}$, and $\chi_{8b}$ uniquely.

%
%where all $\chi \in {\mathrm{Irr}}(P \rtimes  SD_{16})$
%except $\chi_{2b}$ and $\chi_{2c}$ are self-dual, and
%${\chi_{2b}}^\vee = \chi_{2c}$.
%
%   \item
%There is an automorphism 
%$\sigma \in {\mathrm{Aut}}(SD_{16})$
%of order $2$ which induces 
%$\sigma \in {\mathrm{Aut}}(P \rtimes SD_{16})$
%such that $(8A)^\sigma = 8B$ and $(8B)^\sigma = 8A$ and
%$\sigma$ fixes all the other conjugacy classes
%of $P \rtimes SD_{16}$, and hence 
%${2b}^\sigma = {2c}$ and
%${2c}^\sigma = {2b}$, and
%$\sigma$ fixes all the other simple
%$\mathcal O[P \rtimes SD_{16}]$-modules.
%
%   \item
%With the notation in {\bf (ii)}, $\sigma$ induces an
%$\mathcal O$-algebra-automorphism
%$\tilde\sigma \in {\mathrm{Aut}}(\mathcal O[P \rtimes SD_{16})$,
%and, in addition, $\tilde\sigma$ is an interior 
%$P$-algebra-automorphism.
%
%\end{enumerate}

%
    \item 
$B' \cong {\mathcal O}[P \rtimes SD_{16}]$,
as interior $P$-algebras and hence $k$-algebras,
and we can write that
\begin{align*}
{\mathrm{Irr}}(B') 
&= \{ 1_{H'} = \chi_{1a},\chi_{1b},\chi_{1c}, \chi_{1d},
\chi_{2a}, \chi_{2a},\chi_{2c} = {\chi_{2b}}^\vee,
\chi_{8a}, \chi_{8b}  \}, 
\\
{\mathrm{IBr}}(B') 
&= \{ 1a, 1b, 1c, 1d, 2a, 2b, 2c = 2b^\vee \}, \notag
\end{align*}
where the numbers mean the degrees (dimensions) of characters
(modules).
In particular, all simple modules $1a, 1b, 1c, 1d, 2a$ in $B'$ 
except $2b$ and $2c$ are self-dual.
   \item The $3$-decomposition and the Cartan matrices of $B'$,
respectively, are the following:
\begin{center}
{\rm
\begin{tabular}{l|ccccccc}
  & $1a$ & $1b$ & $1c$ & $1d$ & $2a$ & $2b$ & $2c$    \\
\hline
$\chi_{1a}$  & 1  &  .    & .     &  .    &  .  &.   &. \\
$\chi_{1b}$  & .  &  1    & .     &  .    &  .  &.   &. \\
$\chi_{1c}$  & .  &  .    & 1     &  .    &  .  &.   &. \\
$\chi_{1d}$  & .  &  .    & .     &  1    &  .  &.   &. \\
$\chi_{2a}$  & .  &  .    & .     &  .    &  1  &.   &. \\
$\chi_{2b}$  & .  &  .    & .     &  .    &  .  &  1 &. \\
$\chi_{2c}$  & .  &  .    & .     &  .    &  .  &  . & 1 \\
$\chi_{8a}$  & 1  &  1    & .     &  .    &  1  &  1 & 1 \\
$\chi_{8b}$  & .  &  .    & 1     &  1    &  1  &  1 & 1 \\
\end{tabular} 
}
\end{center}
\bigskip
\begin{center}
{\rm
\begin{tabular}{r|ccccccc}
      & $P(1a)$ & $P(1b)$  & $P(1c)$  & $P(1d)$ 
      &$P(2a)$ & $P(2b)$ & $P(2c)$ \\
\hline
$1a$  & 2 & 1 & 0 & 0 & 1 & 1 & 1 \\ 
$1b$  & 1 & 2 & 0 & 0 & 1 & 1 & 1 \\ 
$1c$  & 0 & 0 & 2 & 1 & 1 & 1 & 1 \\ 
$1d$  & 0 & 0 & 1 & 2 & 1 & 1 & 1 \\ 
$2a$ & 1 & 1 & 1 & 1 & 3 & 2 & 2 \\ 
$2b$ & 1 & 1 & 1 & 1 & 2 & 3 & 2 \\
$2c$ & 1 & 1 & 1 & 1 & 2 & 2 & 3 \\
\end{tabular}
}
\end{center}
\end{enumerate}
}

\bigskip\noindent
{\bf Proof.}
(i)
This is found using explicit computation with {\sf GAP} \cite{GAP}.
Using the smallest faithful permutation representation
of $G'$ on $100$ points, available in \cite{ModAtlasRep}, 
$P$ can be computed as a Sylow $3$-subgroup of $G'$, 
and hence the normaliser $H'=N_{G'}(P)$ of $P$ 
is easily determined explicitly.
Now the conjugacy classes of $H'$ can be computed,
and its ordinary character table is found using the Dixon-Schneider algorithm.
Note that there are unique conjugacy classes $2B$ and $4B$
consisting of elements of order $2$ and $4$, respectively,
and having centralisers of order $12$ and $4$, respectively.
%This follows by easy calculations by hand or
%by {\sf GAP} \cite{GAP}, compare with {\bf 3.6(v)}.
%
%(iii) Set $S = SD_{16}$. This is defined by
%$S = \langle a, b \, | \, a^2 = b^8 = 1, a^{-1}ba = b^3 \rangle$.
%Then, we know that there is an automorphism $\sigma$ of $S$
%defined by $a^{\sigma} = a$ and $b^{\sigma} = b^5$.
%Then, the conditions for $\sigma$ are satisfied.
%
%The rest is easy.

(ii)--(iii) Easy from the character table.
\quad$\blacksquare$

\bigskip\noindent
{\bf 5.3.Notation.}
We use the notation 
$1_{H'} = \chi_{1a},\chi_{1b},\chi_{1c}, \chi_{1d},
\chi_{2a}, \chi_{2a},\chi_{2c} = {\chi_{2b}}^\vee,
\chi_{8a}, \chi_{8b}$
and 
$1a, 1b, 1c, 1d, 2a, 2b, 2c = 2b^\vee$,
as in {\bf 5.2}.
Namely, we can write
\begin{align*}
{\mathrm{Irr}}(B') 
&= {\mathrm{Irr}}(H')
= \{ 
1_{H'} = \chi_{1a},\chi_{1b},\chi_{1c}, \chi_{1d},
\chi_{2a}, \chi_{2a},\chi_{2c} = {\chi_{2b}}^\vee,
\chi_{8a}, \chi_{8b}
  \},
\notag \\
{\mathrm{IBr}}(B') 
&= {\mathrm{IBr}}(H')
= \{ 
1a, 1b, 1c, 1d, 2a, 2b, 2c = 2b^\vee
 \}
\end{align*}
Let $f'$ and $g'$ 
be the Green correspondences with respect to $(G', P, H')$.

\bigskip\noindent
{\bf 5.4.Lemma.}
{\it
\begin{enumerate}
  \renewcommand{\labelenumi}{(\roman{enumi})}
   \item 
The radical and socle series of PIMs in $B'$ 
are the following:
%JM: introduced boxes
{\rm
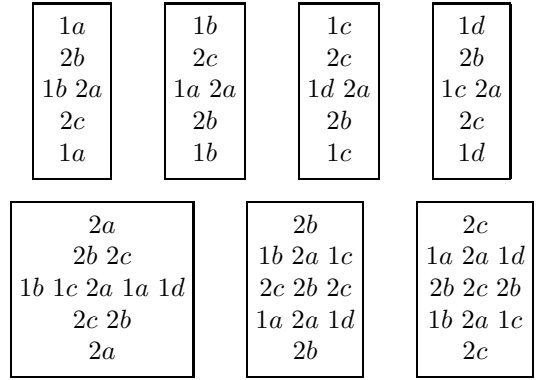
$$
\boxed{
\begin{matrix}
1a \\
2b \\
1b \ 2a \\
2c\\
1a 
\end{matrix}}
\qquad
\boxed{
\begin{matrix}
1b \\
2c \\
1a \ 2a \\
2b\\
1b 
\end{matrix}}
\qquad
\boxed{
\begin{matrix}
1c \\
2c \\
1d \ 2a \\
2b\\
1c 
\end{matrix}}
\qquad
\boxed{
\begin{matrix}
1d \\
2b \\
1c \ 2a \\
2c\\
1d
\end{matrix}}
$$
\smallskip
$$
\boxed{
\begin{matrix}
2a\\
2b \ 2c \\
1b \ 1c \ 2a \ 1a \ 1d\\
2c \ 2b\\
2a
\end{matrix}}
\qquad
\boxed{
\begin{matrix}
2b\\
1b \  2a \ 1c \\
2c \  2b \ 2c\\
1a  \ 2a\ 1d\\
2b
\end{matrix}}
\qquad
\boxed{
\begin{matrix}
2c\\
1a \ 2a \ 1d \\
2b \ 2c \ 2b\\
1b \ 2a \ 1c\\
2c
\end{matrix}}
$$
}
Note that this identifies the simples $2b$ and $2c$ uniquely.

   \item
%JM: For the identification in 8.2
%JM: more cautious ...
An Alperin diagram of the PIM $P(2a)$ is given as follows:
$$ P(2a) \ = \ \begin{matrix} \boxed
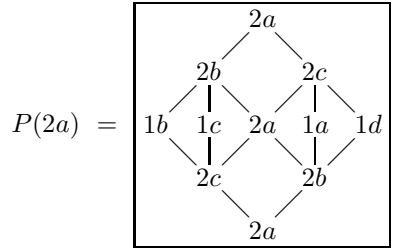
{
\begin{picture}(90,86)(0,0)
\put(40,80){$2a$}
\put(20,60){$2b$}
\put(60,60){$2c$}
\put(0,40){$1b$}
\put(20,40){$1c$}
\put(40,40){$2a$}
\put(60,40){$1a$}
\put(80,40){$1d$}
\put(20,20){$2c$}
\put(60,20){$2b$}
\put(40,0){$2a$}
\put(30,18){\line(1,-1){10}}
\put(50,8){\line(1,1){10}}
\put(10,38){\line(1,-1){10}}
\put(30,28){\line(1,1){10}}
\put(50,38){\line(1,-1){10}}
\put(70,28){\line(1,1){10}}
\put(10,48){\line(1,1){10}}
\put(30,58){\line(1,-1){10}}
\put(50,48){\line(1,1){10}}
\put(70,58){\line(1,-1){10}}
\put(50,78){\line(1,-1){10}}
\put(30,68){\line(1,1){10}}
\put(25,28){\line(0,1){10}}
\put(25,48){\line(0,1){10}}
\put(65,28){\line(0,1){10}}
\put(65,48){\line(0,1){10}}
\end{picture} } \end{matrix} $$
\end{enumerate}
}

\bigskip\noindent
{\bf Proof.}
%JM: More details 
Using the faithful permutation representation of $H'$
obtained in {\bf 5.2}, we have used the {\sf MeatAxe} \cite{MA} 
to construct the PIMs explicitly as matrix representations.
Then we have used the method described in \cite{LuxWie} to 
find the radical and socle series, and the method in \cite{LuxMueRin} 
to compute the whole submodule lattice of $P(2a)$.
%This follows by easy calculations by hand or
%by {\sf GAP} \cite{GAP}, compare with {\bf 3.6(v)}.
%
%(iii) Set $S = SD_{16}$. This is defined by
%$S = \langle a, b \, | \, a^2 = b^8 = 1, a^{-1}ba = b^3 \rangle$.
%Then, we know that there is an automorphism $\sigma$ of $S$
%defined by $a^{\sigma} = a$ and $b^{\sigma} = b^5$.
%Then, the conditions for $\sigma$ are satisfied.
%
%The rest is easy.
\quad$\blacksquare$

\bigskip\noindent
{\bf 5.5.Lemma.}
{\it
\begin{enumerate}
  \renewcommand{\labelenumi}{(\roman{enumi})}
    \item
We can write that
\begin{align*}
{\mathrm{Irr}}(A') 
&=
\{ \chi'_1 = 1_{G'}, \chi'_{154}, \chi'_{22}, \chi'_{1408},
 \chi'_{1925}, \chi'_{770}, \chi'_{3200}, \chi'_{2750},
 \chi'_{1750} \}
\\
{\mathrm{IBr}}(A') 
&=
\{
1_{G'}, 154, 22, 1253, 1176, 748, 321
\}
\end{align*}
     \item
All simples 
$1_{G'}, 154, 22, 1253, 1176, 748, 321$
in $A'$ are self-dual, and have $P$ as their vertices.
     \item
The simples $1_{G'}, 154, 22$ are trivial source $kG'$-modules.
\end{enumerate}
}

\bigskip\noindent
{\bf Proof.}
(i) This was first calculated by %J.F.~
Humphreys \cite[p.329]{Humphreys1982};
see also \cite[$HS$ (mod 3)]{ModularAtlasProject}
and \cite{ModularAtlas}.

(ii) This is obtained by a result of Kn{\"o}rr
\cite[3.7.Corollary]{Knoerr}.

(iii) It follows from \cite{Waki1993} or
\cite[Example 4.8]{Okuyama1997} that the
Green correspondents 
$f'(k_{G'})$, $f'(22)$ and $f'(154)$ are 
$k_{H'} = 1a$, $1b$ and $1c$, respectively,
%Thus we get the assertion, 
see {\bf 5.7} below.
\quad$\blacksquare$

\bigskip\noindent
%JM: avoids an overfull hbox
{\bf 5.6.Notation.} We write %use the notation
$\chi'_1 = 1_{G'}, \chi'_{154}, \chi'_{22}, \chi'_{1408},
 \chi'_{1925}, \chi'_{770}, \chi'_{3200}, \chi'_{2750},
 \chi'_{1750}$,
as well as
$1_{G'}, 154, 22, 1253, 1176, 748, 321$
as in {\bf 5.5}.

\bigskip\noindent
{\bf 5.7.Lemma.}
%JM: introduced boxes
{\it 
$
f'(k_{G'} = 1a) =   \boxed{k_{H'} = 1a} 
\qquad
f'(154)         =   \boxed{1b}
\qquad
f'(22)          =   \boxed{1c}
$
\medskip

$f'(1253) \  = \ 
\begin{matrix}\boxed{
\begin{picture}(50,46)(0,0)
\put(20,40){$2a$}
\put(40,20){$2c$}
\put(0,20){$2b$}
\put(20,0){$2a$}
\put(10,18){\line(1,-1){10}}
\put(30,8){\line(1,1){10}}
\put(10,28){\line(1,1){10}}
\put(30,38){\line(1,-1){10}}
\end{picture}
}\end{matrix}$
\qquad \qquad \qquad
$f'(1176) \ = \ 
\begin{matrix}\boxed{
\begin{picture}(50,46)(0,0)
\put(20,40){$2b$}
\put(40,20){$1c$}
\put(0,20){$1b$}
\put(20,0){$2c$}
\put(10,18){\line(1,-1){10}}
\put(30,8){\line(1,1){10}}
\put(10,28){\line(1,1){10}}
\put(30,38){\line(1,-1){10}}
\end{picture}
}\end{matrix}$

\medskip

$f'(748) \ = \ 
\begin{matrix}\boxed{
\begin{picture}(90,46)(0,0)
\put(80,40){$1d$}
\put(40,40){$2a$}
\put(60,20){$2b$}
\put(20,20){$2c$}
\put(40,0){$2a$}
\put(0,0){$1d$}
\put(50,38){\line(1,-1){10}}
\put(70,28){\line(1,1){10}}
\put(30,28){\line(1,1){10}}
\put(50,8){\line(1,1){10}}
\put(30,18){\line(1,-1){10}}
\put(10,8){\line(1,1){10}}
\end{picture}
}\end{matrix}$
\qquad \quad
$f'(321) \ = \ 
\begin{matrix}\boxed{
\begin{picture}(50,46)(0,0)
\put(20,40){$2c$}
\put(40,20){$1d$}
\put(0,20){$1a$}
\put(20,0){$2b$}
\put(10,18){\line(1,-1){10}}
\put(30,8){\line(1,1){10}}
\put(10,28){\line(1,1){10}}
\put(30,38){\line(1,-1){10}}
\end{picture}
}\end{matrix}$
}

%\begin{center}
%$f'(1253) \quad =$
% \quad\qquad\quad  
%\begin{picture}(100,80)(0,0)
%\put(0,20){$2a$}
%\put(-20,0){$2b$}
%\put(0,-20){$2a$}
%\put(20,0){$2c$}
%
%\put(-12,10){\line(1,1){10}}
%\put(14,21){\line(1,-1){10}}
%\put(-13,-4){\line(1,-1){10}}
%\put(11,-11){\line(1,1){10}}
%\end{picture}
%\end{center}
%
%\bigskip\bigskip\bigskip%
%
%\begin{center}
%$f'(1176)    \ =$ \qquad\qquad
%\begin{picture}(100,80)(0,0)
%\put(0,20){$2b$}
%\put(-20,0){$1b$}
%\put(0,-20){$2c$}
%\put(20,0){$1c$}
%
%\put(-12,10){\line(1,1){10}}
%\put(14,21){\line(1,-1){10}}
%\put(-13,-4){\line(1,-1){10}}
%\put(11,-11){\line(1,1){10}}
%\end{picture}
%\end{center}
%
%\bigskip\bigskip\bigskip
%
%
%\begin{center}
%$f'(748)    \ =$ \qquad\qquad 
%\begin{picture}(100,80)(0,0)
%\put(0,20){$2a$}
%\put(-20,0){$2b$}
%\put(0,-20){$2a$}
%\put(20,0){$2c$}
%\put(40,20){$1d$}
%\put(-40,-20){$1d$}
%
%\put(-12,10){\line(1,1){10}}
%\put(13,18){\line(1,-1){10}}
%\put(-13,-4){\line(1,-1){10}}
%\put(11,-11){\line(1,1){10}}
%\put(28,10){\line(1,1){10}}
%\put(-28,-12){\line(1,1){10}}
%\end{picture}
%\end{center}
%
%\bigskip\bigskip\bigskip%%
%
%
%\begin{center}
%$f'(321) \ =$ \qquad\qquad
%\begin{picture}(100,80)(0,0)
%\put(0,20){$2c$}
%\put(-20,0){$1a$}
%\put(0,-20){$2b$}
%\put(20,0){$1d$}
%
%\put(-12,10){\line(1,1){10}}
%\put(14,21){\line(1,-1){10}}
%\put(-13,-4){\line(1,-1){10}}
%\put(11,-11){\line(1,1){10}}
%\end{picture}
%\end{center}

\bigskip\noindent
{\bf Proof.}
This follows from \cite{Waki1993}, see
\cite[Example 4.8, $HS$]{Okuyama1997}.
\quad$\blacksquare$

\bigskip\noindent
{\bf 5.8.Lemma.}
{\it
The Cartan matrix of $A'$ is the following:}

\medskip
\begin{center}
\begin{tabular}{r|ccccccc}
      & $P(k_{G'})$ & $P(154)$ & $P(22)$ & $P(1253)$ 
      & $P(1176)$ & $P(748)$ & $P(321)$ \\
\hline
$k_{G'}$& 4 & 1 & 1 & 2 & 2 & 2 & 0 \\ 
$154$   & 1 & 3 & 1 & 2 & 0 & 0 & 1 \\ 
$22$    & 1 & 1 & 4 & 2 & 1 & 2 & 1 \\ 
$1253$  & 2 & 2 & 2 & 4 & 2 & 1 & 2 \\ 
$1176$  & 2 & 0 & 1 & 2 & 3 & 2 & 1 \\
$748$   & 2 & 0 & 2 & 1 & 2 & 3 & 0 \\ 
$321$   & 0 & 1 & 1 & 2 & 1 & 0 & 2 \\
\end{tabular}
\end{center}

\bigskip\noindent
{\bf Proof.}
This was first calculated by %J.F.~
Humphreys \cite[p.329]{Humphreys1982};
see also \cite[$HS$ (mod 3)]{ModularAtlasProject}
and \cite{ModularAtlas}.
\qquad$\blacksquare$

\bigskip
%\newpage

\begin{flushleft}
{\bf 6. Stable equivalence between $A$ and $B$ for {\sf HN}}
\end{flushleft}

\bigskip\noindent
{\bf 6.1.Notation.}
First of all, recall the notation
$G$, $A$, $P$, $H$, $B$, $e$, $Q$, $E$, $f$ as in
{\bf 3.3} and {\bf 4.9}.
Let $i$ and $j$ respectively be source idempotents of
$A$ and $B$ with respect to $P$.
As remarked in \cite[pp.821--822]{Linckelmann2001},
we can take $i$ and $j$ such that
${\mathrm{Br}}_P(i){\cdot}e = 
 {\mathrm{Br}}_P(i) \not= 0$ 
and that
${\mathrm{Br}}_P(j){\cdot}e = 
 {\mathrm{Br}}_P(j) \not= 0$.
Set $G_P = C_G(P) = C_H(P) = H_P$,
and set $G_Q = C_G(Q)$ and $H_Q = C_H(Q)$.
By replacing $e_Q$ and $f_Q$ (if necessary), we may assume that
$e_Q$ and $f_Q$ respectively are block idempotents of $kG_Q$ and
$kH_Q$ such that $e_Q$ and $f_Q$ are determined 
by $i$ and $j$, respectively.
Namely, 
${\mathrm Br}_Q(i){\cdot}e_Q =
 {\mathrm Br}_Q(i)$
and
${\mathrm Br}_Q(j){\cdot}f_Q =
 {\mathrm Br}_ Q(j)$.
Let $A_Q = kC_G(Q){\cdot}e_Q$ and 
$B_Q = kC_H(Q){\cdot}f_Q$, so that
$e_Q = 1_{A_Q}$ and $f_Q = 1_{B_Q}$.

\bigskip\noindent
{\bf 6.2.Lemma.} {\it
Let $\mathfrak M_Q$ be a unique (up to isomorphism)
indecomposable direct summand of
$A_Q{\downarrow}^{G_Q \times G_Q}_{G_Q \times H_Q}{\cdot}1_{B_Q}$
with vertex $\Delta P$ (note that such an $\mathfrak M_Q$ always
exists by {\rm\bf{2.5}}).
Then, a pair $(\mathfrak M_Q, \mathfrak M_Q^{\vee})$ induces
a Puig equivalence 
between $A_Q$ and $B_Q$.
}

\bigskip\noindent
{\bf Proof.}
This follows from {\bf 3.4(viii)}, {\bf 2.12(iii)} and
{\bf 2.11(iii)}.
\quad$\blacksquare$

\bigskip\noindent
{\bf 6.3.Lemma.} 
{\it
\begin{enumerate}
  \renewcommand{\labelenumi}{(\roman{enumi})}
    \item
The $(A, B)$-bimodule $1_A{\cdot}kG{\cdot}1_B$ has a
unique (up to isomorphism) indecomposable direct summand
$_A{\mathfrak M}_B$ with vertex $\Delta P$.
Moreover, a functor
$F: {\mathrm{mod}}{\text{-}}A \rightarrow
               {\mathrm{mod}}{\text{-}}B$
defined by $X_A \mapsto (X \otimes_A{\mathfrak M})_B$
induces a splendid stable equivalence of Morita type between
$A$ and $B$. We use the notation $F$ below as well.
   \item
If $X$ is a non-projective trivial source $kG$-module in $A$,
then 
${F}(X) 
= Y \oplus ({\mathrm{proj}})$
for a non-projective indecomposable $kH$-module $Y$ in $B$
such that $Y$ is also a trivial source module, 
and $X$ and $Y$ have a common vertex.
%
%   \item
%
%Suppose that $X$ is an indecomposable $kG$-module in $A$ such that
%a vertex of $X$ belongs to 
%$\mathcal Z (G,P,H)$, see 
%{\rm{\cite[Chapter 4, \S 4]{NagaoTsushima}}}, and that
%$F(X) = Y \oplus 
%({\mathrm{proj}})$
%for a non-projective indecomposable $kH$-module $Y$ in $B$.
%Then, $Y$ is the Fong-Reynolds correspondent of $fX$
%between $B$ and $A_N$.
%
     \item
If $X$ is a non-projective $kG$-module in $A$, then
$F(\Omega X) = \Omega (F(X)) \oplus
({\mathrm{proj}})$.
\end{enumerate}
}

\bigskip\noindent
{\bf Proof.}
This follows just like in 
\cite[Proof of Lemma 6.3]{KoshitaniKunugiWaki2008}.
Namely, we get the assertion by
\cite[Proposition 4.21]{AlperinBroue} and
\cite[Theorem 1.8(i)]{BrouePuig1980} for the morphisms
in the Brauer categories and also 
by
\cite[Theorem]{KoshitaniLinckelmann},
{\bf 6.2}, {\bf 3.4(vi)} and
\cite[Theorem 3.1]{Linckelmann2001},
see \cite[Theorem A.1]{Linckelmann2009}.
\quad$\blacksquare$

\bigskip\noindent
{\bf 6.4.Notation.} We use the notation $\mathfrak M$
and $F$ as in {\bf 6.3}.

\bigskip
%\newpage

\begin{flushleft}
{\bf 7. Images of simples via the functor $F$}
\end{flushleft}

\bigskip\noindent
{\bf 7.1.Lemma.}
{\it
$F(S_1) = 9a$, $F(S_2) = 9b$, $F(S_3) = 9c$.
}

\bigskip\noindent
{\bf Proof.}
These follow from {\bf 4.10}, {\bf 4.11}, {\bf 4.12},
{\bf 2.9} and {\bf 6.3}.
\quad$\blacksquare$

%JM: We need this for the t.s. module in 4.7 as well.
%    We use the same argument as for the t.s. module in 4.8.
\bigskip\noindent
{\bf 7.2.Lemma.}
\begin{enumerate}
  \renewcommand{\labelenumi}{(\roman{enumi})}
    \item
{\it
The trivial source $kG$-module in {\bf 4.8} has $Q \cong C_3$
as its vertex.
}
    \item
{\it
The trivial source $kG$-module in {\bf 4.7} has $Q \cong C_3$
as its vertex.
}
\end{enumerate}

\bigskip\noindent
{\bf Proof.}
(i)
Let $X$ be the trivial source $kG$-module in {\bf 4.8}.
We get by {\bf 6.3(ii)} that
$F(X) = Y \oplus ({\mathrm{proj}})$
for a non-projective indecomposable $B$-module $Y$.
Then, it follows from {\bf 2.7}, {\bf 6.3(i)} and {\bf 7.1}
that
\begin{align*}
0 \ & \not= \ {\mathrm{Hom}}_A(X, S_1) \ 
    \cong \ {\underline{\mathrm{Hom}}}_A(X, S_1) \
    \cong \  {\underline{\mathrm{Hom}}}_B(F(X), F(S_1))  \
\\
   &= \ {\underline{\mathrm{Hom}}}_B(F(X), 9a) \
    = \ {\underline{\mathrm{Hom}}}_B(Y, 9a) \
    \cong \ {\mathrm{Hom}}_B(Y, 9a)
\end{align*}
as $k$-spaces. Clearly, $Y$ is a trivial source $kH$-module
in $B$ by {\bf 6.3(ii)}.

Suppose that $X$ has $P$ as a vertex.
Then, so does $Y$ by {\bf 6.3(ii)}. This yields that
$Y \in \{ 9a, 9b, 9c, 9d, 18a, 18b, 18c\}$ from {\bf 3.8},
and hence $Y \in \{ 9d, 18a, 18b, 18c \}$ by {\bf 7.1}.
But, the above computation shows that
${\mathrm{Hom}}_B(Y, 9a) \not= 0$, a contradiction.

Since $X$ is non-projective, we know that $Q$ is a vertex
of $X$ from {\bf 3.4(vi)}.

(ii)
Let $X'$ be the trivial source $kG$-module in {\bf 4.7}.
We get by {\bf 6.3(ii)} that
$F(X') = Y' \oplus ({\mathrm{proj}})$
for a non-projective indecomposable $B$-module $Y'$.
Then, it follows from {\bf 2.7}, {\bf 6.3(i)} and {\bf 7.1}
that
\begin{align*}
0 \ & \not= \ {\mathrm{Hom}}_A(X', S_3) \ 
    \cong \ {\underline{\mathrm{Hom}}}_A(X', S_3) \
    \cong \  {\underline{\mathrm{Hom}}}_B(F(X'), F(S_3))  \
\\
   &= \ {\underline{\mathrm{Hom}}}_B(F(X'), 9c) \
    = \ {\underline{\mathrm{Hom}}}_B(Y', 9c) \
    \cong \ {\mathrm{Hom}}_B(Y', 9c)
\end{align*}
as $k$-spaces. Clearly, $Y'$ is a trivial source $kH$-module
in $B$ by {\bf 6.3(ii)}.

Suppose that $X'$ has $P$ as a vertex.
Then, so does $Y'$ by {\bf 6.3(ii)}. This yields that
$Y' \in \{ 9a, 9b, 9c, 9d, 18a, 18b, 18c\}$ from {\bf 3.8},
and hence $Y' \in \{ 9d, 18a, 18b, 18c \}$ by {\bf 7.1}.
But, the above computation shows that
${\mathrm{Hom}}_B(Y', 9c) \not= 0$, a contradiction.

Since $X'$ is non-projective, we know that $Q$ is a vertex
of $X'$ from {\bf 3.4(vi)}.

\quad$\blacksquare$

\bigskip\noindent
{\bf 7.3.Lemma.}
{\it
Let $X$ be the trivial source $kG$-module with vertex $Q$
showing up in {\bf 4.8} and {\bf 7.2(i)}.
Then,
$F(X) = V_1 \oplus ({\mathrm{proj}})$,
where $V_1$ is the trivial source $kH$-module in $B$ with
vertex $Q$ given in {\bf 3.8(iii)}. Namely,
$$
F \Biggl( 
\ 
\begin{matrix}
\boxed{
\begin{picture}(51,46)(0,0)
\put(0,0){$S_1$}
\put(40,0){$S_2$}
\put(20,20){$S_4$}
\put(0,40){$S_1$}
\put(40,40){$S_2$}
\put(10,8){\line(1,1){10}}
\put(10,38){\line(1,-1){10}}
\put(30,18){\line(1,-1){10}}
\put(30,28){\line(1,1){10}}
\end{picture}} \end{matrix}
%
%   \boxed
%   {
%  \begin{matrix} S_1 \ \ S_2   \\
%                     S_4       \\
%                 S_1 \ \ S_2 
%  \end{matrix}  
%   }
  \ 
   \Biggr)
\ = \
  \begin{matrix}
\boxed{
\begin{picture}(90,46)(0,0)
\put(40,40){$18a$}
\put(20,20){$18b$}
\put(40,0){$18a$}
\put(60,20){$18c$}
\put(80,40){$9b$}
\put(0,40){$9a$}
\put(0,0){$9b$}
\put(80,0){$9a$}
\put(10,8){\line(1,1){10}}
\put(30,18){\line(1,-1){10}}
\put(50,8){\line(1,1){10}}
\put(70,18){\line(1,-1){10}}
\put(10,38){\line(1,-1){10}}
\put(30,28){\line(1,1){10}}
\put(50,38){\line(1,-1){10}}
\put(70,28){\line(1,1){10}}
\end{picture}}
%  9a \ \ 18a \ \ 9b \\
%     18b \ \ 18c      \\
%  9b \ \ 18a \ \ 9a 
  \end{matrix}
\
\bigoplus
\
({\mathrm{proj}}).
$$
}

\bigskip\noindent
{\bf Proof.}
We know from the proof of {\bf 7.2(i)} that 
$[Y, 9a]^B \not= 0$. Hence we get the assertion 
by {\bf 3.8(iii)}.
\quad$\blacksquare$

\bigskip\noindent
{\bf 7.4.Lemma.}
{\it
Let $X'$ be the trivial source $kG$-module with vertex $Q$
showing up in {\bf 4.7} and {\bf 7.2(ii)}.
Then, 
$F(X') = 
V_2 \oplus ({\mathrm{proj}})$,
where $V_2$ is the trivial source $kH$-module in $B$ 
with vertex $Q$ given
in {\bf 3.8(iii)}. Namely,
$$
F \Biggl( \ \boxed{
          \begin{matrix} S_3 \\
                         S_6 \\
                         S_3 
          \end{matrix}
                 } \ 
   \Biggr)
\ = \ 
\begin{matrix}\boxed{
\begin{picture}(90,46)(0,0)
\put(40,40){$18a$}
\put(20,20){$18c$}
\put(40,0){$18a$}
\put(60,20){$18b$}
\put(80,40){$9d$}
\put(0,40){$9c$}
\put(0,0){$9d$}
\put(80,0){$9c$}
\put(10,8){\line(1,1){10}}
\put(30,18){\line(1,-1){10}}
\put(50,8){\line(1,1){10}}
\put(70,18){\line(1,-1){10}}
\put(10,38){\line(1,-1){10}}
\put(30,28){\line(1,1){10}}
\put(50,38){\line(1,-1){10}}
\put(70,28){\line(1,1){10}}
\end{picture}}
%  9c \ \ 18a \ \ 9d \\
%     18c \ \ 18b      \\
%  9d \ \ 18a \ \ 9c 
  \end{matrix}
\ \bigoplus \
({\mathrm{proj}}).
$$
}

%JM: We now have a simpler proof here.
\bigskip\noindent
{\bf Proof.} 
We know from the proof of {\bf 7.2(ii)} that 
$[Y', 9c]^B \not= 0$. Hence we get the assertion 
by {\bf 3.8(iii)}.
%Just as in the proof of {\bf 7.2},
%we can write that
%$F(X') = Y' \oplus ({\mathrm{proj}})$ for a non-projective
%trivial source $kH$-module $Y'$ in $B$ with vertex $Q$,
%and we get that
%%
%\begin{align*}
%0 \ & \not= \ {\mathrm{Hom}}_A(X', S_3) \ 
%    \cong \ {\underline{\mathrm{Hom}}}_A(X', S_3) \
%    \cong \  {\underline{\mathrm{Hom}}}_B(F(X'), F(S_3))  \
%\\
%   &= \ {\underline{\mathrm{Hom}}}_B(F(X'), 9c) \
%    = \ {\underline{\mathrm{Hom}}}_B(Y', 9c) \
%    \cong \ {\mathrm{Hom}}_B(Y', 9c).
%\end{align*}
%%
%Hence we obtain the assertion from {\bf 3.8(iii)}.
\quad$\blacksquare$

\bigskip\noindent
{\bf 7.5.Lemma.}
{\it It holds that
$F(S_4) \ = \  
%  \boxed{ \begin{matrix}
%                 18a \\
%             18b \ \ 18c  \\
%                 18a
%          \end{matrix}
%        }.$
\begin{matrix}\boxed{
\begin{picture}(55,46)(0,0)
\put(20,40){$18a$}
\put(40,20){$18c$}
\put(0,20){$18b$}
\put(20,0){$18a$}
\put(10,18){\line(1,-1){10}}
\put(30,8){\line(1,1){10}}
\put(10,28){\line(1,1){10}}
\put(30,38){\line(1,-1){10}}
\end{picture}
}\end{matrix}$
}

\bigskip\noindent
{\bf Proof.} 
Let $X$ be the trivial source $kG$-module in $A$
with vertex $Q$ given in {\bf 4.8} and {\bf 7.2(i)}.
By {\bf 7.3}, we can write
$F(X) = V_1 \oplus ({\mathrm{proj}})$, 
where $V_1$ is the trivial source $kH$-module in $B$
given in {\bf 3.8(iii)}.
Then, since $F(S_1) = 9a$ by {\bf 7.1}, it follows
from {\bf 2.8} that
$$
F \Biggl( \
\begin{matrix}\boxed{
\begin{picture}(50,46)(0,0)
\put(40,40){$S_2$}
\put(0,40){$S_1$}
\put(20,20){$S_4$}
\put(20,0){$S_2$}
\put(10,38){\line(1,-1){10}}
\put(30,28){\line(1,1){10}}
\put(25,8){\line(0,1){10}}
\end{picture}
%      \begin{matrix} S_1 \ S_2 \\
%                        S_4    \\
%                        S_2
%      \end{matrix}
}\end{matrix} \
  \Biggr)
 =  
F(X/S_1)  =  ({V_1}/9a) \, \bigoplus \, ({\mathrm{proj}}) 
 = 
\begin{matrix}\boxed{
\begin{picture}(90,46)(0,0)
\put(40,40){$18a$}
\put(20,20){$18b$}
\put(40,0){$18a$}
\put(60,20){$18c$}
\put(80,40){$9b$}
\put(0,40){$9a$}
\put(0,0){$9b$}
\put(10,8){\line(1,1){10}}
\put(30,18){\line(1,-1){10}}
\put(50,8){\line(1,1){10}}
\put(10,38){\line(1,-1){10}}
\put(30,28){\line(1,1){10}}
\put(50,38){\line(1,-1){10}}
\put(70,28){\line(1,1){10}}
\end{picture}
%  9a \ \ 18a \ \ 9b \\
%     18b \ \ 18c      \\
%     9b \ \ 18a  
  } 
  \end{matrix}
\, \bigoplus \,
({\mathrm{proj}}).
$$
Similarly, we get by {\bf 2.8} that 

\medskip
\noindent
%
%\begin{align*}
$
F \Biggl( \ \boxed{\begin{matrix} S_2 \\
                                  S_4 \\
                                  S_2
                  \end{matrix}} \ \Biggr)
=
F \Biggl( {\mathrm{Ker}} \Biggl[ \
              \begin{matrix} 
          \boxed{
\begin{picture}(50,46)(0,0)
\put(40,40){$S_2$}
\put(0,40){$S_1$}
\put(20,20){$S_4$}
\put(20,0){$S_2$}
\put(10,38){\line(1,-1){10}}
\put(30,28){\line(1,1){10}}
\put(25,8){\line(0,1){10}}
\end{picture}
%                             S_1 \ S_2 \\
%                                S_4    \\
%                                S_2 
              }\end{matrix} \twoheadrightarrow  S_1 \Biggr]
  \Biggr)  $ 

\medskip

\hfill
$                               
\cong
 {\mathrm{Ker}} \Biggl( \
\begin{matrix} 
\boxed{
\begin{picture}(90,46)(0,0)
\put(40,40){$18a$}
\put(20,20){$18b$}
\put(40,0){$18a$}
\put(60,20){$18c$}
\put(80,40){$9b$}
\put(0,40){$9a$}
\put(0,0){$9b$}
\put(10,8){\line(1,1){10}}
\put(30,18){\line(1,-1){10}}
\put(50,8){\line(1,1){10}}
\put(10,38){\line(1,-1){10}}
\put(30,28){\line(1,1){10}}
\put(50,38){\line(1,-1){10}}
\put(70,28){\line(1,1){10}}
\end{picture}
%                   9a \ 18a \ 9b \\
%                       18b \ 18c   \\
%                        9a \ 18a
 }     \end{matrix}   \twoheadrightarrow   9a
                \Biggr)
     \bigoplus  ({\mathrm{proj}})          
\ = \  
\begin{matrix} 
\boxed{
\begin{picture}(90,46)(0,0)
\put(40,40){$18a$}
\put(20,20){$18b$}
\put(40,0){$18a$}
\put(60,20){$18c$}
\put(80,40){$9b$}
\put(0,0){$9b$}
\put(10,8){\line(1,1){10}}
\put(30,18){\line(1,-1){10}}
\put(50,8){\line(1,1){10}}
\put(30,28){\line(1,1){10}}
\put(50,38){\line(1,-1){10}}
\put(70,28){\line(1,1){10}}
\end{picture}
%                       18a \ 9b  \\
%                         18b \ 18c \\
%                         9b  \ 18a
         } \end{matrix} 
 \,  \bigoplus \, ({\mathrm{proj}}).            
%\end{align*}
$

\medskip

Then, since $F(S_2) = 9b$ by {\bf 7.1}, we similarly 
obtain by {\bf 2.8} that
$$
F(S_4) \ = \  
\begin{matrix}
  \boxed{ 
\begin{picture}(55,46)(0,0)
\put(20,40){$18a$}
\put(40,20){$18c$}
\put(0,20){$18b$}
\put(20,0){$18a$}
\put(10,18){\line(1,-1){10}}
\put(30,8){\line(1,1){10}}
\put(10,28){\line(1,1){10}}
\put(30,38){\line(1,-1){10}}
\end{picture}

%                 18a \\
%             18b \ \ 18c  \\
%                 18a
        }
         \end{matrix}
\, \bigoplus \, 
({\mathrm{proj}}).
$$
Therefore, {\bf 6.3(i)} and {\bf 2.9} imply the assertion.
\quad$\blacksquare$

\bigskip\noindent
{\bf 7.6.Lemma.}
{\it
It holds that
$F(S_6) \ = \ 
  \begin{matrix} 
  \boxed{
\begin{picture}(90,46)(0,0)
\put(40,40){$18a$}
\put(20,20){$18c$}
\put(40,0){$18a$}
\put(60,20){$18b$}
\put(80,40){$9d$}
\put(0,0){$9d$}
\put(10,8){\line(1,1){10}}
\put(30,18){\line(1,-1){10}}
\put(50,8){\line(1,1){10}}
\put(30,28){\line(1,1){10}}
\put(50,38){\line(1,-1){10}}
\put(70,28){\line(1,1){10}}
\end{picture}
%                &     & 18a &     & 9d \\
%                  & 18c &     & 18b &    \\
%               9d &     & 18a &     &
        }
  \end{matrix}
$.
}

\bigskip\noindent
{\bf Proof.}
Let $X' = \boxed{\begin{matrix} S_3 \\ S_6 \\ S_3 \end{matrix}}$ 
%U(S_3, S_6, S_3)$ 
in {\bf 4.7}, that is, $X'$ is a
trivial source $kG$-module in $A$ with vertex $Q$. Then,
{\bf 7.4} yields that
$$
F(X') \ = \ 
  \begin{matrix}
\boxed{
\begin{picture}(90,46)(0,0)
\put(40,40){$18a$}
\put(20,20){$18c$}
\put(40,0){$18a$}
\put(60,20){$18b$}
\put(80,40){$9d$}
\put(0,40){$9c$}
\put(0,0){$9d$}
\put(80,0){$9c$}
\put(10,8){\line(1,1){10}}
\put(30,18){\line(1,-1){10}}
\put(50,8){\line(1,1){10}}
\put(70,18){\line(1,-1){10}}
\put(10,38){\line(1,-1){10}}
\put(30,28){\line(1,1){10}}
\put(50,38){\line(1,-1){10}}
\put(70,28){\line(1,1){10}}
\end{picture}
%  9c \ \ 18a \ \ 9d \\
%     18c \ \ 18b      \\
%  9d \ \ 18a \ \ 9c 
       }
  \end{matrix}
\ \bigoplus \
({\mathrm{proj}}).
$$
Since $F(S_3) = 9c$ by {\bf 7.1}, we obtain the assertion
from {\bf 6.3(i)} and {\bf 2.9} just as in the proof of
{\bf 7.5}.
\quad$\blacksquare$

\bigskip\noindent
{\bf 7.7.Notation.}
We use the notation $W = F(S_5) \oplus F(S_7)$
in the rest of this paper.

\bigskip\noindent
{\bf 7.8.Lemma.}
{\it We get the following.
\begin{enumerate}
  \renewcommand{\labelenumi}{(\roman{enumi})}
    \item
The module $W$ is self-dual.
    \item
The module $W$ is a direct sum of exactly two 
non-projective non-simple indecomposable $B$-modules, 
and both of them are self-dual.
    \item
It holds that $F(S_5)$ and $F(S_7)$ are neither
simple $B$-modules, and
$2 \leqslant j(W) \leqslant 4$.
    \item
$[9x, W]^B = [W, 9x]^B = 0$ for any
$x \in \{ a, b, c \}$.
    \item
$[18a, W]^B = [W, 18a]^B = 0$.
\end{enumerate}
}

\bigskip\noindent
{\bf Proof.}
(i) This follows from {\bf 4.3(i)} and {\bf 6.3(i)}.

(ii) This follows from {\bf 2.9}, {\bf 6.3(i)}, {\bf 4.3(i)},
{\bf 2.3(i)} and {\bf 4.1}.

(iii) By {\bf (ii)} and {\bf 3.6(vi)}, we get
$j(W) \leqslant 4$. 
Assume that $F(S_5)$ is simple. Then, we know by
{\bf 3.8(ii)} and {\bf 6.3(ii)} that
$S_5$ is a trivial source module, and hence
$S_5$ lifts to a trivial source $\mathcal OG$-module
by {\bf 2.3(i)}. This contradicts the
$3$-decomposition matrix in {\bf 4.1}. 
Hence, $F(S_5)$ is not simple. Similarly, we know that
$F(S_7)$ is not simple. These imply $j(W) \geqslant 2$.

(iv) This is obtained by {\bf 7.1} and {\bf 2.13}.

(v) Set $X = F(S_4)$. By {\bf 7.5}, there is an
epimorphism $X \twoheadrightarrow 18a$.
Hence, we get from {\bf 2.7} and {\bf 7.5} that
\begin{align*}
  {\mathrm{Hom}}_B(18a, W) 
&\cong \
  {\underline{\mathrm{Hom}}}_B(18a, W) 
\  \cong \
  {\underline{\mathrm{Hom}}}_A 
    \Big(F^{-1}(18a), F^{-1}(W) \Big)
\\
&=  \
  {\underline{\mathrm{Hom}}}_A 
    \Big(F^{-1}(18a), S_5 \oplus S_7 \Big)
\  \cong \
  {\mathrm{Hom}}_A 
    \Big(F^{-1}(18a), S_5 \oplus S_7 \Big)
\\
&\subseteq \
 {\mathrm{Hom}}_A 
%JM: use X instead
    \Biggl(F^{-1} (X),
%        \Biggl( \ 
%         \boxed{
%       \begin{matrix} 18a \\ 18b \ 18c \\ 18a \end{matrix}
%               } \ 
%                \Biggr)
      S_5 \oplus S_7 \Biggr)
\  = \   
{\mathrm{Hom}}_A (S_4, S_5 \oplus S_7)
\  =  \ 0.
\quad\blacksquare
\end{align*}

\bigskip\noindent
{\bf 7.9.Notation.}
Let 
$M \ = \ 
 \begin{matrix} 
\boxed{ 
\begin{picture}(70,46)(0,0)
\put(30,40){$S_4$}
\put(30,0){$S_4$}
\put(60,20){$S_7$}
\put(40,20){$S_5$}
\put(20,20){$S_2$}
\put(0,20){$S_1$}
\put(10,18){\line(2,-1){20}}
\put(40,8){\line(2,1){20}}
\put(10,28){\line(2,1){20}}
\put(40,38){\line(2,-1){20}}
\put(25,18){\line(2,-3){7}}
\put(45,18){\line(-2,-3){7}}
\put(25,28){\line(2,3){7}}
\put(45,28){\line(-2,3){7}}
\end{picture}
%S_4 \\ S_1 \ S_2 \ S_5 \ S_7 \\ S_4
         }
 \end{matrix}
$
be the trivial source $kG$-module in $A$ showing up
in {\bf 4.6}, and set $\mathfrak X_B = F(M)$ and we
use the notation $\mathfrak X$ in the rest of this paper.

\bigskip\noindent
{\bf 7.10.Lemma.}
{\it
\begin{enumerate}
  \renewcommand{\labelenumi}{(\roman{enumi})}
    \item
The module ${\mathfrak X}$ has a filtration
$$
{\mathfrak X} \ \ = \ \ 
\begin{matrix}
\begin{matrix} 
    \boxed{ 
\begin{picture}(55,46)(0,0)
\put(20,40){$18a$}
\put(40,20){$18c$}
\put(0,20){$18b$}
\put(20,0){$18a$}
\put(10,18){\line(1,-1){10}}
\put(30,8){\line(1,1){10}}
\put(10,28){\line(1,1){10}}
\put(30,38){\line(1,-1){10}}
\end{picture}
%18a \\ 18b \ 18c \\ 18a
          }
            \end{matrix}
   \\
   \\
\end{matrix}
\
\Biggm|
\begin{matrix}
     \\
    9a \oplus 9b \oplus W 
     \\
\end{matrix}
\
{
\begin{matrix} 
  \\
  \\
\Biggm| \ 
\end{matrix}
}
\begin{matrix}
     \\  
     \\
\begin{matrix}
    \boxed{ 
\begin{picture}(55,46)(0,0)
\put(20,40){$18a$}
\put(40,20){$18c$}
\put(0,20){$18b$}
\put(20,0){$18a$}
\put(10,18){\line(1,-1){10}}
\put(30,8){\line(1,1){10}}
\put(10,28){\line(1,1){10}}
\put(30,38){\line(1,-1){10}}
\end{picture}
% 18a \\ 18b \ 18c \\ 18a
          }
            \end{matrix}
\end{matrix}
{
\begin{matrix}
 \\ \\ \\ \\ \\ 
 \end{matrix}
}
$$
namely, ${\mathfrak X}$ has submodules %$Y$ and $Z$ such that
${\mathfrak X} \supsetneqq Y \supsetneqq Z$ such that
${\mathfrak X}/Y \cong Z \cong
\begin{matrix} 
\boxed{ 
\begin{picture}(55,46)(0,0)
\put(20,40){$18a$}
\put(40,20){$18c$}
\put(0,20){$18b$}
\put(20,0){$18a$}
\put(10,18){\line(1,-1){10}}
\put(30,8){\line(1,1){10}}
\put(10,28){\line(1,1){10}}
\put(30,38){\line(1,-1){10}}
\end{picture}
%18a \\ 18b \ 18c \\ 18a
          }
            \end{matrix}
$
and $Y/Z \cong 9a \oplus 9b \oplus W$.
   \item
It holds
${\mathfrak X} = V \oplus P(18a)$
where $V \in \{ V_3, V_4 \}$.
\end{enumerate}
}

\bigskip\noindent
{\bf Proof.}
(i) This follows from {\bf 4.6}, {\bf 7.1} and {\bf 7.5}.

(ii) We know by {\bf 6.3(ii)} that 
${\mathfrak X} = V \oplus L$ for an indecomposable
$kH$-module $V$ in $B$ with vertex $Q$ and a
projective $kH$-module $L$ in $B$.
Note that $V_i \, {\not|} \, {\mathfrak X}$ for $i = 1, 2$
by {\bf 7.3} and {\bf 7.4}.
Thus, $V \in \{ V_3, V_4 \}$ by {\bf 3.8(iii)}.
Moreover, since $[V_3, 18a]^B = [V_4, 18a]^B = 0$ by {\bf 3.8(iii)},
we know that $[V, 18a]^B = [V, 18a]^B = 0$ again by {\bf 3.8(iii)}.
Thus, we have $P(18a) | L$ by {\bf (i)}, 
and hence $P(18a) | {\mathfrak X}$.

Next, assume that $P(T) | L$ 
for a simple $kH$-module $T$ in $B$ with $T \not\cong 18a$.
Since $Z$ has a unique minimal submodule, and which
is isomorphic to $18a$, we have that 
$P(T) \cap Z = 0$ in ${\mathfrak X}$, and hence that there is a
direct sum $P(T) \oplus Z$ in ${\mathfrak X}$. Set $\bar {\mathfrak X} = {\mathfrak X}/Z$.
Clearly, $\bar {\mathfrak X} \supseteq (P(T) \oplus Z)/Z \cong P(T)$.
Since $P(T)$ is injective, it holds $P(T) | \bar {\mathfrak X}$.
Set $U = (\bar {\mathfrak X})^\vee$. Then, by the dualities, we know
$P(T^\vee) | U$. 
Now, by the filtration of ${\mathfrak X}$, $U$ has a filtration
$$
  U \ = \ 
        9a \oplus 9b \oplus W \ \Biggm| \ 
        \begin{matrix} \\ \\
\begin{matrix} 
           \boxed{ 
\begin{picture}(55,46)(0,0)
\put(20,40){$18a$}
\put(40,20){$18c$}
\put(0,20){$18b$}
\put(20,0){$18a$}
\put(10,18){\line(1,-1){10}}
\put(30,8){\line(1,1){10}}
\put(10,28){\line(1,1){10}}
\put(30,38){\line(1,-1){10}}
\end{picture}
%18a \\ 18b \ 18c \\ 18a 
                 }
                   \end{matrix}
        \end{matrix}
$$
Namely, $U$ has a submodule $Z'$ such that
$$
Z' \cong 
\begin{matrix} 
\boxed{ 
\begin{picture}(55,46)(0,0)
\put(20,40){$18a$}
\put(40,20){$18c$}
\put(0,20){$18b$}
\put(20,0){$18a$}
\put(10,18){\line(1,-1){10}}
\put(30,8){\line(1,1){10}}
\put(10,28){\line(1,1){10}}
\put(30,38){\line(1,-1){10}}
\end{picture}
%18a \\ 18b \ 18c \\ 18a 
                 }
                   \end{matrix}
\quad {\text{and}} \quad 
U/Z' \cong 9a \oplus 9b \oplus W.
$$
We have $T^\vee \not\cong 18a$ by {\bf 3.6(iii)}.
Hence, we get $P(T^\vee) \cap Z' = 0$ in $U$, and hence
there is a direct sum
$P(T^\vee) \oplus Z' \subseteq U$. Then, we have
$$
P(T^\vee) \cong (P(T^\vee) \oplus Z')/Z' \subseteq U/Z' 
\cong 9a \oplus 9b \oplus W.
$$
Since $P(T^\vee)$ is injective, it holds that
$P(T^\vee) | (9a \oplus 9b \oplus W)$, so that
$P(T^\vee) | W$ by {\bf 3.6(vi)}.
This is a contradiction by {\bf 7.8(ii)}.

Now, assume that  $(P(18a) \oplus P(18a)) | {\mathfrak X}$.
Then, since ${\mathrm{soc}}(Z) \cong 18a$, it follows from
{\bf 2.14} that
$$
P(18a) \ \Bigg| \ 
{\mathfrak X}/Z \quad = \quad
     \Biggl( \  
\begin{matrix} 
\boxed{ 
\begin{picture}(55,46)(0,0)
\put(20,40){$18a$}
\put(40,20){$18c$}
\put(0,20){$18b$}
\put(20,0){$18a$}
\put(10,18){\line(1,-1){10}}
\put(30,8){\line(1,1){10}}
\put(10,28){\line(1,1){10}}
\put(30,38){\line(1,-1){10}}
\end{picture}
%18a \\ 18b \ 18c \\ 18a
         }
             \end{matrix}
     \
{\begin{matrix} \\ \\ \Bigg| \end{matrix}}
     \
{\begin{matrix} \\ \\ 
  9a \oplus 9b \oplus W
 \end{matrix}
}
\Biggr)
.   
$$
Then, by taking its dual, we get also that
$$
P(18a) \ \Bigg| \ 
({\mathfrak X}/Z)^\vee \quad = \quad
\Biggl(
{
\begin{matrix} 
9a \oplus 9b \oplus W
\\ \\
 \end{matrix}
}
\
{\begin{matrix} \\ \Bigg| \end{matrix}}
\
{\begin{matrix} \\ \\ Z \end{matrix}}
\Biggr)
$$
where the right-hand-side is a filtration,
by using {\bf 7.8(i)} and {\bf 3.6(iii)}. 
Set $N = ({\mathfrak X}/Z)^\vee$. Then, we may consider that $N$ has a
$B$-submodule $Z$ such that $N/Z \cong 9a \oplus 9b \oplus W$
and $N = P(18a) \oplus N'$ for a $B$-submodule $N'$ of $N$.
Since $j(Z) = 3$, it holds
$Z \subseteq {\mathrm{soc}}_3(N)
   = {\mathrm{soc}}_3(P(18a)) \oplus {\mathrm{soc}}_3(N')$.
This implies that there exists a $B$-epimorphism
$\pi: N/Z \twoheadrightarrow N/{\mathrm{soc}}_3(N)$.
Clearly,
\begin{align*}
  N/{\mathrm{soc}}_3(N) 
\ =& \
    [P(18a) \oplus N']/
    [{\mathrm{soc}}_3(P(18a)) \oplus {\mathrm{soc}}_3(N')]
\\
\ \cong& \
    [P(18a)/{\mathrm{soc}}_3(P(18a)] 
     \oplus [N'/{\mathrm{soc}}_3(N')].
\end{align*}
Since 
$P(18a)/{\mathrm{soc}}_3(P(18a)) = 
    \boxed{ \begin{matrix} 18a \\ 18b \ 18c \end{matrix} 
          }$
by {\bf 3.6(vi)}, we get that
$18a \, | \, [(N/Z)/{\mathrm{rad}}(N/Z)] 
    \cong 9a \oplus 9b \oplus [W/{\mathrm{rad}}(W)]$.
This shows that $[W, 18a]^B \not= 0$, which is a contradiction
by {\bf 7.8(v)}.

Thus, we know $[P(18a) | L]^B = 1$.
Therefore, we get $L \cong P(18a)$. We are done.
\quad$\blacksquare$

\bigskip\noindent
{\bf 7.11.Lemma.}
{\it 
$W/{\mathrm{rad}}(W) \cong {\mathrm{soc}}(W) \cong 18b \oplus 18c$.
}

\bigskip\noindent
{\bf Proof.}
By {\bf 7.8(i)} and {\bf 3.6(iii)}, it suffices to show
only $W/{\mathrm{rad}}(W) \cong 18b \oplus 18c$.
By {\bf 7.10(ii)}, we have
\begin{equation}
\mathfrak X \ = \ 
V \oplus P(18a) \ = \
        \begin{matrix} 
\boxed{
\begin{picture}(90,46)(0,0)
\put(60,40){$18c$}
\put(20,40){$18b$}
\put(80,20){$9y$}
\put(40,20){$18a$}
\put(0,20){$9x$}
\put(60,0){$18b$}
\put(20,0){$18c$}
\put(10,18){\line(1,-1){10}}
\put(50,18){\line(1,-1){10}}
\put(30,8){\line(1,1){10}}
\put(70,8){\line(1,1){10}}
\put(10,28){\line(1,1){10}}
\put(50,28){\line(1,1){10}}
\put(30,38){\line(1,-1){10}}
\put(70,38){\line(1,-1){10}}
\end{picture}

%   18b \ \ 18c     \\
%                         9x \ 18a \ 9y    \\
%                          18c \ \ 18b 
      }
        \end{matrix}
      \ \oplus \ P(18a),
\quad
{\text{where}} \ 
%JM: replaced (d,c) by (c,d)
(x,y) \in \{ (b,a), (c,d) \}.
\end{equation}
By {\bf 7.10(i)}, $\mathfrak X$ has a filtration
\begin{equation}
{\mathfrak X} \ \ = \ \ 
\begin{matrix}
\begin{matrix} 
    \boxed{ 
\begin{picture}(55,46)(0,0)
\put(20,40){$18a$}
\put(40,20){$18c$}
\put(0,20){$18b$}
\put(20,0){$18a$}
\put(10,18){\line(1,-1){10}}
\put(30,8){\line(1,1){10}}
\put(10,28){\line(1,1){10}}
\put(30,38){\line(1,-1){10}}
\end{picture}
%18a \\ 18b \ 18c \\ 18a
          }
            \end{matrix}
   \\
   \\
\end{matrix}
\
\Biggm|
\begin{matrix}
     \\
    9a \oplus 9b \oplus W 
     \\
\end{matrix}
\
{
\begin{matrix} 
  \\
  \\
\Biggm| \ 
\end{matrix}
}
\begin{matrix}
     \\  
     \\
\begin{matrix} 
   \boxed{ 
\begin{picture}(55,46)(0,0)
\put(20,40){$18a$}
\put(40,20){$18c$}
\put(0,20){$18b$}
\put(20,0){$18a$}
\put(10,18){\line(1,-1){10}}
\put(30,8){\line(1,1){10}}
\put(10,28){\line(1,1){10}}
\put(30,38){\line(1,-1){10}}
\end{picture}
%18a \\ 18b \ 18c \\ 18a
          }
            \end{matrix}
\end{matrix}
{
\begin{matrix}
 \\ \\ \\ \\ \\ 
 \end{matrix}
}
\end{equation}
Set $L_i(W) = {\mathrm{rad}}^i(W)/{\mathrm{rad}}^{i+1}(W)$
for each $i = 0, 1, ...$.
Then (13) and (14) show that $(18b \oplus 18c) | L_1(W)$.
%
%Thus, by {\bf 3.6(vi)} and {\bf 7.8(iv)-(v)}, it holds that
%$L_1(W) \cong (m \times 18b) \oplus (n \times 18c)$
%for some positive integers $m$ and $n$.
%
Recall $(18b)^\vee \cong 18c$ by {\bf 3.6(iii)}.

Suppose that $(18b \oplus 18b) \, | \, L_1(W)$. Then, by (14),
$\mathfrak X$ has a factor module 
$\bar{\mathfrak X}$ 
which has a filtration
$$
\begin{matrix} 
\boxed{ 
\begin{picture}(55,46)(0,0)
\put(20,40){$18a$}
\put(40,20){$18c$}
\put(0,20){$18b$}
\put(20,0){$18a$}
\put(10,18){\line(1,-1){10}}
\put(30,8){\line(1,1){10}}
\put(10,28){\line(1,1){10}}
\put(30,38){\line(1,-1){10}}
\end{picture}
%18a \\ 18b \ 18c \\ 18a 
      }
\end{matrix}
\ {\begin{matrix} \\ \Biggm|  \end{matrix}} \
      \begin{matrix} \\ \\  \\ 18b \oplus 18b \end{matrix}
$$
Since $[\mathfrak X, 18b]^B = 1$ by (13), and since there 
do not exist modules of forms
$\boxed{ \begin{matrix} 18b \\ 18b \end{matrix} }$
nor
$\boxed{ \begin{matrix} 18c \\ 18b \end{matrix} }$
by {\bf 3.6(vi)}, there must be a $kH$-module
having radical and socle series
$$
  \boxed{
       \begin{matrix} 18a \\ 18b \ 18c \\ 18a \\ 18b
       \end{matrix}
        }
$$
But this is a contradiction by {\bf 3.10(i)}.

Similarly, we get a contradiction by using {\bf 3.10(ii)}
if $(18c \oplus 18c) | L_1(W)$.

%JM: Delete a , to avoid an overfull hbox
Thus it holds that
$[W, 18b]^B = [W, 18c]^B = 1$ and 
$[W, T]^B = 0$ for any $T \in \{ 9a, 9b, 9c, 18a \}$
by {\bf 3.6(iii)} and {\bf 7.8(iv)-(v)}.
However, we have to investigate for $9d$.

Assume, first, that the case $(x, y) = (b, a)$ happens
in (13). Then, (14) and (13) imply that 
$W = 9a + 9b + 9c + 9d + 2 \times 18b + 2 \times 18c$,
as composition factors.

Suppose that $9d \, | \, L_1(W)$. Then, since $c_W(9d) = 1$ and
since $W$ and $9d$ are both self-dual by {\bf 3.6(iii)} and 
{\bf 7.8(i)}, we get that $9d \, | \, W = F(S_5) \oplus F(S_7)$.
Recall that $F(S_5)$ and $F(S_7)$ are both 
non-projective indecomposable $kH$-modules by {\bf 2.9}
and {\bf 6.3(i)-(ii)}.
Since $9d$ is a trivial source $kH$-module by {\bf 3.8(ii)},
we know by {\bf 6.3(ii)} that $S_5$ or $S_7$ is a
trivial source module, and hence that $S_5$ or $S_7$
lifts from $k$ to $\mathcal O$ by {\bf 2.3(i)}.
This is a contradiction by the $3$-decomposition matrix in
{\bf 4.1}. 

Hence, $9d {\not|} \, L_1(W)$. 
This yields $L_1(W) \cong 18b \oplus 18c$.

%JM: Interchanged c and d
Next, assume that the case $(x,y) = (c,d)$ in (13) happens.
Then, (13) and (14) imply that
\begin{align}
W = 2 \times 9c + 2 \times 9d + 2 \times 18b + 2 \times 18c,
\qquad {\text{as \ composition \ factors}}.
\end{align}

Suppose that $(9d \oplus 9d) \, | \, L_1(W)$. 
Then, the self-dualities
of $9d$ and $W$ in {\bf 3.6(iii)} and {\bf 7.8(i)} imply that
$(9d \oplus 9d) \, | \, W = F(S_5) \oplus F(S_7)$.
Hence, $W \cong 9d \oplus 9d$ by {\bf 7.8(ii)},
contradicting {\bf 7.7} and {\bf 7.8}.

Thus, 
\begin{align}
[W, 9d]^B \leqslant 1.
\end{align}

Assume, next, that $[W, 9d]^B =1$.
Hence, by the dualities in {\bf 3.6(iii)}, we have
\begin{align}
  L_1(W) \cong {\mathrm{soc}}(W) \cong 18b \oplus 18c \oplus 9d. 
\end{align}
We get by {\bf 7.8} that $W = W_1 \oplus W_2$ where $W_i$ is
a non-simple non-projective indecomposable self-dual $B$-module
for $i = 1, 2$. 
Thus, by (17) and by interchanging $W_1$ and $W_2$,
we may assume that
$L_1 (W_1) \cong 18b, 18c$ or $9d$.

{\bf Case 1:} $L_1 (W_1) \cong 18b$. Then,
${\mathrm{soc}}(W_1) \cong 18c$ 
since $(18b)^{\vee} \cong 18c$ by
{\bf 3.6(iii)} and since $W_1$ is self-dual. Hence, the structure
of $P(18b)$ in {\bf 3.6(vi)} yields that
$W_1 = \boxed{\begin{matrix} 18b \\ 9c \\ 18c \end{matrix}}$.
%U(18b, 9c, 18c)$. 
Hence, (15) and (17) imply that
$L_1(W_2) \cong 18c \oplus 9d$ and $L_2(W_2) \cong 9c$.
But this is a contradiction since 
%JM: Replaced kG by B
${\mathrm{Ext}}^1_{B}(18c, 9c) = 0 =
 {\mathrm{Ext}}^1_{B}(9d, 9c)$ by {\bf 3.6(vi)}.

{\bf Case 2:} $L_1(W_1) \cong 18c$.
As in {\bf Case 1}, we know that 
$W_1 = \boxed{\begin{matrix} 18c \\ 9c \\ 18b \end{matrix}}$.
%U(18c, 9c, 18b)$.
Then we get a contradiction by {\bf 3.6(vi)} as in {\bf Case 1}.

{\bf Case 3:} $L_1(W_1) \cong 9d$.
By the self-dualities of $W_1$ in {\bf 7.8(ii)} 
and simple $B$-modules in {\bf 3.6(iii)}, we get that
${\mathrm{soc}}(W_1) \cong 9d$.
It follows by {\bf 2.16} that 
${\mathrm{soc}}(W_1) \subseteq {\mathrm{rad}}(W_1)$.
Hence $c_{W_1}(9d) = 2$ by (15). Thus, the structure of $P(9d)$
in {\bf 3.6(vi)} yields that $W_1 \cong P(9d)$,
a contradiction.

Therefore $[W, 9d]^B \not= 1$, and hence
$[W, 9d]^B = 0$ by (16).
So that we have $L_1(W) \cong 18b \oplus 18c$.
\quad$\blacksquare$

\bigskip\noindent
{\bf 7.12.Lemma.}
{\it
${\mathfrak X} = V_3 \oplus P(18a)$.
}

\bigskip\noindent
{\bf Proof.}
Suppose that ${\mathfrak X} = V_4 \oplus P(18a)$. Then, we get by
{\bf 7.10(i)-(ii)} and {\bf 3.6(iv)} that
$W = 
2 \times 9c + 2 \times 9d + 2 \times 18b + 2 \times 18c$,
as composition factors. 
We use the same notation $L_i(W)$ as in the proof or {\bf 7.11}.
By {\bf 7.11}, $L_1(W) \cong 18b \oplus 18c$. 
Since $c_W(9c) = 2$, it follows from {\bf 3.6(vi)} and
{\bf 7.8(iii)} that $j(W) = 4$ and $9c \, | \, L_4(W)$.
This means $9c \, | \, {\mathrm{soc}}(W)$, 
contradicting {\bf 7.11}.
Therefore, we get the assertion by {\bf 7.10(ii)}.
\quad$\blacksquare$

\bigskip\noindent
{\bf 7.13.Lemma.}
{\it
$W \ = \ 
\begin{matrix} 
    \boxed{ 
\begin{picture}(50,46)(0,0)
\put(20,40){$18b$}
\put(40,20){$9c$}
\put(0,20){$9b$}
\put(20,0){$18c$}
\put(10,18){\line(1,-1){10}}
\put(30,8){\line(1,1){10}}
\put(10,28){\line(1,1){10}}
\put(30,38){\line(1,-1){10}}
\end{picture}
%18b \\ 9b \ 9c \\ 18c
          }
            \end{matrix}
    \oplus
 \begin{matrix} 
    \boxed{
\begin{picture}(50,46)(0,0)
\put(20,40){$18c$}
\put(40,20){$9d$}
\put(0,20){$9a$}
\put(20,0){$18b$}
\put(10,18){\line(1,-1){10}}
\put(30,8){\line(1,1){10}}
\put(10,28){\line(1,1){10}}
\put(30,38){\line(1,-1){10}}
\end{picture}
%18c \\ 9a \ 9d \\ 18b
          }
            \end{matrix}
\ = \ 
F(S_5) \oplus F(S_7)$.
%JM: delete newline to avoid underfull hbox 

\smallskip\noindent
Namely, either one of the following two cases occurs:
$$
\begin{matrix}
{\mathrm{Case}} \ ({\mathrm{a}}):& & &      F(S_5) &= \ \ 
\begin{matrix} 
          \boxed{ 
\begin{picture}(50,46)(0,0)
\put(20,40){$18b$}
\put(40,20){$9c$}
\put(0,20){$9b$}
\put(20,0){$18c$}
\put(10,18){\line(1,-1){10}}
\put(30,8){\line(1,1){10}}
\put(10,28){\line(1,1){10}}
\put(30,38){\line(1,-1){10}}
\end{picture}
%18b \\ 9b \ 9c \\ 18c
                }
                  \end{matrix}
& \mathrm{and} &
                           F(S_7) &= \ \  
\begin{matrix} 
           \boxed{ 
\begin{picture}(50,46)(0,0)
\put(20,40){$18c$}
\put(40,20){$9d$}
\put(0,20){$9a$}
\put(20,0){$18b$}
\put(10,18){\line(1,-1){10}}
\put(30,8){\line(1,1){10}}
\put(10,28){\line(1,1){10}}
\put(30,38){\line(1,-1){10}}
\end{picture}
%18c \\ 9a \ 9d \\ 18b
                }
                   \end{matrix}
\\     
& & & & \\
{\mathrm{Case}} \ ({\mathrm{b}}):&  & &    F(S_5) &= \ \
\begin{matrix}
          \boxed{ 
\begin{picture}(50,46)(0,0)
\put(20,40){$18c$}
\put(40,20){$9d$}
\put(0,20){$9a$}
\put(20,0){$18b$}
\put(10,18){\line(1,-1){10}}
\put(30,8){\line(1,1){10}}
\put(10,28){\line(1,1){10}}
\put(30,38){\line(1,-1){10}}
\end{picture}
% 18c \\ 9a \ 9d \\ 18b
                }
                  \end{matrix}
& \mathrm{and} &
            F(S_7) & = \ \ 
\begin{matrix} 
           \boxed{ 
\begin{picture}(50,46)(0,0)
\put(20,40){$18b$}
\put(40,20){$9c$}
\put(0,20){$9b$}
\put(20,0){$18c$}
\put(10,18){\line(1,-1){10}}
\put(30,8){\line(1,1){10}}
\put(10,28){\line(1,1){10}}
\put(30,38){\line(1,-1){10}}
\end{picture}
%18b \\ 9b \ 9c \\ 18c
                }
                   \end{matrix}
\end{matrix}
$$
}

\bigskip\noindent
{\bf Proof.}
Here as well we use the notation $L_i(W)$ for $i = 1, 2, ...$
just as in the proof of {\bf 7.11}.
It follows from {\bf 7.10(i)} that
$\mathfrak X$ has a filtration
\begin{align}
{\mathfrak X} \ \ = \ \ 
\begin{matrix}
\begin{matrix} 
    \boxed{ 
\begin{picture}(55,46)(0,0)
\put(20,40){$18a$}
\put(40,20){$18c$}
\put(0,20){$18b$}
\put(20,0){$18a$}
\put(10,18){\line(1,-1){10}}
\put(30,8){\line(1,1){10}}
\put(10,28){\line(1,1){10}}
\put(30,38){\line(1,-1){10}}
\end{picture}
%18a \\ 18b \ 18c \\ 18a
          }
            \end{matrix}
   \\
   \\
\end{matrix}
\
\Biggm|
\begin{matrix}
     \\
    9a \oplus 9b \oplus W 
     \\
\end{matrix}
\
{
\begin{matrix} 
  \\
  \\
\Biggm| \ 
\end{matrix}
}
\begin{matrix}
     \\  
     \\
 \begin{matrix} 
    \boxed{
\begin{picture}(55,46)(0,0)
\put(20,40){$18a$}
\put(40,20){$18c$}
\put(0,20){$18b$}
\put(20,0){$18a$}
\put(10,18){\line(1,-1){10}}
\put(30,8){\line(1,1){10}}
\put(10,28){\line(1,1){10}}
\put(30,38){\line(1,-1){10}}
\end{picture}
%18a \\ 18b \ 18c \\ 18a
          }
            \end{matrix}
\end{matrix}
{
\begin{matrix}
 \\ \\ \\ \\ \\ 
 \end{matrix}
}
\end{align}
namely, ${\mathfrak X}$ has submodules $Y$ and $Z$ such that
${\mathfrak X} \supsetneqq Y \supsetneqq Z$,
${\mathfrak X}/Y \cong Z \cong
\begin{matrix} 
\boxed{ 
\begin{picture}(55,46)(0,0)
\put(20,40){$18a$}
\put(40,20){$18c$}
\put(0,20){$18b$}
\put(20,0){$18a$}
\put(10,18){\line(1,-1){10}}
\put(30,8){\line(1,1){10}}
\put(10,28){\line(1,1){10}}
\put(30,38){\line(1,-1){10}}
\end{picture}
%18a \\ 18b \ 18c \\ 18a
          }
            \end{matrix}
$
and $Y/Z \cong 9a \oplus 9b \oplus W$.
On the other hand, {\bf 7.12} says that
\begin{equation}
{\mathfrak X} =
       \begin{matrix} 
 \boxed{
\begin{picture}(90,46)(0,0)
\put(60,40){$18c$}
\put(20,40){$18b$}
\put(80,20){$9a$}
\put(40,20){$18a$}
\put(0,20){$9b$}
\put(60,0){$18b$}
\put(20,0){$18c$}
\put(10,18){\line(1,-1){10}}
\put(50,18){\line(1,-1){10}}
\put(30,8){\line(1,1){10}}
\put(70,8){\line(1,1){10}}
\put(10,28){\line(1,1){10}}
\put(50,28){\line(1,1){10}}
\put(30,38){\line(1,-1){10}}
\put(70,38){\line(1,-1){10}}
\end{picture}
%         18b \ \ 18c  \\
%                      9b \ 18a \ 9a \\
%                       18c \ \ 18b
          }
       \end{matrix}
\ \bigoplus \ 
P(18a).
\end{equation} 
Then, we know by (18), (19) and {\bf 3.6(iv)} that
\begin{equation}
W = 9a + 9b + 9c + 9d + 2 \times 18b + 2 \times 18c,
\quad
{\text{as \ composition \ factors}}.
\end{equation}
By {\bf 7.11} and (20), we know $j(W) \geqslant 3$.

Assume that $j(W) \geqslant 4$. Then, $j(W) = 4$ by {\bf 7.8(iii)}.
Since $L_1(W) \cong 18b \oplus 18c$ by {\bf 7.11},
we get by {\bf 3.6(vi)} that
$$
    L_4(W) \ \Big| \ 
    L_4(P(18b)) \bigoplus L_4(P(18c)) \ = \ 
    (9a \oplus 18a \oplus 9d) \bigoplus (9b \oplus 18a \oplus 9c)
$$
and
$$
    L_4(W) \ \Big| \ {\mathrm{soc}}(W) \ = \ 18b \oplus 18c.
$$
This is a contradiction. 

Hence $j(W) = 3$.
Thus, again by {\bf 7.11}, (20) and {\bf 3.6(vi)}, we know that
$W$ has radical and socle series
\begin{equation}
  \boxed{  \begin{matrix}  18b \ 18c  \\
                           9b \ 9c \ 9a \ 9d  \\
                           18c \ 18b
           \end{matrix}
        }
\end{equation}
Now, as in the proof of {\bf 7.11}, we get by {\bf 7.8} that
$W = W_1 \oplus W_2$ where $W_i$ is
a non-simple non-projective indecomposable self-dual $B$-module
for $i = 1, 2$.
Then, by (21), we may assume that
$L_1(W_1) \cong 18b$, ${\mathrm{soc}}(W_1) \cong 18c$,
$L_1(W_2) \cong 18c$ and ${\mathrm{soc}}(W_2) \cong 18b$
since $(18b)^{\vee} \cong 18c$ by {\bf 3.6(iii)}.
Hence the structures of $P(18b)$ and $P(18c)$ in {\bf 3.6(vi)}
yield that
$$
W_1 \ = \ 
       \boxed{   \begin{matrix} 18b \\ 9b \ 9c \\ 18c
                 \end{matrix}
             }
\qquad {\text{and}} \qquad 
W_2 \ = \ 
       \boxed{   \begin{matrix} 18c \\ 9a \ 9d \\ 18b
                 \end{matrix}
             }
\qquad\blacksquare
$$

%\bigskip\noindent
%{\bf 7.14.Lemma.}
%{\it It holds that
%$$
%F(S_5) \ = \ 
%\boxed{ 
%  \begin{matrix}   18b  \\
%                 9b \ \ 9c \\
%                   18c
%  \end{matrix}
%      },
%\qquad
%F(S_7) \ = \ 
%\boxed{ 
%  \begin{matrix}   18c  \\
%                 9a \ \ 9d \\
%                   18b
%  \end{matrix}
%      },     
%$$
%}
%
%\bigskip\noindent
%{\bf Proof.}
%{\bf }
%\quad$\blacksquare$

\bigskip

%\newpage

\begin{flushleft}
{\bf 8. Proof of main results}
\end{flushleft}

\bigskip\noindent
{\bf 8.1.Notation.}
We still keep  the notation
$F$, $j$, $B'$, 
%$\tilde\sigma$, 
$f'$ and $g'$,
see {\bf 6.4}, {\bf 3.7} and {\bf 5.2}--{\bf 5.4}.
Set $E = SD_{16}$, and let $P \rtimes E$ be the
canonical semi-direct product such that
$E$ acts on $P$ faithfully.
Recall that
${\mathrm{Aut}}(P) \cong {\mathrm{GL}}_2(3)$
since $P = C_3 \times C_3$, and hence
$SD_{16}$ is a Sylow $2$-subgroup of 
${\mathrm{GL}}_2(3)$.

\bigskip\noindent
{\bf 8.2.Lemma.}
{\it
 The non-principal block algebra $A$ of $\sf HN$ and
the principal block algebra $A'$ of $\sf HS$ are 
Puig equivalent.
}

\bigskip\noindent
{\bf Proof.}
Let $j$ be the same as in {\bf 3.6(ii)}.
Since $jBj \cong {\mathcal O}[P \rtimes E] = B'$
as interior $P$-algebras by {\bf 3.6(ii)}, 
we can identify $jBj$ and $B'$.
Define a functor 
$F': {\mathrm{mod}}{\text{-}}B \rightarrow 
     {\mathrm{mod}}{\text{-}}B'$
via $F'(-) = - \otimes_B Bj$.
By {\bf 3.6(ii)}, $F'$ induces a Puig equivalence
(which is stronger than a Morita equivalence) between $B$ and $B'$.
In the following we use the information on the structures
of PIMs in $B$ and $B'$ described in
{\bf 3.6(vi)} and {\bf 5.2(iii)}, respectively, 
without quoting these statements.

Then, first of all, we know that
$F'(18a) = 2a$ by looking at the PIMs
$P(18a)$ and $P(2a)$. Similarly, we know at least that
$\{ F'(9a), F'(9b), F'(9c), F'(9d) \}
 = \{ 1a = k_{H'}, 1b, 1c, 1d \}$.
It follows from {\bf 5.4} that 
$1x \otimes 1x = 1a$ for any $x \in \{ a, b, c, d\}$
since they are just in ${\mathrm{Irr}}(E)$.
Hence a technique of self-Puig equivalence
in \cite[2.8.Lemma]{KoshitaniKunugiWaki2008} can be
used just as in the proof of
\cite[6.8.Lemma]{KoshitaniKunugiWaki2008}.
Namely, we can assume that $F'(9a) = 1a$. 
Hence, by comparing the second Loewy layers
of $P(9a)$ and $P(1a)$, we get $F'(18b) = 2b$.
Similarly, by looking at the third Loewy layers of
$P(9a)$ and $P(1a)$, we have $F'(9b) = 1b$.
If we look at the fourth Loewy layers of these PIMs,
then we know $F'(18c) = 2c$.
Thus, by looking at the second Loewy layers of
$P(18c)$ and $P(2c)$, we know also that $F'(9d) = 1d$.
These mean that $F'(9c) = 1c$.
Namely, we can assume that

\begin{gather}
  \begin{aligned}
F'(9a) = 1a, \ F'(9b) = 1b, \ &F'(9c) = 1c, \ F'(9d) = 1d,
\\
F'(18a) = 2a, \ F'(18b) &= 2b, \ F'(18c) = 2c. 
  \end{aligned}
\end{gather}

We know by {\bf 7.13} that Case(a) or Case(b) happens.

\bigskip\noindent
Assume, first, that {\bf Case(b)} occurs.
Then, by bunching up 
{\bf 2.2}, {\bf 7.1}, {\bf 7.5}, {\bf 7.6}, {\bf 7.13} and {\bf 5.7}, 
we get the %following 
diagram shown in Table \ref{caseb}.

%
%\smallskip\noindent
%
%JM: introduced table, otherwise pagebreak is absolutely nasty
\begin{table}\caption{Case(b).}\label{caseb}
{
\small
\begin{center}
\begin{tabular}{c c c c c c c}
mod-$A$ & $\overset{F}{\longrightarrow}$    
        & mod-$B$ 
%JM: replaced tensor functor by F', deleted identifications
        & $\overset{F'}{\longrightarrow}$ 
%$\overset{-\otimes_B Bj}{\longrightarrow}$ 
        & $\mathrm{mod}{\text{-}}B'$ 
%{$  \begin{matrix} {\mathrm{mod}}{\text{-}}jBj
%          \\ \cong {\mathrm{mod}}{\text{-}}B'
%          \\    =  {\mathrm{mod}}{\text{-}}{\mathcal O[P \rtimes E]}
%                   \end{matrix}$ }
        & $\overset{f'^{-1}}{\longrightarrow}$  
        & mod-$A'$  \\
\hline
\\
$S_1$ & $\mapsto$ & $\boxed{9a}$ & $\mapsto$ & $\boxed{1a}$ 
      & $\mapsto$ & $k_{G'}$          
\\  \\
$S_2$ & $\mapsto$ & $\boxed{9b}$ & $\mapsto$ & $\boxed{1b}$ 
      & $\mapsto$ & $154$          \\  \\
$S_3$ & $\mapsto$ & $\boxed{9c}$ & $\mapsto$ & $\boxed{1c}$ 
      & $\mapsto$ & $22$          \\  \\
$S_4$ & $\mapsto$ & 
            $
 \begin{matrix} 
\boxed{
\begin{picture}(55,46)(0,0)
\put(20,40){$18a$}
\put(40,20){$18c$}
\put(0,20){$18b$}
\put(20,0){$18a$}
\put(10,18){\line(1,-1){10}}
\put(30,8){\line(1,1){10}}
\put(10,28){\line(1,1){10}}
\put(30,38){\line(1,-1){10}}
\end{picture}
%18a \\ 18b \ 18c \\ 18a
                   }
                     \end{matrix}
$ 
      & $\mapsto$ & 
            $
 \begin{matrix} 
\boxed{
\begin{picture}(50,46)(0,0)
\put(20,40){$2a$}
\put(40,20){$2c$}
\put(0,20){$2b$}
\put(20,0){$2a$}
\put(10,18){\line(1,-1){10}}
\put(30,8){\line(1,1){10}}
\put(10,28){\line(1,1){10}}
\put(30,38){\line(1,-1){10}}
\end{picture}
%2a \\ 2b \ 2c \\ 2a
                   }
                     \end{matrix}
$  
      & $\mapsto$ & $1253$          \\  \\
$S_5$ & $\mapsto$ & 
            $
 \begin{matrix} 
\boxed{
\begin{picture}(50,46)(0,0)
\put(20,40){$18c$}
\put(40,20){$9d$}
\put(0,20){$9a$}
\put(20,0){$18b$}
\put(10,18){\line(1,-1){10}}
\put(30,8){\line(1,1){10}}
\put(10,28){\line(1,1){10}}
\put(30,38){\line(1,-1){10}}
\end{picture}
%18c \\ 9a \ 9d \\ 18b
                   }
                     \end{matrix}
$ 
      & $\mapsto$ & 
            $
 \begin{matrix} 
\boxed{
\begin{picture}(50,46)(0,0)
\put(20,40){$2c$}
\put(40,20){$1d$}
\put(0,20){$1a$}
\put(20,0){$2b$}
\put(10,18){\line(1,-1){10}}
\put(30,8){\line(1,1){10}}
\put(10,28){\line(1,1){10}}
\put(30,38){\line(1,-1){10}}
\end{picture}
%2c \\ 1a \ 1d \\ 2b
                   }
                     \end{matrix}
$  
      & $\mapsto$ & $321$          \\ \\
$S_6$ & $\mapsto$ 
      &  
                $  \begin{matrix}  
  \boxed{
\begin{picture}(90,46)(0,0)
\put(40,40){$18a$}
\put(20,20){$18c$}
\put(40,0){$18a$}
\put(60,20){$18b$}
\put(80,40){$9d$}
\put(0,0){$9d$}
\put(10,8){\line(1,1){10}}
\put(30,18){\line(1,-1){10}}
\put(50,8){\line(1,1){10}}
\put(30,28){\line(1,1){10}}
\put(50,38){\line(1,-1){10}}
\put(70,28){\line(1,1){10}}
\end{picture}
%&     & 18a &     & 9d \\
%                                  & 18c &     & 18b &    \\
%                               9d &     & 18a &     &
        }
                  \end{matrix} $
& $\mapsto$
 & 
    $              \begin{matrix} 
   \boxed{
\begin{picture}(90,46)(0,0)
\put(40,40){$2a$}
\put(20,20){$2c$}
\put(40,0){$2a$}
\put(60,20){$2b$}
\put(80,40){$1d$}
\put(0,0){$1d$}
\put(10,8){\line(1,1){10}}
\put(30,18){\line(1,-1){10}}
\put(50,8){\line(1,1){10}}
\put(30,28){\line(1,1){10}}
\put(50,38){\line(1,-1){10}}
\put(70,28){\line(1,1){10}}
\end{picture}
% &     & 2a &     & 1d \\
%                                  & 2c &     & 2b &    \\
%                               1d &     & 2a &     &
            }
                  \end{matrix} $
      & $\mapsto$ & $748$          \\ \\
$S_7$ & $\mapsto$ & 
            $
\begin{matrix}
\boxed{ 
\begin{picture}(50,46)(0,0)
\put(20,40){$18b$}
\put(40,20){$9c$}
\put(0,20){$9b$}
\put(20,0){$18c$}
\put(10,18){\line(1,-1){10}}
\put(30,8){\line(1,1){10}}
\put(10,28){\line(1,1){10}}
\put(30,38){\line(1,-1){10}}
\end{picture}
% 18b \\ 9b \ 9c \\18c
                   }
                     \end{matrix}
$ 
      & $\mapsto$ & 
            $
 \begin{matrix} 
\boxed{
\begin{picture}(50,46)(0,0)
\put(20,40){$2b$}
\put(40,20){$1c$}
\put(0,20){$1b$}
\put(20,0){$2c$}
\put(10,18){\line(1,-1){10}}
\put(30,8){\line(1,1){10}}
\put(10,28){\line(1,1){10}}
\put(30,38){\line(1,-1){10}}
\end{picture}
%2b \\ 1b \ 1c \\ 2c
                   }
                     \end{matrix}
$  
      & $\mapsto$ & $1176$ 
\end{tabular}
\end{center}
}
\hrulefill
\end{table}

%\bigskip\bigskip\bigskip
%\bigskip\noindent
First, all the three functors above
are given by bimodules which are $p$-permutation
modules over $\mathcal O[G_1 \times H_1]$ 
for corresponding two finite
groups $G_1$ and $H_1$, which are $\Delta P$-projective,
and also which induce a stable equivalence of Morita type
at each step,
if we indentify the source algebra $jBj$ as
$\mathcal O[P \rtimes E]$.

%\medskip
%\noindent
Secondly, it has to be noted that
all non-simple modules in the above diagram are
uniquely determined (up to isomorphism)
by just the diagrams given in the above boxes:
%JM: Added the argument from an earlier e-mail
This is clear for $F(S_1)$, $F(S_2)$, $F(S_1)$,
$f'(k_{G'})$, $f'(154)$, and $f'(22)$ anyway,
as well as for $F(S_4)$ and $f'(1253)$ by the structure
of $P(18a)$ and $P(2a)$ given in {\bf 3.6(vi)} and {\bf 5.2(iii)}.

To tackle $F(S_6)$, the structure of $P(18a)$ specified
in {\bf 3.6(vii)} shows that $P(18a)$ has a unique quotient
with composition factors $9d+ 2\times 18a+ 18b+ 18c$.
Moreover, $P(9d)$ has a unique quotient
with composition factors $9d+18a+18b$.
Since they both have a unique submodule with composition factors
$18a+18b$, the glueing to yield $F(S_6)$ also is
uniquely defined, and thus $F(S_6)$ is uniquely determined by
the diagram given. For $f'(748)$ we argue similarly using {\bf 5.2(iv)}.

We consider $F(S_7)$: 
Note first that for $P(18b)$ there is no Alperin diagram defined.
By {\bf 3.6(vi)}, let $X$ be the unique quotient module of $P(18b)$ 
having radical and socle series
$\boxed{\begin{matrix} 18b \\ 9b \ 9c 
\end{matrix}}$.
By the structure of $P(18b)$ given in {\bf 3.6(vi)} we have 
$[\Omega(X),18a]^B=1$, hence using {\bf 3.6(vii)} there is a homomorphism
$\varphi\in\mathrm{Hom}_B(P(18a),\Omega(X))$ such that
$$ \mathrm{Im}(\varphi)  \ = \
\begin{matrix} \boxed{
\begin{picture}(75,66)(0,0)
\put(40,60){$18a$}
\put(20,40){$18c$}
\put(60,40){$18b$}
\put(0,20){$9a$}
\put(20,20){$9d$}
\put(40,20){$18a$}
\put(20,0){$18b$}
\put(10,18){\line(1,-1){10}}
\put(30,8){\line(1,1){10}}
\put(25,8){\line(0,1){10}}
\put(10,28){\line(1,1){10}}
\put(25,28){\line(0,1){10}}
\put(30,38){\line(1,-1){10}}
\put(50,28){\line(1,1){10}}
\put(30,48){\line(1,1){10}}
\put(50,58){\line(1,-1){10}}
\end{picture} } \end{matrix} $$
This implies $\Omega(X)/\mathrm{Im}(\varphi)\cong 18c$.
Since $18c$ occurs exactly twice as a composition factor of $\Omega(X)$,
and also is a composition factor of $\mathrm{Im}(\varphi)$,
we conclude that $[\Omega(X),18c]^B=1$, thus
$\mathrm{dim}_k[\mathrm{Ext}^1_B(X,18c)]=1$.
Hence a module having radical and socle series
$\boxed{\begin{matrix} 18b \\ 9b \ 9c \\ 18c
\end{matrix}}$
is uniquely defined.
For $F(S_5)$, $f'(1176)$, and $f'(321)$ we argue similarly.

%\medskip
%\noindent
Then, it follows from {\bf 2.15} that 
$A$ and $A'$ are splendidly stable equivalent of Morita type,
that is, $A$ and $A'$ are stable equivalent which is realized
by an $\mathcal O[G \times G']$-bimodule which is a $p$-permutation
module and $\Delta P$-projective.
Hence, first of all, the stable equivalence actually
gives a Morita equivalence by a result of Linckelmann
\cite[Theorem 2.1(iii)]{Linckelmann1996MathZ}.
Then, if we look at the proof of 
\cite[Theorem 2.1(iii)]{Linckelmann1996MathZ}
which is actually given in
\cite[Remark 2.7]{Linckelmann1996MathZ}, 
we know that the Morita equivalence between $A$ and $A'$
gives a bijection such as
%$S_1 \longleftrightarrow k_{G'} \qquad {\text{and}} \qquad$
$S_5 \leftrightarrow 321$.
Hence, we must have equalities between the corresponding
Cartan invariants, namely,
$c(S_5, S_5) = c(321, 321)$.
However, we get 
that $c(S_5, S_5) = 3$ by {\bf 4.1}, and on the other hand,
that $c(321, 321) = 2$ by {\bf 5.8}.
This is a contradiction.
Thus, {\bf Case(b)} cannot happen.

%JM: introduced table
\begin{table}\caption{Case(a).}\label{casea}
{\small
\begin{center}
\begin{tabular}{c c c c c c c}
mod-$A$ & $\overset{F}{\longrightarrow}$    
        & mod-$B$ 
%JM: replaced tensor functor by F', deleted identifications
        & $\overset{F'}{\longrightarrow}$ 
         %$\overset{-\otimes_B Bj}{\longrightarrow}$ 
        & $\mathrm{mod}{\text{-}}B'$ 
%{$  \begin{matrix} {\mathrm{mod}}{\text{-}}jBj
%          \\ \cong {\mathrm{mod}}{\text{-}}B'
%          \\    =  {\mathrm{mod}}{\text{-}}{\mathcal O[P \rtimes E]}
%                   \end{matrix}$ }
        & $\overset{f'^{-1}}{\longrightarrow}$  
        & mod-$A'$  \\
\hline
\\
$S_1$ & $\mapsto$ & $\boxed{9a}$ & $\mapsto$ & $\boxed{1a}$ 
      & $\mapsto$ & $k_{G'}$          
\\  \\
$S_2$ & $\mapsto$ & $\boxed{9b}$ & $\mapsto$ & $\boxed{1b}$ 
      & $\mapsto$ & $154$          \\  \\
$S_3$ & $\mapsto$ & $\boxed{9c}$ & $\mapsto$ & $\boxed{1c}$ 
      & $\mapsto$ & $22$          \\  \\
$S_4$ & $\mapsto$ & 
            $
\begin{matrix} 
\boxed{ 
\begin{picture}(55,46)(0,0)
\put(20,40){$18a$}
\put(40,20){$18c$}
\put(0,20){$18b$}
\put(20,0){$18a$}
\put(10,18){\line(1,-1){10}}
\put(30,8){\line(1,1){10}}
\put(10,28){\line(1,1){10}}
\put(30,38){\line(1,-1){10}}
\end{picture}
%18a \\ 18b \ 18c \\ 18a
                   }
                     \end{matrix}
$ 
      & $\mapsto$ & 
            $
\begin{matrix} 
\boxed{ 
\begin{picture}(50,46)(0,0)
\put(20,40){$2a$}
\put(40,20){$2c$}
\put(0,20){$2b$}
\put(20,0){$2a$}
\put(10,18){\line(1,-1){10}}
\put(30,8){\line(1,1){10}}
\put(10,28){\line(1,1){10}}
\put(30,38){\line(1,-1){10}}
\end{picture}
%2a \\ 2b \ 2c \\ 2a
                   }
                     \end{matrix}
$  
      & $\mapsto$ & $1253$          \\  \\
$S_5$ & $\mapsto$ & 
            $
\begin{matrix} 
\boxed{ 
\begin{picture}(50,46)(0,0)
\put(20,40){$18b$}
\put(40,20){$9c$}
\put(0,20){$9b$}
\put(20,0){$18c$}
\put(10,18){\line(1,-1){10}}
\put(30,8){\line(1,1){10}}
\put(10,28){\line(1,1){10}}
\put(30,38){\line(1,-1){10}}
\end{picture}
%18b \\ 9b \ 9c \\ 18c
                   }
                     \end{matrix}
$ 
      & $\mapsto$ & 
            $
 \begin{matrix} 
\boxed{
\begin{picture}(50,46)(0,0)
\put(20,40){$2b$}
\put(40,20){$1c$}
\put(0,20){$1b$}
\put(20,0){$2c$}
\put(10,18){\line(1,-1){10}}
\put(30,8){\line(1,1){10}}
\put(10,28){\line(1,1){10}}
\put(30,38){\line(1,-1){10}}
\end{picture}
%2b \\ 1b \ 1c \\ 2c
                   }
                     \end{matrix}
$  
      & $\mapsto$ & $1176$          \\ \\
$S_6$ & $\mapsto$ 
      &   
 $                 \begin{matrix}  
 \boxed{
\begin{picture}(90,46)(0,0)
\put(40,40){$18a$}
\put(20,20){$18c$}
\put(40,0){$18a$}
\put(60,20){$18b$}
\put(80,40){$9d$}
\put(0,0){$9d$}
\put(10,8){\line(1,1){10}}
\put(30,18){\line(1,-1){10}}
\put(50,8){\line(1,1){10}}
\put(30,28){\line(1,1){10}}
\put(50,38){\line(1,-1){10}}
\put(70,28){\line(1,1){10}}
\end{picture}
%&     & 18a &     & 9d \\
%                                  & 18c &     & 18b &    \\
%                               9d &     & 18a &     &
        }
                  \end{matrix} $
& $\mapsto$
 & 
   $               \begin{matrix} 
   \boxed{
\begin{picture}(90,46)(0,0)
\put(40,40){$2a$}
\put(20,20){$2c$}
\put(40,0){$2a$}
\put(60,20){$2b$}
\put(80,40){$1d$}
\put(0,0){$1d$}
\put(10,8){\line(1,1){10}}
\put(30,18){\line(1,-1){10}}
\put(50,8){\line(1,1){10}}
\put(30,28){\line(1,1){10}}
\put(50,38){\line(1,-1){10}}
\put(70,28){\line(1,1){10}}
\end{picture}
% &     & 2a &     & 1d \\
%                                  & 2c &     & 2b &    \\
%                               1d &     & 2a &     &
            }
                  \end{matrix} $
      & $\mapsto$ & $748$          \\ \\
$S_7$ & $\mapsto$ & 
            $
\begin{matrix} 
\boxed{ 
\begin{picture}(50,46)(0,0)
\put(20,40){$18c$}
\put(40,20){$9d$}
\put(0,20){$9a$}
\put(20,0){$18b$}
\put(10,18){\line(1,-1){10}}
\put(30,8){\line(1,1){10}}
\put(10,28){\line(1,1){10}}
\put(30,38){\line(1,-1){10}}
\end{picture}
%18c \\ 9a \ 9d \\18b
                  }
                   \end{matrix}
$ 
      & $\mapsto$ & 
            $
\begin{matrix} 
\boxed{ 
\begin{picture}(50,46)(0,0)
\put(20,40){$2c$}
\put(40,20){$1d$}
\put(0,20){$1a$}
\put(20,0){$2b$}
\put(10,18){\line(1,-1){10}}
\put(30,8){\line(1,1){10}}
\put(10,28){\line(1,1){10}}
\put(30,38){\line(1,-1){10}}
\end{picture}
%2c \\ 1a \ 1d \\ 2b
                   }
                     \end{matrix}
$  
      & $\mapsto$ & $321$ 
\end{tabular}
\end{center}
}
\hrulefill
\end{table}

\bigskip\noindent
This means that only {\bf Case(a)} occurs, as
is shown in Table \ref{casea}.
Then, again the same argument given above still works.
Namely, we have a Morita equivalence between $A$ and $A'$,
and hence the Morita equivalence is 
a Puig equivalence by a result of Puig
(and, independently, of Scott)
\cite[Remark 7.5]{Puig1999}, see
\cite[Theorem 4.1]{Linckelmann2001}.
\quad$\blacksquare$

\bigskip\noindent
{\bf 8.3.Proofs of 1.3 and 1.4}.
Recall that a Puig equivalence lifts from $k$ to $\mathcal O$
by a result of Puig \cite[7.8.Lemma]{Puig1988Inv}
(see \cite[(38.8)Proposition]{Thevenaz}, and that
so does a splendid Rickard equivalence by a result of
Rickard \cite[Theorem 5.2]{Rickard1996},
see \cite[P.75, lines $-17 \sim -16$]{Harris}.
Thus, it is enough to consider blocks $A$, $B$,
$A'$ and $B'$ only over $k$. Thus, we get {\bf 1.4}
by {\bf 8.2}.

By results of Okuyama
\cite[Example 4.8]{Okuyama1997} and
\cite[Corollary 2]{Okuyama2000}, the conjectures
{\bf 1.1} and {\bf 1.2} hold for $A'$. Namely, we get the
following diagram:

$$
\begin{CD}
     A   @> {\text{\rm{\small{Puig \ equiv.}}}} >>  A'     
\\
     @.    @VV {\text{\rm{\small{splendid \ Rickard \ equiv.}}}}V
\\
     B   @<< {\text{\rm{\small{Puig \ equiv.}}}} <    B'
\end{CD}
$$

\medskip\noindent
Therefore, we finally get that $A$ and $B$ are
splendidly Rickard equivalent. That is, the proof of
{\bf 1.3} is completed.
\quad$\blacksquare$

\bigskip\noindent
{\bf 8.4.Proof of 1.5}.
We get {\bf 1.5} from {\bf 3.2} and {\bf 1.3}.
\quad$\blacksquare$

%\bigskip\bigskip
%\newpage
%%
%
%
%\begin{flushleft}
%{\bf 9. Appendix}
%\end{flushleft}%
%
%By using almost the same method in {\bf 8.2},
%we can answer a question mentioned in our
%previous paper
%\cite[6.14.Remark and question]{KoshitaniKunugiWaki2008}.
%Namely, in the statements of
%\cite[Lemma 6.12]{KoshitaniKunugiWaki2008},
%only the first {\bf Case 1} happens. Therefore, we can
%improve our previous theorem
%\cite[Theorem 1.4]{KoshitaniKunugiWaki2008}, namely,
%the non-principal block algebra of $\mathcal OJ_4$
%with defect group $C_3 \times C_3$ and  the
%principal block algebra of $\mathcal O\mathcal A_8$
%are actually Puig equivalent, where $J_4$ is the
%Janko simple group and $\mathcal A_8$ is the alternating
%group on 8 letters.
%
%{\LARGE\bf
%(Give much more detailed proof!!)
%} 

\bigskip

%SK
\newpage

\begin{center}{\bf Acknowledgements}
\end{center}
\smallskip\noindent
The first author thanks Gabriel Navarro 
for informing \cite[(3.2)Lemma]{Navarro2004}.

A part of this work was done while the first author was
staying in Braunschweig Technical University and 
RWTH Aachen University in 2007 and 2009. He is grateful to
Bettina Eick and Gerhard Hiss for their kind hospitality.
For this research the first author was partially
supported by the Japan Society for Promotion of Science (JSPS),
Grant-in-Aid for Scientific Research 
(C)17540010, 2005--2007, and
(C)20540008, 2008--2010.

\end{document}